\pdfoutput=1
\documentclass[11pt]{amsart}
\usepackage[margin=2.5cm]{geometry}
\usepackage{amsmath,amssymb,amsthm,mathtools,booktabs,enumitem,array,longtable}
\usepackage{graphicx}
\usepackage{tikz}
\usetikzlibrary{arrows.meta,positioning,decorations.markings,patterns,calc}
\usepackage{placeins}
\usepackage[hyphens]{url}
\usepackage[colorlinks=true,linkcolor=blue,citecolor=blue,urlcolor=blue]{hyperref}

\newtheorem{theorem}{Theorem}[section]
\newtheorem{proposition}[theorem]{Proposition}
\newtheorem{lemma}[theorem]{Lemma}
\newtheorem{corollary}[theorem]{Corollary}
\newtheorem{inthm}{Theorem}

\theoremstyle{definition}
\newtheorem{definition}[theorem]{Definition}
\theoremstyle{remark}
\newtheorem{remark}[theorem]{Remark}
\newtheorem*{remark*}{Remark}

\newcommand{\E}{\mathbb{E}}
\newcommand{\Pp}{\mathbb{P}}
\newcommand{\R}{\mathbb{R}}
\newcommand{\Z}{\mathbb{Z}}
\newcommand{\Var}{\operatorname{Var}}
\newcommand{\Cov}{\operatorname{Cov}}
\newcommand{\KL}{\operatorname{KL}}
\newcommand{\Insp}{\operatorname{Insp}}
\newcommand{\dist}{\operatorname{dist}}
\newcommand{\one}{\mathbf 1}
\newcommand{\cV}{\mathcal V}
\newcommand{\cF}{\mathcal F}
\newcommand{\cG}{\mathcal G}
\newcommand{\cH}{\mathcal H}
\newcommand{\cD}{\mathcal D}
\newcommand{\cX}{\mathcal X}
\newcommand{\SLE}{\mathrm{SLE}}
\newcommand{\Dob}{\mathrm{Dob}}
\newcommand{\Dab}{D_{x_1x_2}}
\newcommand{\pp}{\mathfrak p}
\newcommand{\ppc}{\mathfrak p_c}
\newcommand{\Osc}{\operatorname{Osc}}
\newcommand{\piD}{\pi^{\rm D}_4}
\newcommand{\AD}{A^{\rm D}_4}
\newcommand{\rD}{r^{\rm D}_4}
\newcommand{\epsL}{\varepsilon_L}
\newcommand{\Kmix}{K_{\rm mix}}
\newcommand{\Rr}{\mathsf R}
\newcommand{\Dd}{\mathsf D}
\newcommand{\Ww}{\mathsf W}

\DeclareRobustCommand{\tagref}[2]{\hyperlink{tag:#1}{\textup{#2}}}
\DeclareRobustCommand{\tagdef}[2]{\hypertarget{tag:#1}{#2}}
\begin{document}
\sloppy

\title[Near-critical Ising interfaces in the thermal direction: absolute continuity]{Near-critical Ising interfaces in the thermal direction: Absolute continuity with respect to SLE$_3$}
\author[Yu Feng]{Yu Feng$^{\dagger}$}
\thanks{$^{\dagger}$yufeng\_proba@163.com, yuproba@umich.edu}
\makeatletter
\g@addto@macro\@setauthors{\par\vspace{4pt}{\centering\normalfont\small University of Michigan, Ann Arbor, USA\par}}
\makeatother

\begin{abstract}
The scaling limits of interfaces in the critical two-dimensional Ising model are described by the
Schramm--Loewner evolution $\SLE_3$. While, in the scaling limit, the law of the whole near-critical
FK-Ising configuration is known to be singular with respect to the critical law in the thermal
direction, and the same is expected for the spin configuration, we prove in this paper that the law
of a single interface is not. More precisely, we consider the Ising model with mesh $\delta$ in the
Dobrushin square at $\beta=\beta_c-m\delta$, where $m\ne0$ is fixed, and we prove that the
Radon--Nikodym derivative of the law of the interface with respect to its critical law has moments
of every positive and negative order, bounded uniformly in $\delta$, so that the law of every
subsequential scaling limit of the interface is mutually absolutely continuous with respect to that
of $\SLE_3$. Our analysis explains the mechanism behind this result: the logarithmically growing
specific heat (that is, the variance of the Hamiltonian, normalized by the number of edges) cancels
between the full domain and the domain slit along the interface.
\end{abstract}


\date{\today}
\renewcommand*{\shortauthors}{}
\maketitle

\setcounter{tocdepth}{1}
\tableofcontents
\enlargethispage{2pt}

\section{Introduction}\label{sec:intro}
\subsection{Background and motivation}\label{ss:background} The Ising model was introduced by Lenz
and studied by Ising \cite{Isi25} as a model of ferromagnetism. Spins $\pm1$ are placed at the
vertices of a lattice, and neighboring spins interact so as to favor alignment, with a strength
measured by the inverse temperature $\beta$. On the square lattice, the model undergoes a continuous
phase transition at the critical inverse temperature $\beta_c=\frac12\log(1+\sqrt2)$
\cite{KW41,Ons44,Yan52,BDC12a}, which corresponds to the self-dual point of the associated
random-cluster model. At high temperature, $\beta<\beta_c$, spin correlations decay exponentially,
with a rate called the inverse correlation length; at low temperature, $\beta>\beta_c$, the model is
ordered. As $\beta$ approaches $\beta_c$, the correlation length tends to infinity, and the
large-scale behavior of the model is expected to become scale invariant and, in fact, conformally
invariant, in the sense of conformal field theory \cite{BPZ84}.

Following the introduction by Schramm \cite{Sch00} of the Schramm--Loewner evolution $\SLE_\kappa$,
a one-parameter family of random curves characterized by conformal invariance and a domain Markov
property, it became possible to describe the scaling limits of interfaces of critical lattice models
in a mathematically precise way. For the Ising model, conformal invariance at $\beta=\beta_c$ was
established through Smirnov's fermionic observables \cite{Smi10,CS12}. Building on these and on the
precompactness criteria of Kemppainen and Smirnov \cite{KS17}, Chelkak, Duminil-Copin, Hongler,
Kemppainen and Smirnov \cite{CDHKS14} proved that the interface of the Ising model with Dobrushin
boundary conditions converges to $\SLE_3$, and that the interface of the associated random-cluster
(FK-Ising) model converges to $\SLE_{16/3}$. These convergences were then extended to more general
boundary conditions, domains and lattices~\cite{HK13,Izy15,Izy17,CHI21,Izy22,Che24,FPW24,Par25}.

Away from the critical point, the Ising correlation length is finite, and it is comparable to
$|\beta-\beta_c|^{-1}$ as $\beta\to\beta_c$ (\cite[Thm.~2]{BDC12} and \cite[Thms.~1.3
and~1.13]{DCM}). In the \emph{near-critical} regime, the inverse temperature approaches $\beta_c$ as
the mesh $\delta$ tends to zero, at the rate $\beta=\beta_c-m\delta$, $m\ne0$, for which the
correlation length is comparable to the size $n=1/\delta$ of the domain; in the scaling limit, one
then expects a nontrivial ``massive'' field theory, with mass parameter $m$. In the
Fortuin--Kasteleyn (random-cluster) representation \cite{FK72}, the edge parameter at inverse
temperature $\beta$ is given by $\pp(\beta):=1-e^{-2\beta}$, and its critical value is
$\ppc:=\pp(\beta_c)=\sqrt2/(1+\sqrt2)$. Duminil-Copin, Garban and Pete \cite{DGP14} studied the
near-critical FK-Ising model. They determined the correlation length \cite[Thm.~1.2]{DGP14}, defined
through crossing probabilities, and proved that these remain bounded away from $0$ and $1$ for
rectangles below this length, uniformly in the boundary conditions \cite[Thm.~1.3]{DGP14}. They also
highlighted that, in contrast with Bernoulli percolation, the correlation length is not governed by
the number of pivotal edges at criticality. Their analysis relies on Smirnov's fermionic observable
away from criticality, studied by Beffara and Duminil-Copin \cite{BDC12}. Park \cite{Par18,Par21}
and Chelkak, Izyurov and Mahfouf \cite{CIM21} extended the discrete complex analysis of Smirnov and
Chelkak--Smirnov \cite{Smi10,CS12} to the massive setting: the fermionic observables at
$\beta-\beta_c$ of order $\delta$ are no longer s-holomorphic, but they satisfy a massive version of
s-holomorphicity (cf.~\cite{BDC12}). Rigorous results on the scaling limits of near-critical Ising
correlations in the thermal direction go back to Palmer and Tracy~\cite{PT81} (see
also~\cite{Pal07}), and have seen significant progress in recent
years~\cite{Par18,Par22b,CIM21,PVW25,Mah25}. In the external field direction, the near-critical
scaling limit of the magnetization field was constructed by Camia, Garban and Newman \cite{CGN16},
and Camia, Jiang and Newman \cite{CJN20} proved its exponential decay.

For the scaling limits of near-critical interfaces, Makarov and Smirnov \cite{MS10} proposed a
description in terms of ``massive'' versions of $\SLE$. In particular, they stated, with a sketch of
proof, that massive loop-erased random walk converges to a massive $\SLE_2$ whose law is absolutely
continuous with respect to that of $\SLE_2$ \cite[Thm.~2.1]{MS10} (a complete proof was given by
Chelkak and Wan \cite{CW21}). They also asked whether, for $\kappa\le4$, the laws of massive and
usual $\SLE_\kappa$ are mutually absolutely continuous for all perturbations
\cite[Question~4.5]{MS10}, and whether they are singular for $\kappa>4$ \cite[Question~4.6]{MS10}.
For $\kappa=4$, Papon proved that the massive harmonic explorer converges to a massive $\SLE_4$
whose law is absolutely continuous with respect to that of $\SLE_4$ \cite{Pap23a}, and related
massive $\SLE_4$ and massive $\mathrm{CLE}_4$ to the massive planar Gaussian free field
\cite{Pap23b}. For the Ising interface in a magnetic field of order $\delta^{15/8}$, Papon
\cite{Pap24} proved convergence to a massive $\SLE_3$ whose law is absolutely continuous with
respect to that of $\SLE_3$, with an explicit density. Huang, Papon and Xia \cite{HPX26} proved a
quenched analogue for the random-field Ising model, and showed that, in contrast, the laws of the
subsequential limits of the outermost spin loops are almost surely singular with respect to the law
of $\mathrm{CLE}_3$ \cite[Thm.~1.2]{HPX26}. Related results exist for near-critical dimers
\cite{BHS22}. In the external field direction, even the law of the whole magnetization field in a
bounded domain is absolutely continuous with respect to the critical law (see the proof of
\cite[Prop.~1.5]{CGN16}).

The situation is fundamentally different for percolation. Nolin and Werner \cite{NW09} proved that
the laws of the subsequential scaling limits of near-critical percolation interfaces are singular
with respect to the law of $\SLE_6$, and Garban, Pete and Schramm \cite{GPS18} constructed the
near-critical full scaling limit of percolation. In this sense, a near-critical perturbation of
percolation is already visible on a single interface.

The thermal direction of the Ising model lies, in some sense, between these two cases. The
logarithmic growth of the \emph{specific heat} (that is, the variance of the Hamiltonian, normalized
by the number of edges) makes the law of the whole FK-Ising configuration singular: for every
$m\ne0$, the law of every subsequential limit of the near-critical FK-Ising configurations at
$\beta=\beta_c-m\delta$, in the quad-crossing topology, is singular with respect to the law of the
critical limit \cite[Thm.~2.5 and Remark~4]{GK25}. However, the perturbation is detected there
through a function of the whole configuration, which is not measurable with respect to a single
interface. As Garban and Kupiainen note, their result suggests that the law of the whole
near-critical spin configuration, described in the scaling limit by a near-critical
$\mathrm{CLE}_3$, is also singular, but this is far from an easy consequence of it
\cite[\S2.4]{GK25}. They also point out that the law of near-critical $\SLE_3$ is expected to be
absolutely continuous with respect to that of $\SLE_3$ \cite[\S2.4]{GK25}, and they conjecture the
analogous statement for the FK-Ising interface and $\SLE_{16/3}$ \cite[Conj.~2]{GK25}. This leads to
the following question. Is the law of a single near-critical interface also singular with respect to
that of its critical counterpart, as for percolation, or is it absolutely continuous, as in the
external field direction?

In this paper, we answer this question for the spin interface in the Dobrushin square (defined in
\S\ref{ss:model}), on both sides of the critical point. We prove that, for every fixed
$m\in\R\setminus\{0\}$, the law of every subsequential scaling limit of the interface at
$\beta=\beta_c-m\delta$ is mutually absolutely continuous with respect to the law of $\SLE_3$, with
densities in every $L^p$ in both directions (Theorem~\ref{thm:main}). This confirms, for the
Dobrushin square and every subsequential limit, the expectation expressed in \cite[\S2.4]{GK25} and
\cite[Question~1(2)]{CGS26} (see also \cite[Question~4.5]{MS10}) for the Ising interface in the
thermal direction, without identifying the limit with a massive $\SLE_3$ in the sense of
\cite{MS10}. The uniqueness of the subsequential limit and other domains will be discussed in
\S\ref{ss:outlook}. In \S\ref{ss:strategy}, we will also explain the mechanism that leads to
absolute continuity: the logarithmically growing specific heat, which is responsible for the
singularity of the law of the whole FK-Ising configuration, cancels \emph{edge by edge} between the
full domain and the domain slit along the interface, and what remains is controlled by a geometric
functional of the interface.

In this paper, we do not use the convergence of any massive fermionic observable. A recent preprint
by OpenAI~\cite[Thm.~2.1]{OAI26} proves the convergence of near-critical FK-Ising interfaces to a
massive $\SLE_{16/3}$ using the convergence of the massive FK-Ising observable of \cite{Par21} (more
details in \S\ref{ss:outlook}).

\subsection{The model, the interface and the space of curves}\label{ss:model} Let $n\ge 2$ be an
integer and $\delta:=1/n$. Let $\Omega_n:=\{1,\dots,n\}^2$. The \emph{outer boundary}
$\partial\Omega_n$ is the set of vertices of $\Z^2\setminus\Omega_n$ that have a nearest neighbor in
$\Omega_n$, and $E$ is the set of nearest-neighbor edges of $\Z^2$ with at least one endpoint in
$\Omega_n$; note that $E$ contains no edge between two vertices of $\partial\Omega_n$. We consider
the Ising model on the graph $(\Omega_n\cup\partial\Omega_n,E)$, with the spins of the outer
boundary vertices fixed by a boundary condition. The \emph{Dobrushin boundary condition} fixes the
spin of an outer boundary vertex $(i,j)$ to $+1$ if $i\le n/2$ and to $-1$ if $i>n/2$. The four
corner vertices of $[0,n+1]^2$ are not outer boundary vertices and carry no edge of $E$; since the
exploration defined below may inspect them in the four corner faces, we assign them spins by the
same rule. For a spin configuration $\sigma\in\{\pm1\}^{\Omega_n}$, extended to the outer boundary
by this boundary condition, we put
\begin{equation}\label{eq:H}
\varepsilon_e:=\sigma_x\sigma_y,\enspace \text{for }e=xy\in E,\qquad H(\sigma):=\sum_{e\in E}\varepsilon_e,\qquad
P_\beta(\sigma):=\frac{e^{\beta H(\sigma)}}{Z^\Dob_\beta},\enspace \beta\in\R.
\end{equation}
Here, $H$ is the negative of the Hamiltonian. The boundary condition changes sign at two places:
between the outer boundary vertices $(\lfloor n/2\rfloor,0)$ and $(\lfloor n/2\rfloor+1,0)$ at the
bottom, and between $(\lfloor n/2\rfloor,n+1)$ and $(\lfloor n/2\rfloor+1,n+1)$ at the top. We
denote these two places by $x_1=x_1^\delta$ and $x_2=x_2^\delta$. After scaling by $\delta$,
$(\Omega_n;x_1,x_2)$ approximates the \emph{Dobrushin square}
\[
(\Omega;x_1,x_2):=\bigl((0,1)^2;(\tfrac12,0),(\tfrac12,1)\bigr).
\]

\emph{The interface.} Faces are the unit squares of the lattice. A dual path is a sequence of faces
in which consecutive faces share an edge, and we draw it as a polygonal line through the midpoints
of the shared edges (\S\ref{ss:interface}). The interface $\gamma=\gamma(\sigma)$ is the dual path
from $x_1$ to $x_2$ that keeps the $+$ spins on its left. It is built by an exploration. At each
face, entered through one of its edges, let $F_L$ and $F_R$ be the two corners of the face ahead of
the path, on the left and on the right. If $\sigma_{F_L}=-1$, the path turns left; if
$\sigma_{F_L}=+1$ and $\sigma_{F_R}=-1$, it goes straight; and if $\sigma_{F_L}=\sigma_{F_R}=+1$, it
turns right. Every edge crossed by $\gamma$ therefore joins a $+$ vertex to a $-$ vertex. The
\emph{inspected vertices} $\Insp(\gamma)$ are the corners $F_L$ at every step, together with the
corners $F_R$ at the steps with $\sigma_{F_L}=+1$. A dual path $g$ is \emph{admissible} if
$g=\gamma(\sigma)$ for some $\sigma$. The turning rule only matters at faces whose four corners
alternate in sign, and, as we will explain in Remark~\ref{rem:turning}, the results below hold for
either of the two deterministic rules.

\begin{figure}[ht]\centering
\newcommand{\simpanel}[2]{\begin{minipage}[t]{0.32\textwidth}\centering\includegraphics[width=\textwidth]{figures/#1.png}\\[2pt]{\small #2}\end{minipage}}
\simpanel{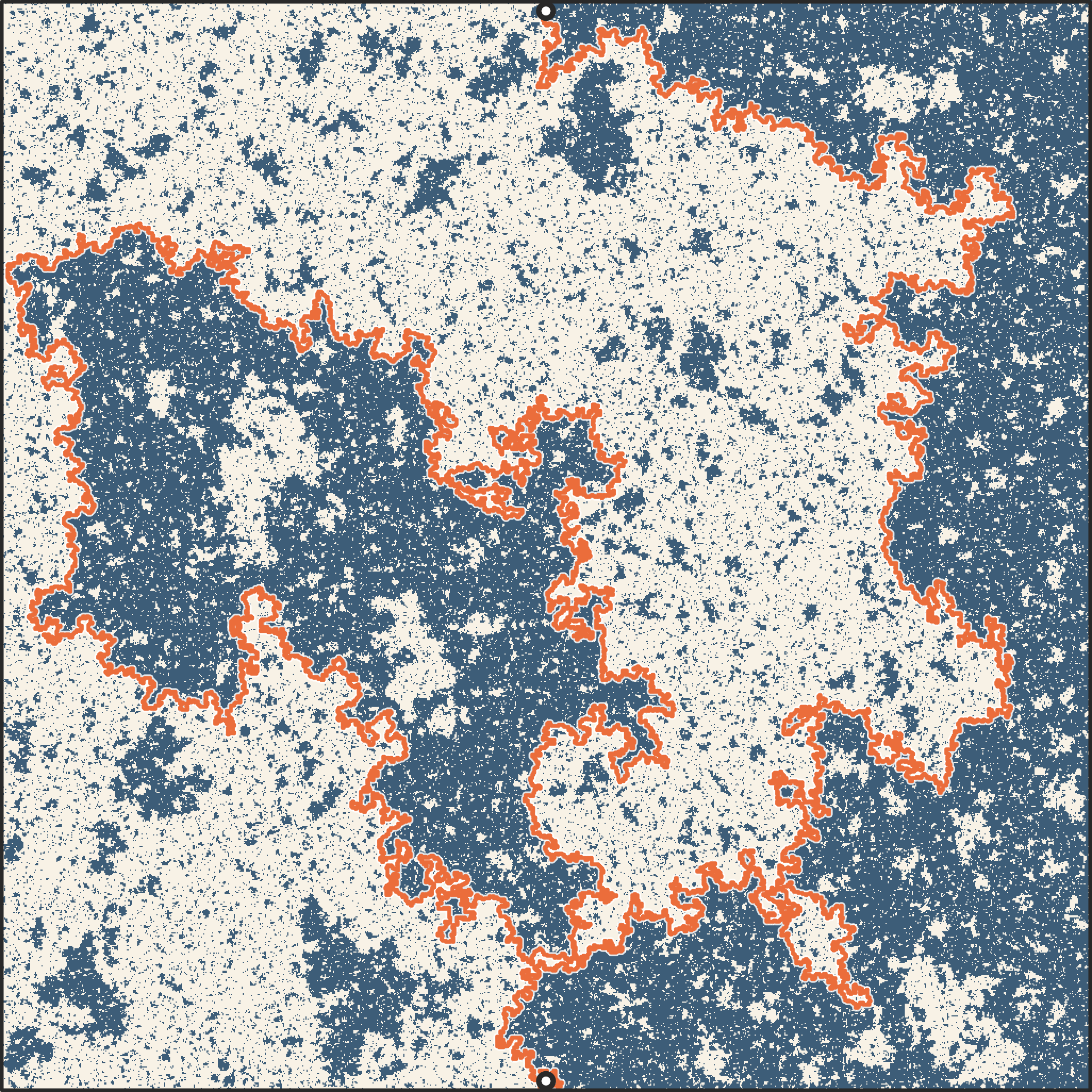}{(a) $\beta=\beta_c-m\delta$, $m=1.5$}\hfill
\simpanel{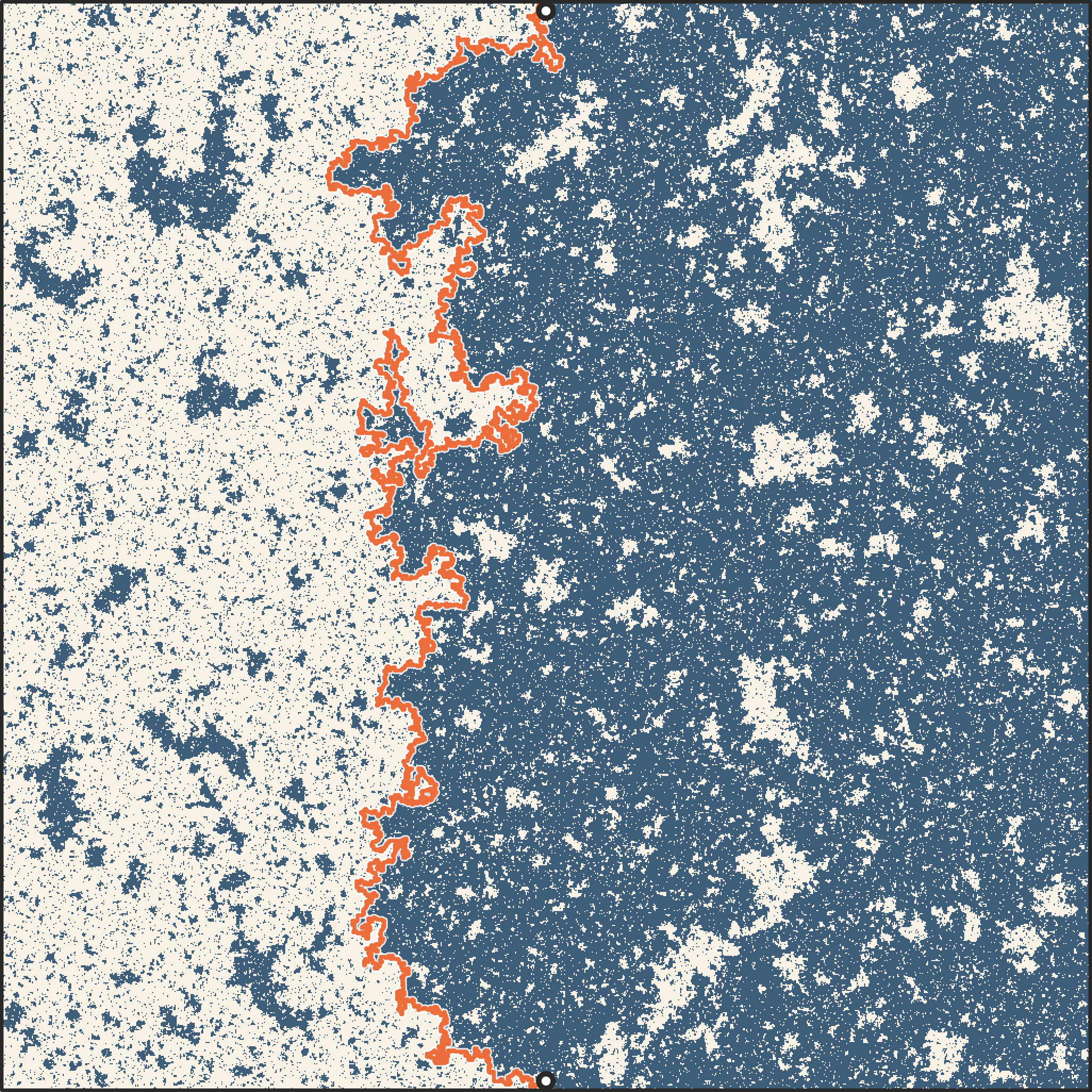}{(b) $\beta=\beta_c-m\delta$, $m=-1.5$}\hfill
\simpanel{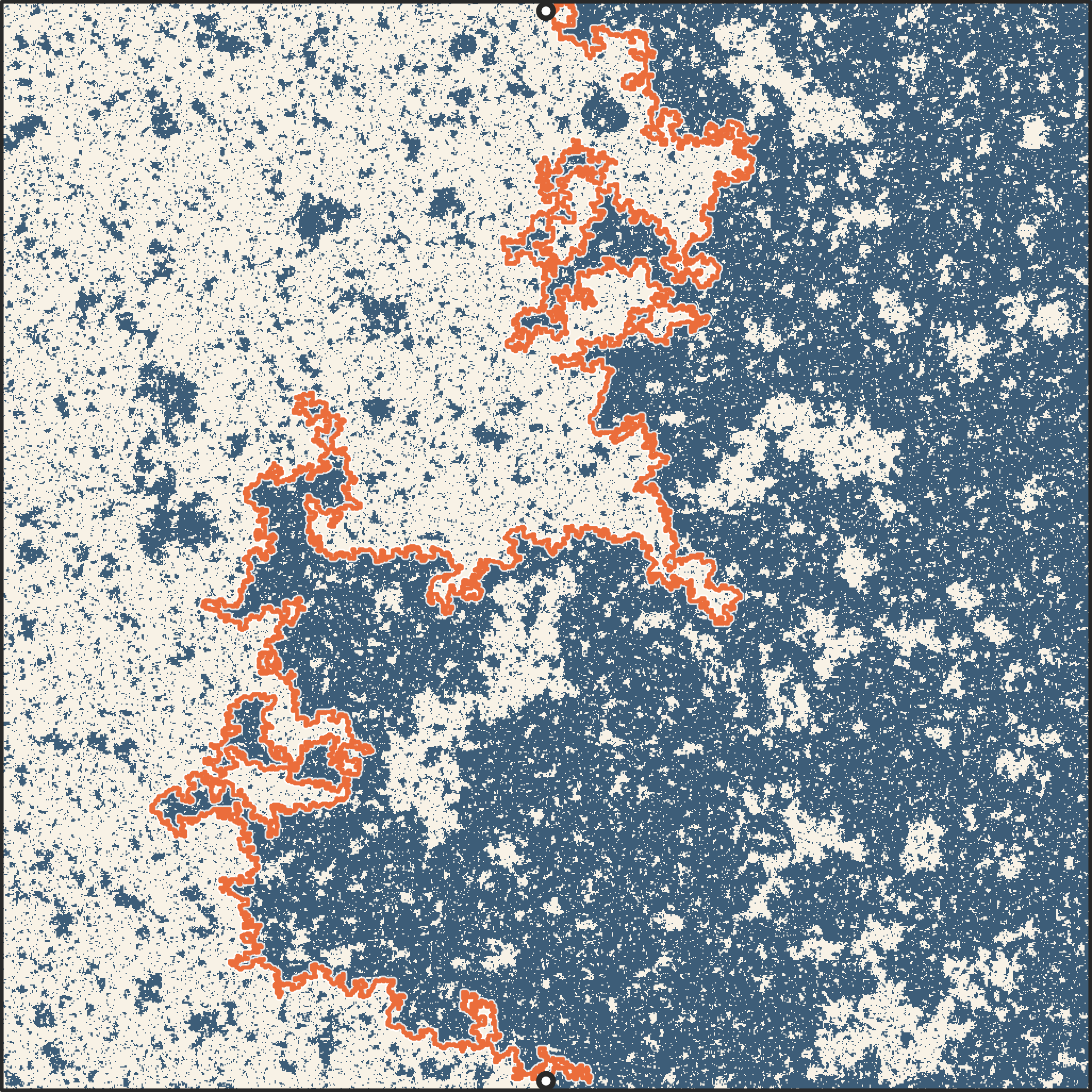}{(c) $\beta=\beta_c$}\\[8pt]
\simpanel{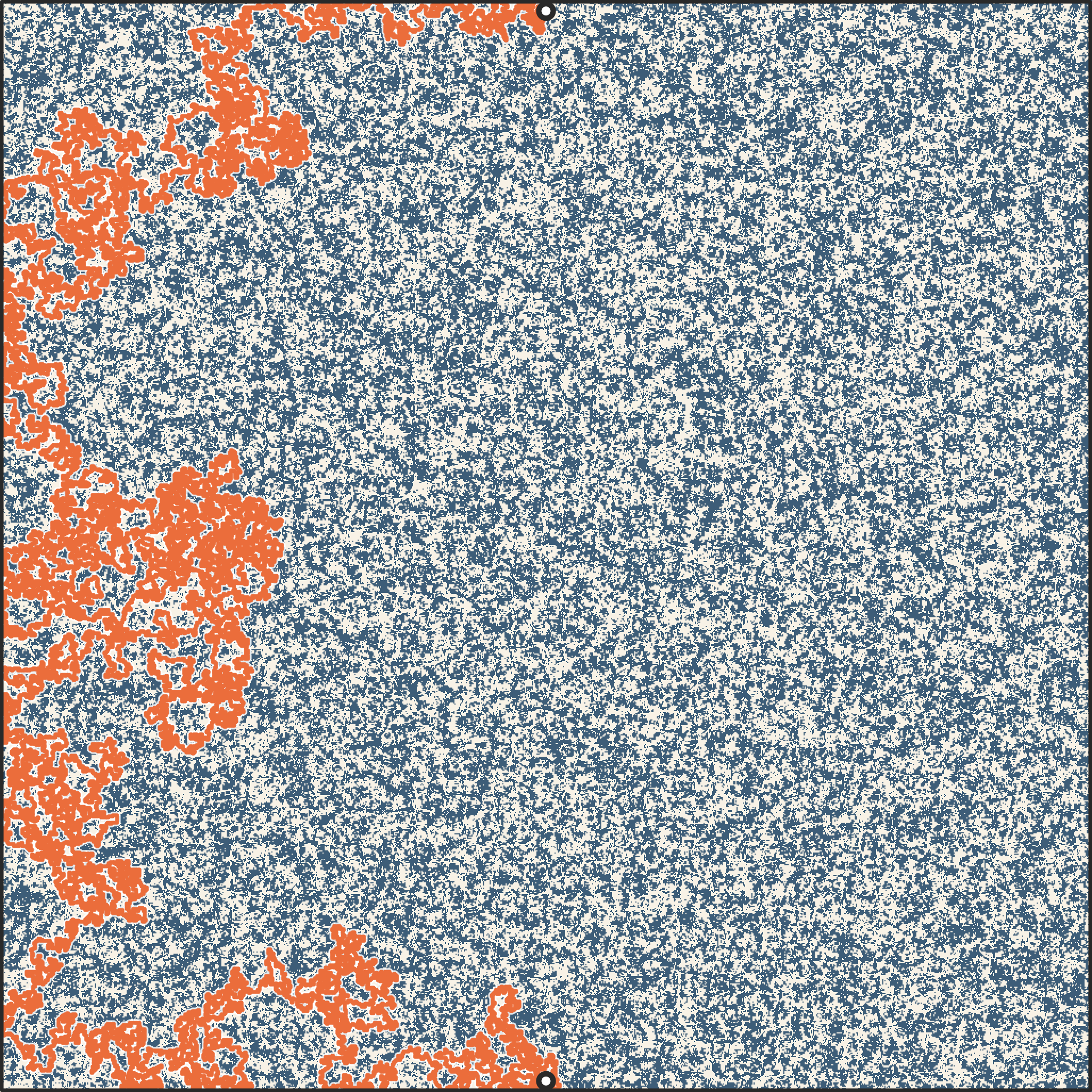}{(d) fixed $\beta=0.35<\beta_c$}\qquad
\simpanel{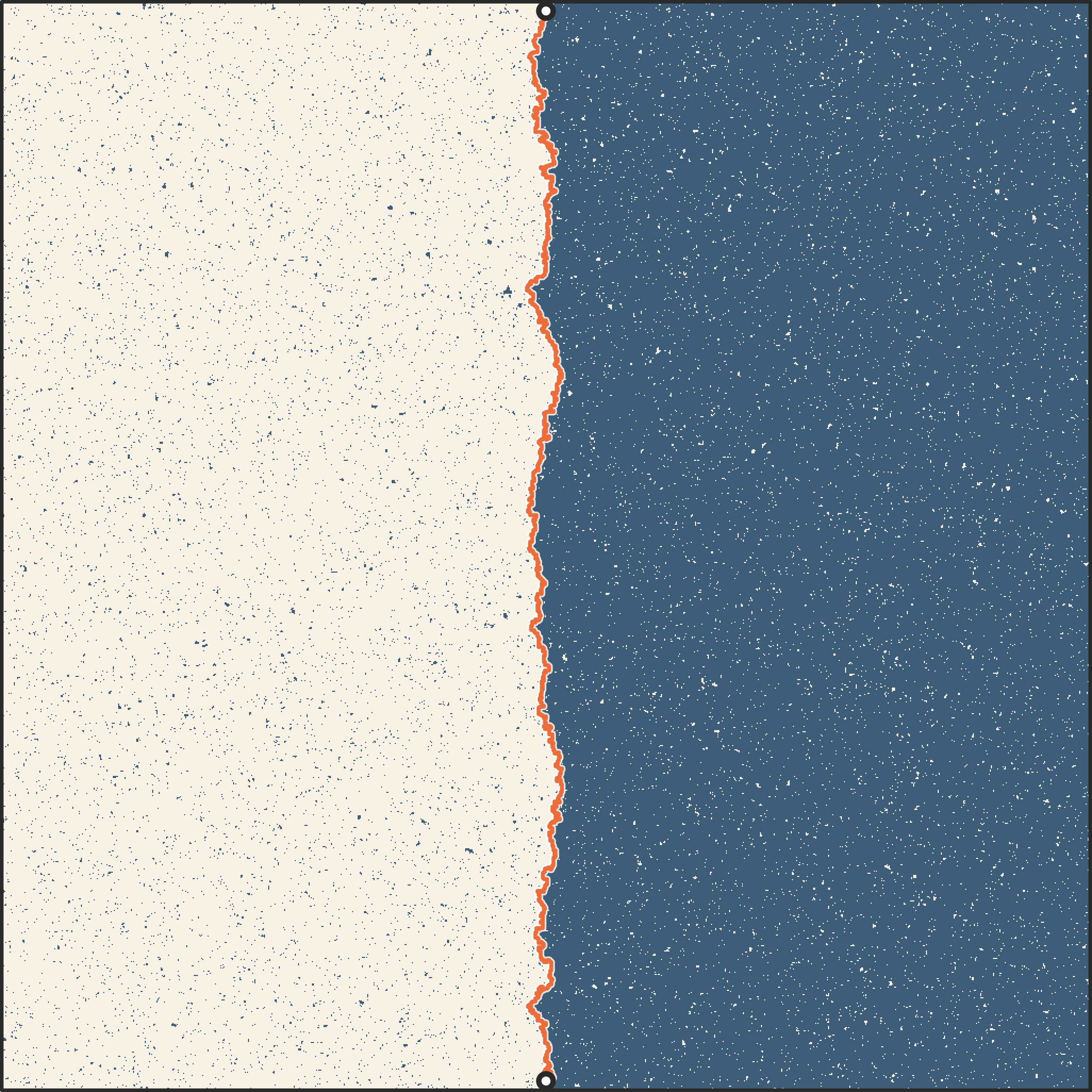}{(e) fixed $\beta=0.55>\beta_c$}
\caption{Spin configurations of the Ising model simulated in the Dobrushin square with $n=1024$ (so
that $\delta=1/1024$), with plus boundary spins on the left arc and minus boundary spins on the
right arc. Plus spins are light and minus spins are dark; the interface $\gamma$ from $x_1$ (bottom)
to $x_2$ (top) is drawn in orange, and the marked points are indicated by small circles. (a)
Near-critical, on the high-temperature side; (b) near-critical, on the low-temperature side; (c)
critical; (d) a fixed high temperature, independent of $\delta$; (e) a fixed low temperature,
independent of $\delta$.}\label{fig:sim}
\end{figure}

\emph{Temperatures.} Let $\beta_c=\frac12\log(1+\sqrt2)$ be the critical inverse temperature and fix
$m\in\R\setminus\{0\}$, independent of $\delta$. We consider the \emph{near-critical} regime
$\beta=\beta_c-m\delta$, which lies on the high-temperature side for $m>0$ and on the
low-temperature side for $m<0$. Figure~\ref{fig:sim} shows typical configurations in this regime, at
criticality, and at fixed temperatures away from criticality. Let $P^\delta_m$ and $P^\delta_c$ be
the laws of the curve $\delta\gamma$ under $P_{\beta_c-m\delta}$ and under $P_{\beta_c}$,
respectively. Both assign positive probability to every admissible path, so the Radon--Nikodym
derivative
\[
R_\delta(g):=\frac{P_{\beta_c-m\delta}(\gamma=g)}{P_{\beta_c}(\gamma=g)},\qquad \text{for }g\text{ admissible},
\]
is well defined and positive. We write $\E_\beta$ and $\Var_\beta$ for the expectation and the
variance under $P_\beta$, respectively, and $\E_c:=\E_{\beta_c}$.

\emph{The space of curves.} To compare the laws of the interface with $\SLE_3$, we view them as laws
of curves modulo reparametrization. Let $\cX$ be the set of continuous maps $[0,1]\to\R^2$ modulo
reparametrization, with the metric
\[
d(\gamma_1,\gamma_2):=\inf_{\varphi_1,\varphi_2}\ \sup_{s\in[0,1]}|\gamma_1(\varphi_1(s))-\gamma_2(\varphi_2(s))|,
\]
where the infimum is over increasing homeomorphisms of $[0,1]$. Here, two maps are identified when
$d(\gamma_1,\gamma_2)=0$. This identifies maps that differ by an increasing reparametrization, and
also maps that differ by inserting or removing intervals of constancy (which an increasing
homeomorphism alone cannot do), so that $d$ is a metric on $\cX$. The space $(\cX,d)$ is complete
and separable (see, e.g., \cite[\S2.1.1]{KS17}). Hence, every weakly convergent sequence of Borel
probability measures on $\cX$ is tight (see, e.g., \cite[Thm.~5.2]{Bil99}).

Let $\nu$ be the law of chordal $\SLE_3$ in $\Omega$ from $x_1$ to $x_2$, viewed as a probability
measure on $\cX$. By \cite[Thm.~1]{CDHKS14} (Theorem~\ref{thm:CDHKS} below), $P^\delta_c\to\nu$
weakly as $\delta\to0$.

\subsection{Absolute continuity of the law of the near-critical interface with respect to that of
\texorpdfstring{$\SLE_3$}{SLE3}}\label{ss:mainresult} Our main result states that the R\'enyi
divergences of every order $p>1$ and the Kullback--Leibler divergences, in both directions, between
the laws of the near-critical and critical interfaces remain bounded uniformly in $\delta$, and that
this property passes to the scaling limit.

For probability measures $\mu\ll\nu$, the \emph{R\'enyi divergence} of order $p>1$ and the
\emph{Kullback--Leibler divergence} are defined by
\[
\frac1{p-1}\log\E_\nu\Bigl[\Bigl(\frac{d\mu}{d\nu}\Bigr)^p\Bigr]\qquad\text{and}\qquad\KL(\mu\|\nu):=\E_\mu\Bigl[\log\frac{d\mu}{d\nu}\Bigr].
\]
In particular, $\KL(P^\delta_c\|P^\delta_m)=-\E_c[\log R_\delta]$ and
$\KL(P^\delta_m\|P^\delta_c)=\E_c[R_\delta\log R_\delta]$.
\begin{theorem}\label{thm:main}
Let $(\Omega;x_1,x_2)$ be the Dobrushin square and let $m\in\R\setminus\{0\}$ be fixed.
\begin{enumerate}[label=\textup{(\roman*)}]
\item For every $p>1$ and $s>0$,
\[
\sup_{n\ge2}\E_c\bigl[R_\delta^{\,p}\bigr]<\infty\qquad\text{and}\qquad\sup_{n\ge2}\E_c\bigl[R_\delta^{-s}\bigr]<\infty .
\]
Moreover, $\sup_{n\ge2}\KL(P^\delta_c\|P^\delta_m)<\infty$ and
$\sup_{n\ge2}\KL(P^\delta_m\|P^\delta_c)<\infty$.
\item The family $\{P^\delta_m\}_{\delta}$ is tight in $(\cX,d)$, and every subsequential weak limit
$\mu$ in $(\cX,d)$ is mutually absolutely continuous with respect to the law $\nu$ of chordal
$\SLE_3$. Moreover, $\frac{d\mu}{d\nu}\in L^p(\nu)$ and $\frac{d\nu}{d\mu}\in L^p(\mu)$ for every
$p<\infty$; in particular, $\KL(\mu\|\nu)<\infty$ and $\KL(\nu\|\mu)<\infty$.
\end{enumerate}
\end{theorem}

By Theorem~\ref{thm:main}(ii), every almost sure property of $\SLE_3$, such as being a simple curve
\cite{RS05} of Hausdorff dimension $\frac{11}8$ \cite{Bef08}, is shared by every subsequential limit
$\mu$. Part~(i) is particularly interesting because the corresponding R\'enyi and Kullback--Leibler
divergences for the whole spin configuration grow like $m^2\log n$, by the logarithmic growth of the
specific heat in a box of side $n$ \cite[Lemma~2.1]{GK25}, which goes back to Onsager~\cite{Ons44};
in this sense, the law of a single interface is less sensitive to the perturbation than the law of
the whole configuration. Part~(ii) follows from part~(i) and the convergence at $\beta_c$ by
standard measure-theoretic arguments (\S\ref{ss:limit}). The negative moments are proved in the same
way as the positive ones, with a kernel that is supported on both sides of $\beta_c$
(\S\ref{ss:renyiid}).

\subsection{Organization of the rest of the paper and discussion of the cancellation
mechanism}\label{ss:strategy} In Section~\ref{sec:prelim}, which can be consulted as needed, we will
fix the notation, recall the FK-Ising model and the conditional law of the spins given the
interface, and state the inputs. In Section~\ref{sec:renyi}, we will prove exact R\'enyi and
Kullback--Leibler identities for the law of the interface, reduce Theorem~\ref{thm:main}(i) to the
four estimates \tagref{E1}{(E1)}--\tagref{E4}{(E4)} below and deduce Theorem~\ref{thm:main} from
them (\S\S\ref{ss:assembly}--\ref{ss:limit}). Sections~\ref{sec:geom}, \ref{sec:cov}, \ref{sec:VP}
and~\ref{sec:mumu} will prove \tagref{E1}{(E1)}, \tagref{E2}{(E2)}, \tagref{E3}{(E3)} and
\tagref{E4}{(E4)}, respectively; Section~\ref{sec:cov} also contains the covariance estimates for
FK-Ising, uniform over rough domains, on which Section~\ref{sec:VP} relies. The two appendices will
deduce four-arm stability from \cite{DCM} (Appendix~\ref{app:R1}) and derive the strip identity from
\cite[\S3, Step~1]{BDC12}, as in the proof of \cite[Prop.~2.5]{DGP14} (Appendix~\ref{app:strip}).
Apart from \S\S\ref{ss:assembly}--\ref{ss:limit}, which assume \tagref{E1}{(E1)}--\tagref{E4}{(E4)},
each section uses only earlier ones and the appendices.

In the rest of this subsection, we sketch the main arguments of the proof of Theorem~\ref{thm:main}.
We do this for the reader's convenience and because the arguments themselves might be of independent
interest: they explain why the R\'enyi and Kullback--Leibler divergences remain bounded for a single
interface, although they are unbounded for the whole configuration.

\emph{The exact R\'enyi identity.} The proof controls the decrease of the variance of $H$ when the
interface is fixed,
\[
\cV_\beta(g):=\Var_\beta(H)-\Var_\beta(H\mid\gamma=g),\qquad\text{for }g\text{ admissible},
\]
whose expectation $\E_\beta[\cV_\beta(\gamma)]=\Var_\beta(\E_\beta[H\mid\gamma])$ is nonnegative,
although $\cV_\beta(g)$ may be negative for some $g$. Since $P_\beta(\sigma)$ is proportional to
$e^{\beta H(\sigma)}$, the second derivative of $\beta\mapsto\log P_\beta(\gamma=g)$ is
$-\cV_\beta(g)$. As a result, $\E_c[R_\delta^{\,a}]$ is exactly an exponential moment of an integral
$\int K_a(\beta)\,\cV_\beta(\gamma)\,d\beta$ for every real $a$, and $\KL(P^\delta_c\|P^\delta_m)$
is the expectation of such an integral. The kernels $K_a$ are explicit, nonnegative for
$a\notin[0,1]$ and of mass proportional to $(m\delta)^2=m^2/n^2$ (Proposition~\ref{prop:renyi}). The
first derivative of $\log P_\beta(\gamma=g)$ in $\beta$ cancels exactly in these identities; this is
essential, since it admits no useful bound for individual paths.

\emph{The cancellation of the bulk logarithm.} Both variances in $\cV_\beta$ grow like $n^2\log n$.
The key observation is that the logarithm is \emph{local}: it is carried by the \emph{local specific
heat} $\chi_\beta(e;r)$, that is, the sum of the covariances of $\varepsilon_e$ with the
$\varepsilon_f$ such that the midpoints of $e$ and $f$ are at Euclidean distance $|f-e|\le r$, in
the infinite-volume measure at $\beta$; $\chi_\beta(e;r)$ grows like $\log r$. In $\Omega_n$, this
sum is cut off at the Euclidean distance $d_\Omega(e)$ from the midpoint of $e$ to the outer
boundary; in the domain slit along a path $g$, it is cut off at the Euclidean distance $d_g(e)$ from
the midpoint of $e$ to $g$, viewed as a polygonal line. We compare both domains with the \emph{bulk
reference} $X_\beta:=\sum_e\chi_\beta(e;r_\Omega(e))$, where $r_\Omega(e):=1\vee\frac13d_\Omega(e)$.
The main ingredient is that the local specific heat at radius $r$, computed in a domain with plus
boundary condition, depends on the domain only up to an \emph{additive} constant, as long as the
domain contains the disk of radius $2r$ around $e$ (Theorem~\ref{thm:LSS}). In $\Omega_n$, this
applies at radius $r_\Omega(e)$, so that the logarithm carried by each edge is the same as in
$X_\beta$ and cancels; in the slit domain, it applies only at a radius comparable to
$1\vee\min(d_g(e),d_\Omega(e))$, and the loss is at most $C(1+\log_+[d_\Omega(e)/(1\vee d_g(e))])$
(Figure~\ref{fig:cancel}). Summing over the edges, we find that what remains is bounded by $Cn^2$
times the \emph{proximity functional}
\[
\Pi_\delta(g):=1+\Phi(g),\quad\text{where}\quad\Phi(g):=\delta^2\sum_{e\in E}\log_+\frac{d_\Omega(e)}{1\vee d_g(e)},
\]
which is large only if $g$ comes close to many edges. In macroscopic units, $\Phi$ is, up to a
factor $2$, a Riemann sum for the integral over $\Omega$ of the positive part of the logarithm of
the ratio of the distances to $\partial\Omega$ and to the curve. Since the limiting curve has
Hausdorff dimension $\frac{11}8<2$, $\Pi_\delta$ is typically bounded, and Russo--Seymour--Welsh
(RSW) estimates control its exponential moments.

\begin{figure}[ht]\centering
\begin{tikzpicture}[scale=0.5]
\draw[line width=1.1pt] (0,0) rectangle (12,12);
\draw[line width=1.4pt] plot[smooth,tension=0.7] coordinates {(6,0) (5.2,1.8) (6.0,3.2) (4.7,4.6) (4.5,6.3) (5.6,7.8) (7.2,9.0) (6.2,10.6) (6,12)};
\fill (6,6) circle (2.5pt); \node[above] at (6,6.45) {\small $e$};
\draw[dashed] (6,6) circle (2.0);
\draw[densely dotted,thick] (6,6) circle (0.5);
\draw[->] (6,6) -- (8,6); \node[right] at (8,6) {\scriptsize $r_\Omega(e)$};
\draw[->] (6,6) -- (6,5.5); \node[below,inner sep=1pt] at (6,5.45) {\scriptsize $\asymp d_g(e)$};
\node[left] at (6.3,10.0) {\small $g$};
\node at (6,-0.7) {\small $x_1$}; \node at (6,12.7) {\small $x_2$};
\end{tikzpicture}
\caption{The cancellation of the bulk logarithm. The local specific heat of $e$ grows like the
logarithm of the radius. In the bulk reference $X_\beta$, it is cut off at $r_\Omega(e)\asymp
d_\Omega(e)$. Given the interface $g$, the lower bound \tagref{E2}{(E2)} keeps only the edges $f$
within a third of the distance from $e$ to the pinned vertices (inspected or outer boundary
vertices, \S\ref{ss:LSS}), which is comparable to $d_g(e)$. The loss, at most
$C(1+\log_+[d_\Omega(e)/(1\vee d_g(e))])$, is the corresponding summand of $n^2\Pi_\delta(g)$, and
the bulk logarithm $\log d_\Omega(e)$ cancels.}\label{fig:cancel}
\end{figure}
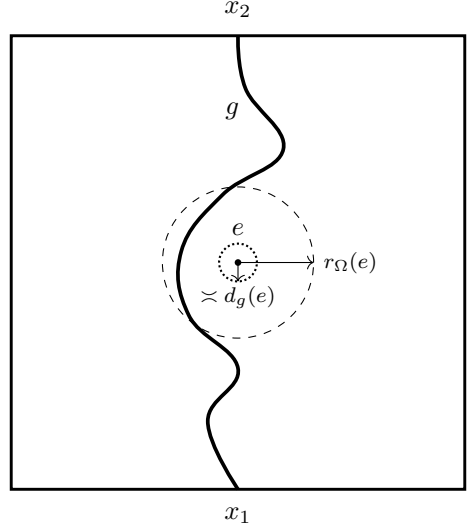

\begin{figure}[ht]\centering
\begin{tikzpicture}[x=1cm,y=1cm,font=\scriptsize,
  nd/.style={draw,rounded corners=2pt,align=flush center,inner sep=2.5pt,text width=2.75cm,minimum height=1.0cm,
    execute at begin node={\hyphenpenalty=10000\exhyphenpenalty=10000}},
  ok/.style={nd,thick,fill=white},
  ct/.style={nd,densely dashed,fill=gray!15},
  cl/.style={nd,thick,densely dotted,fill=white},
  ar/.style={-{Stealth[length=1.6mm]},black!75}]
\node[ok,text width=6.2cm,minimum height=0.85cm] (main) at (0,0)
  {\textbf{Theorem~\ref{thm:main}}\\{}discrete bounds and limit, \S\S\ref{ss:assembly}--\ref{ss:limit}};
\node[ok,text width=6.2cm,minimum height=0.85cm] (S1) at (0,-1.5)
  {\textbf{Step 1}: exact identities\\{}and the decomposition \eqref{eq:reduction}\\{}Section~\ref{sec:renyi}};
\node[ok] (S2) at (-4.8,-3.2) {\textbf{Step 2} \tagref{E1}{(E1)}\\{}exponential moments of $\Pi_\delta$\\{}Prop.~\ref{prop:expmom}};
\node[ok] (S3) at (-1.6,-3.2) {\textbf{Step 3} \tagref{E2}{(E2)}\\{}lower bound on the conditional variance\\{}Prop.~\ref{prop:Bminus}};
\node[ok] (S4) at (1.6,-3.2) {\textbf{Step 4} \tagref{E3}{(E3)}\\{}specific heat of the plus domain\\{}Thm.~\ref{thm:VP}};
\node[ok] (S5) at (4.8,-3.2) {\textbf{Step 5} \tagref{E4}{(E4)}\\{}boundary disorder correlation\\{}Thm.~\ref{thm:mumu}};
\node[ok,text width=2.45cm] (S6) at (-2.3,-5.7) {\textbf{Step 6}\\{}local specific-heat stability\\{}Thm.~\ref{thm:LSS}};
\node[cl,text width=4.5cm] (S7) at (2.06,-5.7) {\textbf{four-arm stability}\\{}deduced from \mbox{\cite{DCM}}\\{}(Thm.~\ref{thm:R1}, App.~\ref{app:R1})};
\node[cl,text width=3.7cm] (Q2) at (6.99,-5.7) {\textbf{strip identity}\\{}row relation in \mbox{\cite[\S3, Step~1]{BDC12}}; duality in the proof of \mbox{\cite[Prop.~2.5, p.~15]{DGP14}}\\{}(Thm.~\ref{thm:strip}, App.~\ref{app:strip})};
\node[ct,text width=7.0cm,anchor=north] (P1) at (0.89,-7.4) {\textbf{mixing rate and arm events}\\{}ratio mixing \mbox{\cite[Prop.~2.9]{DCM}};\\{}stability of the mixing rate \mbox{\cite[Thm.~1.6(iii)]{DCM}};\\{}$\Delta_{\ppc}(R)\asymp R^{-1}$ \mbox{\cite[Thm.~1.13]{DCM}};\\{}quasi-multiplicativity \mbox{\cite[Thm.~1.3]{CDH16}}, \mbox{\cite[Thm.~1.6(ii)]{DCM}};\\{}comparison to pivotality and derivative bounds \mbox{\cite[Thm.~1.6(iv), Cor.~7.1]{DCM}}\\{}(Thm.~\ref{thm:mix}, Lemma~\ref{lem:Delta}, App.~\ref{app:R1})};
\node[ct,text width=3.5cm,anchor=north] (P3) at (-4.8,-7.4) {\textbf{RSW near criticality}\\{}\mbox{\cite[Thm.~1.3]{DGP14}};\\{}at $\ppc$: \mbox{\cite[Thm.~1]{DHN11}}\\{}(Thm.~\ref{thm:RSW})};
\node[ct,text width=3.7cm] (P4) at (6.99,0) {\textbf{convergence at $\beta_c$}\\{}\mbox{\cite[Thm.~1]{CDHKS14}}\\{}(Thm.~\ref{thm:CDHKS})};
\coordinate (bot) at (0,-11.0);
\coordinate (ra) at (4.75,0);
\draw[ar] (S1) -- (main);
\draw[ar] (P4) -- (main);
\draw[ar] (S2.north) |- (S1.west);
\draw[ar] (S3.north) -- (S3.north |- S1.south);
\draw[ar] (S4.north) -- (S4.north |- S1.south);
\draw[ar] (S5.north) |- (S1.east);
\draw[ar] (S6.north) -- (S3.south);
\draw[ar] (S6.north east) -- (S4.south);
\draw[ar] (S7) -- (S6);
\draw[ar] (S7.north) -- (S4.south east);
\draw[ar] (P1.north -| S6.south) -- (S6.south);
\draw[ar] (P1.north -| S7.south) -- (S7.south);
\draw[ar] (Q2.north) -- (S5.south east);
\draw[ar] (P3.south) -- (P3.south |- bot) -- (ra |- bot) -- (ra |- S5.south);
\draw[ar] (P3.north) -- (S2.south);
\end{tikzpicture}
\caption{Roadmap of the proof. Arrows point from an ingredient to the statement that uses it.
\emph{Solid boxes}: results proved in this paper. \emph{Dashed gray boxes}: results stated and
proved in the literature; their statements and hypotheses will be recalled in
\S\S\ref{ss:FK}--\ref{ss:inputs} and Appendix~\ref{app:R1}. \emph{Dotted boxes}: results deduced in
Appendices~\ref{app:R1} and~\ref{app:strip} from the literature. Only the principal arrows are
drawn; RSW will also be used in Steps~4 and~6, and the mixing rate in Step~4.}\label{fig:diagram}
\end{figure}
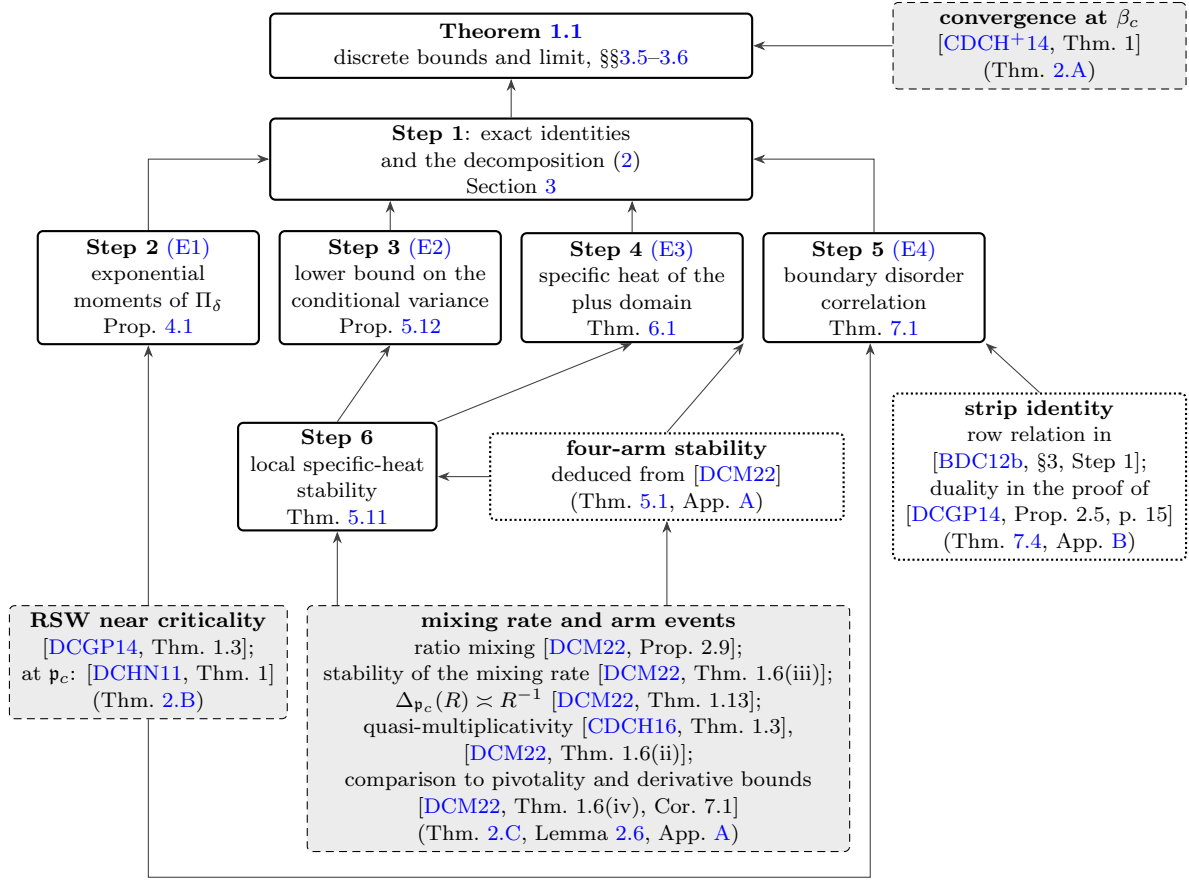

\enlargethispage{\baselineskip}
\emph{The reduction to four estimates.} To implement this idea, we first separate the effect of the
Dobrushin boundary condition. Let $\Var^+_\beta$ and $Z^+_\beta$ be the variance and the partition
function with plus boundary condition, let $\phi^1_\beta$ be the FK measure on
$(\Omega_n\cup\partial\Omega_n,E)$ with wired boundary condition at the parameter $\pp(\beta)$ and
let $\Dab$ be the event that no open path joins the two boundary arcs. Then
$\phi^1_\beta(\Dab)=Z^\Dob_\beta/Z^+_\beta$ is the \emph{boundary disorder correlation}, and
$\ell(\beta):=\log\phi^1_\beta(\Dab)$ satisfies $\Var_\beta(H)=\Var^+_\beta(H)+\ell''(\beta)$. In
Section~\ref{sec:renyi}, we will deduce the following exact decomposition: for every admissible $g$
and each kernel $K_a$, where $a\in\{p,-s\}$ is the order of the moment in Theorem~\ref{thm:main}, or
$a={\rm KL}$ for the Kullback--Leibler divergence,
\begin{equation}\label{eq:reduction}
\begin{split}
\int K_a(\beta)\,\cV_\beta(g)\,d\beta&=\int K_a(\beta)\bigl[\Var^+_\beta(H)-X_\beta\bigr]\,d\beta+\Delta^2_a\ell\\
&\quad+\int K_a(\beta)\bigl[X_\beta-\Var_\beta(H\mid\gamma=g)\bigr]\,d\beta,
\end{split}
\end{equation}
where $\Delta^2_a\ell$ is the second difference of $\ell$ defined in \eqref{eq:seconddiff}. The
proof then reduces to four estimates, uniform in $n$, in admissible $g$ and in
\[
\beta\in W_M:=[\beta_c-M\delta,\beta_c+M\delta],
\]
where $M:=\max(p,s,1)\,|m|$. Two of them concern the curve, and two concern the fixed square:
\begin{equation}\label{eq:inputs}
\begin{aligned}
&\text{\tagdef{E1}{(E1)}}\qquad \E_\beta\bigl[e^{\eta\Pi_\delta(\gamma)}\bigr]\le C(M,\eta),\quad\text{for every }\eta>0;\\
&\text{\tagdef{E2}{(E2)}}\qquad \Var_\beta(H\mid\gamma=g)\ge X_\beta-Cn^2\,\Pi_\delta(g);\\
&\text{\tagdef{E3}{(E3)}}\qquad \Var^+_\beta(H)\le X_\beta+Cn^2;\\
&\text{\tagdef{E4}{(E4)}}\qquad \phi^1_\beta(\Dab)\le C\,\phi^1_{\beta_c}(\Dab),\quad\text{for }\beta\le\beta_c.
\end{aligned}
\end{equation}
For $m>0$, since $\ell$ is nonincreasing, $\Delta^2_a\ell$ is bounded above by a combination of
increments $\ell(\beta)-\ell(\beta_c)$ with $\beta<\beta_c$, also for the negative moments, whose
kernel straddles $\beta_c$. Hence, \tagref{E4}{(E4)}, a comparison of $\ell$ itself up to an
additive constant, which involves no derivative of $\ell$, bounds $\Delta^2_a\ell$ by a constant.
Together with \tagref{E2}{(E2)} and \tagref{E3}{(E3)}, this gives the \emph{pathwise bound}, valid
for every admissible path $g$,
\[
\int K_a(\beta)\,\cV_\beta(g)\,d\beta\ \le\ C+C\,\|K_a\|_1\,n^2\,\Pi_\delta(g),\qquad \text{with }\|K_a\|_1\,n^2\asymp m^2 .
\]
By \tagref{E1}{(E1)}, this gives the discrete bounds of Theorem~\ref{thm:main}(i), and the passage
to the limit uses only standard arguments (\S\ref{ss:limit}). The more challenging part is to prove
\tagref{E2}{(E2)} and \tagref{E3}{(E3)}, that is, to compare the slit domain and the plus domain
with the same bulk reference $X_\beta$, to accuracy $O(n^2\Pi_\delta(g))$ and $O(n^2)$,
respectively, and uniformly over rough domains.

\emph{Structure of the proof.} Figure~\ref{fig:diagram} summarizes the proof in six steps, together
with the inputs on which they rely. Step~1 is the reduction of Section~\ref{sec:renyi}, and
Steps~2--5 prove \tagref{E1}{(E1)}--\tagref{E4}{(E4)}, respectively. Steps~3 and~4 rest on two
FK-Ising tools. The first is Step~6, \emph{local specific-heat stability} (Theorem~\ref{thm:LSS}):
the sum of the covariances of $\varepsilon_e$ with the $\varepsilon_f$, $|f-e|\le r$, under the plus
boundary condition in a domain $D$ differs from $\chi_\beta(e;r)$ by at most $C$ whenever $D$
contains the disk of radius $2r$ around $e$, uniformly over domains, rough slit domains included.
The second is \emph{four-arm stability} (Theorem~\ref{thm:R1}): the probability of four alternating
arms, open and dual-open, across the annulus between radii $r$ and $R$ is at most
$C(r/R)^{1+\varepsilon_1}$, for a universal constant $\varepsilon_1>0$, whenever the FK parameter is
within $\lambda/R$ of its critical value, with $C$ depending on $\lambda$. The exponent greater than
$1$ keeps the dyadic sums of Steps~4 and~6 free of logarithms.

\emph{The low-temperature side.} The same proof gives Theorem~\ref{thm:main} for $m<0$, that is, at
$\beta=\beta_c+|m|\delta$. Estimates \tagref{E1}{(E1)}, \tagref{E2}{(E2)} and \tagref{E3}{(E3)} hold
on both sides of $\beta_c$. The control of $\Delta^2_a\ell$ then also requires the lower bound
$\phi^1_\beta(\Dab)\ge C^{-1}\phi^1_{\beta_c}(\Dab)$ for $\beta>\beta_c$, which
Theorem~\ref{thm:mumu} gives as well, and $\KL(P^\delta_c\|P^\delta_m)$ is bounded through the
negative moments by Jensen's inequality. For definiteness, the proofs are written for $m>0$.

We end this subsection by emphasizing that the constants in Theorem~\ref{thm:main} depend on $m$,
$p$ and $s$, and that they might be unbounded as $|m|\to\infty$.

\subsection{Outlook}\label{ss:outlook} We conclude this section with remarks on related questions
and on possible extensions of our results.

\emph{Uniqueness of the subsequential limit.} Theorem~\ref{thm:main} does not identify $\mu$ or
assert that the limit is unique. The convergence in law of the interface and the description of its
driving function are the subject of a paper in preparation \cite{Fen26}.

\emph{Other domains.} It is natural to ask how far the assumption that the domain is a square can be
relaxed. Several ingredients do not involve the square: the convergence at $\beta_c$
(Theorem~\ref{thm:CDHKS}) holds for general Dobrushin domains; four-arm stability
(Theorem~\ref{thm:R1}) is a statement about the full plane; the strip identity
(Theorem~\ref{thm:strip}) concerns an auxiliary strip; and the argument of Section~\ref{sec:cov}
holds in arbitrary domains. The square enters the proof mainly in two places.

First, Section~\ref{sec:VP} uses the bounds \tagref{Geo}{(Geo)} and \tagref{Geop}{(Geo$'$)} on the
number of edges near the boundary, the closed dual half-circuits of \tagref{Rb}{(Rb)} and the fact
that every outer boundary vertex has a single edge in $E$ (\S\ref{ss:geo}). The bounds and the
half-circuits extend to \emph{rectilinear} domains, that is, domains bounded by a fixed number of
axis-parallel segments of length comparable to $n$, any two nonadjacent ones being at distance
comparable to $n$, and we expect RSW in boxes meeting the boundary to give the half-circuits for
smooth domains, for which the bounds are elementary. The single-edge property fails at the reflex
corners of the boundary, where an outer boundary vertex has two edges in $E$, and the separation
argument of \S\ref{ss:geo} would have to be adapted there; however, this does not seem to be an
essential obstruction.

Second, and this is the main obstruction, the gluing at the marked points (\S\ref{ss:glue}) reduces
the boundary disorder correlation to connection probabilities in a strip, compared across parameters
through the strip identity. This requires the boundary to be straight and axis-parallel within
distance comparable to $n$ of $x_1$ and $x_2$. We therefore expect our results to hold, with the
same framework, for every rectilinear or smooth (or piecewise smooth) Dobrushin domain with this
property; for a rectilinear domain, this means that $x_1$ and $x_2$ are at distance comparable to
$n$ from its corners. General smooth domains would require a version of the strip identity, or of
Theorem~\ref{thm:mumu}, near a marked point on a curved or tilted boundary, where the translation
invariance and the self-duality of the strip (Appendix~\ref{app:strip}) are lost; we do not pursue
this here.

In summary, besides the convergence at $\beta_c$, the proof uses the following four properties of
the domain:
\begin{enumerate}[label=\textup{(\roman*)}]
\item the bounds \tagref{Geo}{(Geo)} and \tagref{Geop}{(Geo$'$)} on the number of edges near the
boundary hold;
\item the closed dual half-circuits of \tagref{Rb}{(Rb)} exist near the boundary;
\item every outer boundary vertex has a single edge in $E$;
\item the boundary is straight and axis-parallel within distance comparable to $n$ of $x_1$ and
$x_2$.
\end{enumerate}
The square has all four. Rectilinear domains have (i) and (ii), and (iv) when $x_1$ and $x_2$ are at
distance comparable to $n$ from their corners, but (iii) fails at their reflex corners. Smooth
domains have (i), and we expect them to have (ii), but (iii) fails wherever their boundary is not
axis-parallel.

\emph{The FK-Ising interface.} A recent preprint by OpenAI~\cite[Thm.~2.1]{OAI26} proves the
convergence of near-critical FK-Ising interfaces on isoradial lattices, for either sign of the mass,
to a massive $\SLE_{16/3}$ characterized by its drift; the absolute continuity of the law of the
complete curve conjectured in \cite[Conj.~2]{GK25} is not part of that result (see
\cite[Remark~6.4]{OAI26}). Its proof uses the convergence of the massive FK-Ising observable of
\cite{Par21} (see \cite[\S3.1]{OAI26}), whereas the present paper uses the massive observable only
through the exact relation between consecutive rows of a straight strip (Appendix~\ref{app:strip}).

\emph{Singularity of the whole configuration.} The law of the whole near-critical FK-Ising
configuration is singular with respect to the critical law \cite[Thm.~2.5]{GK25}, and the same is
expected for the whole spin configuration~\cite[\S2.4]{GK25} and \cite[Open Problem~2]{CGS26}. This
does not contradict Theorem~\ref{thm:main}: the (conjectured) singularity of the law of the whole
spin configuration may come from the logarithmic growth of the specific heat, to which all the edges
of the domain contribute, whereas, for a single interface, this bulk logarithm cancels edge by edge
between the full domain and the slit domain (Theorem~\ref{thm:VP} and
Proposition~\ref{prop:Bminus}), and only the edges close to the interface contribute. For the spin
configuration, the separating events of \cite{GK25} cannot be used directly, since they are
measurable with respect to the FK-Ising configuration, which is not a function of the spin
configuration.

\FloatBarrier

\section{Preliminaries}\label{sec:prelim} In this section, we collect the notation and inputs on
which the proof relies, so that the reader may consult it as needed. Until the end of
Section~\ref{sec:mumu}, we assume that $m>0$; for the case $m<0$, see the discussion at the end of
\S\ref{ss:strategy}. The mass $m$ and the moment orders $p>1$ and $s>0$ are fixed, and we recall
from \S\ref{ss:strategy} that
\[M=\max(p,s,1)\,|m|=\max(p,s,1)\,m.\]
All the inverse temperatures that occur lie in the window
\[W_M=[\beta_c-M\delta,\beta_c+M\delta],\]
where $\delta=1/n$ is the mesh. Constants $C,c\in(0,\infty)$ may change from line to line. Unless
stated otherwise, they depend only on $M$ (and on $p,s$ where indicated). We write $a\lesssim b$, or
$b\gtrsim a$, if $a\le Cb$, and $a\asymp b$ if $a\lesssim b$ and $b\lesssim a$. We write $|\cdot|_2$
and $|\cdot|_\infty$ for the Euclidean and the $\ell^\infty$ norms on $\R^2$, and $\dist_2$ and
$\dist_\infty$ for the corresponding distances between subsets of $\R^2$. In Theorem~\ref{thm:R1}
and Appendix~\ref{app:R1}, the constants are universal, except that they may depend on $\lambda$
where indicated. Notation introduced inside a proof is local to that proof.

\subsection{Notation and the parameter window}\label{ss:notation} For an edge $e$ of $\Z^2$, we
denote by $x_e$ its left endpoint if $e$ is horizontal and its bottom endpoint if $e$ is vertical,
and by $\mathrm{mid}(e)$ its midpoint. We use the distances
$|e-f|:=|\mathrm{mid}(e)-\mathrm{mid}(f)|_2$ and $|e-f|_\infty:=|x_e-x_f|_\infty$, which satisfy
$|e-f|_\infty\le|e-f|+1$ and $|e-f|\le\sqrt2\,|e-f|_\infty+1$. Boxes are the sets
\[
\Lambda_r(x):=\{y\in\Z^2:|y-x|_\infty\le r\},
\]
with inner boundary $\partial\Lambda_r(x):=\{y:|y-x|_\infty=r\}$, and
$\Lambda_r(e):=\Lambda_r(x_e)$, $\Lambda_r:=\Lambda_r(0)$. For a set $V$ of vertices, $E(V)$ is the
set of edges with both endpoints in $V$, and $\cF(V)$ is the $\sigma$-algebra generated by the
states of the edges of $E(V)$. For $r<R$, we put
\[
E_{r,R}(x):=E(\Lambda_R(x))\setminus E(\Lambda_r(x)).
\]
For a vertex $x$, we put $a_x:=\dist_\infty(x,\Z^2\setminus\Omega_n)$, so that $a_x=0$ off
$\Omega_n$, and $\hat a_x:=1\vee a_x$; for an edge, $a_e:=a_{x_e}$ and $\hat a_e:=\hat a_{x_e}$. The
Euclidean distances of an edge to the outer boundary and to a dual path $g$, viewed as a polygonal
line, and the radius of \S\ref{ss:strategy} are
\[
d_\Omega(e)=\dist_2\bigl(\mathrm{mid}(e),\Z^2\setminus\Omega_n\bigr),\qquad r_\Omega(e)=1\vee\tfrac13d_\Omega(e),\qquad
d_g(e)=\dist_2\bigl(\mathrm{mid}(e),g\bigr).
\]
In the square, every vertex of $\Z^2\setminus\Omega_n$ nearest to $x_e$ is reached along a
coordinate axis and $|\mathrm{mid}(e)-x_e|=\frac12$, hence $|d_\Omega(e)-a_e|\le\frac12$ for $e\in
E$.

\emph{The window $W_M$ and the largest scale $N$.} Recall that $\pp(\beta)=1-e^{-2\beta}$ and that
$\ppc=\pp(\beta_c)$ is the critical point of the FK-Ising model \cite{BDC12a}. Since
$0<\pp'(\beta)\le2$ on $[0,\infty)$, we have $|\pp(\beta)-\ppc|\le2|\beta-\beta_c|$ for $\beta\ge0$.
We put
\begin{equation}\label{eq:N}
N:=4n,\qquad n_0(M):=\lceil2M/\beta_c\rceil .
\end{equation}
For $n\ge n_0(M)$ and $\beta\in W_M$, we have $\beta\in[\frac12\beta_c,\frac32\beta_c]$, hence
$\pp(\beta)\in[\frac13,\frac34]$ and
\begin{equation}\label{eq:window}
|\pp(\beta)-\ppc|\le 2M\delta=\frac{8M}N .
\end{equation}
The number $N$ is the largest scale of the problem: all distances between edges of $E$ and all
scales used below are at most $N$. For the finitely many $n<n_0(M)$, $R_\delta$ is a positive
function on a finite set, since \eqref{eq:H} defines $P_\beta$ for every real $\beta$, so the
discrete bounds of Theorem~\ref{thm:main}(i) are trivial. From now on, we assume that $n\ge n_0(M)$,
and the same remark applies to the larger thresholds, such as $n_1(M)$, that will be introduced
later.

\subsection{FK-Ising toolbox}\label{ss:FK} Let $G$ be a finite graph with vertex set
$V(G)\subseteq\Z^2$, and let $\partial G\subseteq V(G)$ contain every vertex incident to an edge of
$\Z^2$ not in $E(G)$. We write $G\supseteq\Lambda_s(x)$ if $\Lambda_s(x)$, with all its edges, is a
subgraph of $G$ and no vertex of $\Lambda_s(x)\setminus\partial\Lambda_s(x)$ lies in $\partial G$. A
\emph{boundary condition} $\xi$ is a partition of $\partial G$. The wired one, $\xi=1$, has one
block; the free one, $\xi=0$, has only singletons; and $\xi\le\xi'$ means that $\xi'$ is coarser
than $\xi$. The FK-Ising measure is
\[
\phi^\xi_{G,\pp}(\omega)\ \propto\ \pp^{o(\omega)}(1-\pp)^{c(\omega)}\,2^{k(\omega^\xi)},\qquad\omega\in\{0,1\}^{E(G)},
\]
where $o$ and $c$ count the open and closed edges, and $k(\omega^\xi)$ is the number of clusters
once the blocks of $\xi$ are identified. We use the following classical facts (see, e.g.,
\cite[Ch.~3--4]{Gri06}).
\begin{enumerate}[leftmargin=3.6em]
\item[\tagdef{DMP}{(DMP)}] \emph{Domain Markov property.} For a subgraph $G'$ of $G$, given
$\omega|_{E(G)\setminus E(G')}$, the law of $\omega|_{E(G')}$ is $\phi^\zeta_{G',\pp}$, where
$\zeta$ is the partition of $\partial G'$ (the vertices of $G'$ incident to $E(G)\setminus E(G')$ or
in $\partial G$) induced by the outside configuration and $\xi$.
\item[\tagdef{FKG}{(FKG)}] \emph{FKG inequality.} Increasing events are positively correlated, and
conditional laws given the states of some edges are stochastically increasing in those states.
\item[\tagdef{CBC}{(CBC)}] \emph{Comparison of boundary conditions.} If $\xi\le\xi'$, then
$\phi^\xi_{G,\pp}\le_{\rm st}\phi^{\xi'}_{G,\pp}$; hence, for a subgraph $G'$,
\[
\phi^0_{G',\pp}\le_{\rm st}\phi^\xi_{G,\pp}|_{E(G')}\le_{\rm st}\phi^1_{G',\pp}.
\]
Moreover, $\phi^\xi_{G,\pp}$ is stochastically increasing in $\pp$.
\item[\tagdef{FE}{(FE)}] \emph{Finite-energy property.} Given the states of all other edges, the
edge $e$ is open with probability in $[\pp/(2-\pp),\pp]$.
\end{enumerate}

For every $\pp\in(0,1)$, the infinite-volume FK-Ising measure \[\phi_\pp:=\lim_{S\to\infty}\phi^1_{\Lambda_S,\pp}=\lim_{S\to\infty}\phi^0_{\Lambda_S,\pp}\]
is well defined (see, e.g., \cite[Thms.~4.19 and~6.17(b)]{Gri06} and \cite{DST17}).

Let $D\subset\Z^2$ be finite. Let $G_D$ be the graph of the edges with at least one endpoint in $D$,
with $\partial G_D:=V(G_D)\setminus D$, and let $\langle\cdot\rangle^+_{D,\beta}$ denote the
expectation with respect to the Ising measure on $G_D$ at inverse temperature $\beta$ with plus
boundary condition, that is, with the spins of the vertices of $\partial G_D$ equal to $+1$.
Throughout the paper, angle brackets $\langle\cdot\rangle$, with sub- and superscripts indicating
the domain, the inverse temperature and the boundary condition, denote the expectation with respect
to the corresponding Ising measure, and we write $\langle X;Y\rangle:=\langle XY\rangle-\langle
X\rangle\langle Y\rangle$ for the corresponding covariance.
\begin{lemma}[\cite{ES88} and {\cite[Proof of Lemma~2.1]{GK25}}]\label{lem:ES}
Let $\pp=\pp(\beta)$. For $e\ne f\in E(G_D)$,
\[
\langle\varepsilon_e;\varepsilon_f\rangle^+_{D,\beta}=\frac4{\pp^2}\Cov_{\phi^1_{G_D,\pp}}(\omega_e,\omega_f),\qquad
\langle\varepsilon_e\rangle^+_{D,\beta}=\frac2\pp\,\phi^1_{G_D,\pp}[\omega_e]-1 .
\]
The same holds in infinite volume with the plus state (the unique Gibbs state if $\beta\le\beta_c$)
and $\phi_\pp$.
\end{lemma}
By \tagref{FKG}{(FKG)}, all these covariances are nonnegative; in the spin language, this is
Griffiths' second inequality (GKS~II) \cite{Gri67}.

The last FK tool is a coupling of two FK measures along an adapted exploration, up to a stopping
time. It will be used in Sections~\ref{sec:cov} and~\ref{sec:VP} and in Appendix~\ref{app:strip},
where two boundary conditions are compared while edges are revealed one at a time.
\begin{lemma}[monotone exploration coupling; see, e.g., {\cite[Prop.~2.6, Remarks~2.7 and~2.8]{DCM}}]\label{lem:coupling}
Let $G$ be finite and $\xi\le\xi'$ be boundary conditions on $G$. Let $U_1,U_2,\dots$ be i.i.d.\
uniform on $[0,1]$, $\cG_k:=\sigma(U_1,\dots,U_k)$ and let $e_1,e_2,\dots$ be an \emph{adapted
exploration} of $E(G)$: distinct edges, with each $e_k$ measurable with respect to $\cG_{k-1}$ (for
instance a function of the states already revealed in $\omega$, in $\omega'$ or in both). Then there
is a coupling of $\omega\sim\phi^\xi_{G,\pp}$ and $\omega'\sim\phi^{\xi'}_{G,\pp}$, with
$\omega_{e_k},\omega'_{e_k}$ being $\cG_k$-measurable and $\omega\le\omega'$ on the revealed edges,
such that, for every stopping time $\tau$, conditionally on $\cG_\tau$, the unrevealed parts of
$\omega$ and $\omega'$ have the laws $\phi^\xi_{G,\pp}(\cdot\mid\omega_{e_1},\dots,\omega_{e_\tau})$
and $\phi^{\xi'}_{G,\pp}(\cdot\mid\omega'_{e_1},\dots,\omega'_{e_\tau})$; these are FK measures on
the unrevealed edges with the induced boundary conditions. After $\tau$, the unrevealed parts may be
coupled in any manner, for instance identically where these two laws coincide, or monotonically.
\end{lemma}
In this coupling, both configurations are revealed with the same uniform variables:
\[
\omega_{e_k}:=\one\{U_k\le a_k\}\qquad\text{and}\qquad\omega'_{e_k}:=\one\{U_k\le a'_k\},
\]
where $a_k\le a'_k$ are the conditional probabilities that $e_k$ is open given the revealed states,
by \tagref{FKG}{(FKG)} and \tagref{CBC}{(CBC)}.

\subsection{The interface and its convergence at criticality}\label{ss:interface} Recall from
\S\ref{ss:model} that $\gamma$ is the dual path from $x_1$ to $x_2$ that separates the $+$ spins, on
its left, from the $-$ spins, on its right, with a fixed turning rule. Drawn through the midpoints
of the crossed edges, $\gamma$ is a simple curve, and every inspected vertex lies within distance
$2$ of $\gamma$. For an admissible $g$, the set $\Insp(g)$ of inspected vertices and their spins are
determined by $g$, the spins being $+1$ on the left of $g$ and $-1$ on its right. Conversely,
$\{\gamma=g\}$ is the event that the vertices of $\Insp(g)$ carry these spins, since the exploration
reads no other spins in $\Omega_n$, and both endpoints of every edge crossed by $g$ lie in
$\Insp(g)\cup\partial\Omega_n$. Let $\cD(g)$ be the set of nearest-neighbor components of
$\Omega_n\setminus\Insp(g)$. Each $D\in\cD(g)$ lies on one side of $g$, so that the vertices of
$\Insp(g)\cup\partial\Omega_n$ adjacent to $D$ all carry the same spin $\epsilon_D$. Hence, given
$\gamma=g$, the spins in different components are independent and the law on $D$ is the Ising model
on $D$ with boundary condition given by $\epsilon_D$; moreover, the energies of the edges with no
endpoint in $\bigcup_DD$ are deterministic. By the global spin flip and GKS~II, all energy
covariances under this conditional law are nonnegative.

We view $\delta\gamma$ as an element of $\cX$; the map $g\mapsto\delta g$ is injective on admissible
paths. We use continuum results only through the following theorem. It is stated in \cite{CDHKS14}
for discrete Dobrushin domains that approximate a bounded simply connected domain with two marked
points. It applies here, since the Dobrushin squares converge to $(\Omega;x_1,x_2)$ in the
Carath\'eodory sense, with converging marked points, and since drawing $\gamma$ on the
square--octagon lattice, as in \cite{CDHKS14}, changes $\delta\gamma$ by $O(\delta)$ in the metric
$d$.

\begin{inthm}[{\cite[Thm.~1]{CDHKS14}}; see also {\cite{CS12,KS17}}]\label{thm:CDHKS}
As $\delta\to0$, the law $P^\delta_c$ of $\delta\gamma$ under $P_{\beta_c}$ converges weakly in
$(\cX,d)$ to the law $\nu$ of chordal $\SLE_3$ in $\Omega$ from $x_1$ to $x_2$.
\end{inthm}

\begin{remark}[turning rule]\label{rem:turning}
The paper \cite{CDHKS14} works with the leftmost or the rightmost interface, and our rule (``$+$
nearest-neighbor connected, $-$ $*$-connected'') produces one of them. Either convention can be
used: the reflection of the square in its vertical axis, composed with the global spin flip, maps
our interface to the one with the opposite convention, for the Dobrushin boundary condition if $n$
is even, and for the boundary condition with $x_1,x_2$ moved one lattice step to the right if $n$ is
odd. Our proofs apply verbatim to these boundary conditions, since the position of the sign changes
enters only through Lemma~\ref{lem:factor}, which uses only that $x_1^*$ and $x_2^*$
(\S\ref{ss:strip}) have the same abscissa, at distance at least $\lfloor n/4\rfloor$ from the
vertical sides of the dual box. Theorem~\ref{thm:CDHKS} applies to both boundary conditions, and
chordal $\SLE_3$ is invariant under the reflection, so that Theorem~\ref{thm:main} holds for either
convention.
\end{remark}

\subsection{RSW estimates and the mixing rate near criticality}\label{ss:inputs} Besides
Theorem~\ref{thm:CDHKS}, the proof uses RSW near criticality and the mixing rate $\Delta_\pp$ of
\cite{DCM} (see Definition~\ref{def:Delta} below), whose square will bound the edge covariances in
Section~\ref{sec:cov}. The inputs used only in the proof of Theorem~\ref{thm:R1} will be stated in
Appendix~\ref{app:R1}, and the strip identity in \S\ref{ss:DGP}. A rectangle is a set
\[
\mathcal R=([u,u+r]\times[v,v+r'])\cap\Z^2,
\]
with $u,v\in\Z$ (or in a translate of $\Z^2$, such as the dual lattice) and integer side lengths
$r,r'\ge1$, together with its edges. Let $\mathcal C_v(\mathcal R)$ be the event that an open path
in $\mathcal R$ joins its bottom and top sides. By planar duality, its complement is the event that
a closed dual path crosses $\mathcal R$ in the other direction, so the bounds below hold for dual
crossings as well.

\begin{inthm}[RSW near criticality; {\cite[Thm.~1.3]{DGP14}}, and {\cite[Thm.~1]{DHN11}} at $\ppc$]\label{thm:RSW}
For all $\lambda>0$ and $\varrho\ge1$, there are $0<c_-\le c_+<1$ such that, for every rectangle
$\mathcal R$ with side lengths $r$ and $r'\in[r/\varrho,\varrho r]$, every $\pp\in[\frac13,\frac34]$
with $|\pp-\ppc|\le\lambda/r$ and every boundary condition $\xi$ on $\mathcal R$,
\[
c_-\le\phi^\xi_{\mathcal R,\pp}[\mathcal C_v(\mathcal R)]\le c_+ .
\]
\end{inthm}
We write $c_\pm(\lambda,\varrho)$ for these constants. For $\pp=\pp(\beta)$ with $\beta\in W_M$, all
the rectangles that we use have sides at most $20n$, so that \eqref{eq:window} gives
$|\pp(\beta)-\ppc|\le40M/r$ for both sides, and we use $\lambda=40M$.

\begin{definition}\label{def:Delta}
For an edge $e$ incident to $0$ and $1\le r\le R$,
\[
\Delta_\pp(R):=\phi^1_{\Lambda_R,\pp}[\omega_e]-\phi^0_{\Lambda_R,\pp}[\omega_e]
\]
and
\[
\Delta_\pp(r,R):=\phi^1_{\Lambda_R,\pp}[\mathcal C(\Lambda_r)]-\phi^0_{\Lambda_R,\pp}[\mathcal C(\Lambda_r)];
\]
by symmetry, $\Delta_\pp(R)$ does not depend on the choice of $e$. The characteristic length is
defined to be
\begin{equation}\label{eq:Ldef}
L(\pp):=\inf\bigl\{R\ge1:\ \phi_\pp[\mathcal C(\Lambda_R)]\notin[\epsL,1-\epsL]\bigr\}\in\{1,2,\dots\}\cup\{\infty\},
\end{equation}
where $\mathcal C(\Lambda_R)$ is the event that an open path in $\Lambda_R$ joins its left and right
sides, and $\epsL$ is the threshold fixed below. We adopt the notation in \cite{DCM}.
\end{definition}
We fix the threshold $\epsL\in(0,\frac12)$ small enough for all the results of \cite{DCM} (see
\cite[\S1.3 and Remark~2.14]{DCM}), and smaller than $c_-(1,1)$ and $1-c_+(1,1)$; every result of
\cite{DCM} is used with this threshold. By Theorem~\ref{thm:RSW} with $\lambda=\varrho=1$ and
\tagref{DMP}{(DMP)}, $\phi_\pp[\mathcal C(\Lambda_R)]\in[\epsL,1-\epsL]$ whenever
$\pp\in[\frac13,\frac34]$ and $2R|\pp-\ppc|\le1$; in particular, $L(\ppc)=\infty$. By
\tagref{DMP}{(DMP)} and \tagref{CBC}{(CBC)}, $\Delta_\pp(\cdot)$ is nonincreasing; we put
$\Delta_\pp(0):=1$, and $\Delta_\pp(s):=\Delta_\pp(\lfloor s\rfloor)$ for real $s\ge0$.

The constants of \cite{DCM} are uniform only for $\pp$ away from $0$ and $1$ \cite[\S1.2]{DCM},
which is why the inputs below are stated for $\pp\in[\frac13,\frac34]$, an interval that contains
$\pp(\beta)$ for every $\beta\in W_M$ (\S\ref{ss:notation}) and every other parameter used below.
The same crossing bound gives, with $c_L:=\frac12$, the analogue for \eqref{eq:Ldef} of the lower
bound of \cite[Thm.~1.2]{DGP14} (see also \cite[Thm.~1.13, (1.30)]{DCM}):
\begin{equation}\label{eq:Lbound}
L(\pp)\ge\frac{c_L}{|\pp-\ppc|}\qquad\text{for every }\pp\in[\tfrac13,\tfrac34]\setminus\{\ppc\} .
\end{equation}
\begin{inthm}[ratio mixing; {\cite[Prop.~2.9, (2.4)]{DCM}}]\label{thm:mix}
There are $c_{\rm mix}>0$ and $K_0<\infty$ such that, for every $\pp\in[\frac13,\frac34]$, every
vertex $x$, every $r\ge1$ and $R\ge K_0r$, every finite $G\supseteq\Lambda_R(x)$, every event $A$
depending only on the edges of $\Lambda_r(x)$ (not necessarily monotone) and all boundary conditions
$\xi,\psi$ on $G$,
\[
|\phi^\psi_{G,\pp}[A]-\phi^\xi_{G,\pp}[A]|\le(r/R)^{c_{\rm mix}}\,\phi^\psi_{G,\pp}[A].
\]
\end{inthm}

\begin{definition}[four-arm event]\label{def:A4}
For $1\le r<R$, define $A_4(r,R)$ to be the event that four disjoint arms of alternating types
(open, dual-open, open, dual-open, as in \cite{DCM}) cross the annulus $\Lambda_R\setminus
\Lambda_r$.
\end{definition}
We put $\pi_4(\pp;r,R):=\phi_\pp[A_4(r,R)]$, we let $A_4(r,R)$ be the sure event for $r\ge R$ and we
write $A^x_4(r,R)$ for its translate by $x$. For $r\le r'\le R'\le R$, restricting the arms to a
smaller annulus gives
\begin{equation}\label{eq:A4mono}
A_4(r,R)\subseteq A_4(r',R)\cap A_4(r,R').
\end{equation}

We now record the forms in which RSW and ratio mixing will enter
Sections~\ref{sec:geom}--\ref{sec:VP}. We fix $\lambda>0$. There is $c_1=c_1(\lambda)\in(0,\frac12]$
such that the following hold for every $\pp\in[\frac13,\frac34]$.
\begin{enumerate}[leftmargin=3em]
\item[\tagdef{Ra}{(Ra)}] Let $r\ge1$ and $|\pp-\ppc|\le\lambda/r$, and let $\xi$ be a boundary
condition on a graph $G\supseteq\Lambda_{2r}(x)$. Under $\phi^\xi_{G,\pp}$, conditionally on any
configuration outside $E(\Lambda_{2r}(x))\setminus E(\Lambda_r(x))$, with probability at least
$c_1$, an open circuit with edges in $E(\Lambda_{2r}(x))\setminus E(\Lambda_r(x))$ surrounds
$\Lambda_r(x)$. The same holds for a closed dual circuit crossing only edges of this set.
\item[\tagdef{Rb}{(Rb)}] Let $Q$ be a translate of the box $\{0,\dots,S\}^2$, with inner boundary
$\partial Q$ and any boundary condition; let $y\in\partial Q$, $2\le r\le S/8$ and
$|\pp-\ppc|\le\lambda/r$; and let
\[
\mathrm{Ann}^+_r(y):=(\Lambda_{2r}(y)\setminus\Lambda_r(y))\cap Q
\]
be the \emph{half-annulus} of radii $r$ and $2r$ around $y$. Conditionally on any configuration
outside $E(\mathrm{Ann}^+_r(y))$, with probability at least $c_1$, a closed dual path crossing only
edges of $E(\mathrm{Ann}^+_r(y))$, with both ends on $\partial Q$ (a \emph{closed dual
half-circuit}), separates $\Lambda_r(y)\cap Q$ from $Q\setminus\Lambda_{2r}(y)$.

\end{enumerate}
These two properties are standard applications of Theorem~\ref{thm:RSW}, \tagref{DMP}{(DMP)},
\tagref{FKG}{(FKG)} and \tagref{FE}{(FE)}; we omit the proofs.

Next, we fix an integer $\Kmix\ge\max(K_0,4)$ with $\Kmix^{-c_{\rm mix}}\le\frac12$
(Theorem~\ref{thm:mix}). Then, for every $\pp\in[\frac13,\frac34]$:
\begin{enumerate}[leftmargin=3em]
\item[\tagdef{Mb}{(Mix)}] If $r\ge1$, $A\in\cF(\Lambda_r(x))$ and $\xi$ is a boundary condition on a
graph $G\supseteq\Lambda_{\Kmix r}(x)$, then
$\frac23\phi_\pp[A]\le\phi^\xi_{G,\pp}[A]\le\frac32\phi_\pp[A]$.
\end{enumerate}
This is a standard application of Theorem~\ref{thm:mix} and \tagref{DMP}{(DMP)}; we omit the proof.

The last consequence is the bound $\Delta_\pp(s)\lesssim1/s$ for all $\pp$ with
$|\pp-\ppc|\le\lambda/N$ and all $s\le N$, including the scales above $L(\pp)$, which will be used
in the dyadic sums of Sections~\ref{sec:cov} and~\ref{sec:VP}; the lower bound in part~(i) of the
following lemma will be used in Appendix~\ref{app:R1}.
\begin{lemma}[the mixing rate near criticality]\label{lem:Delta}
\begin{enumerate}[label=\textup{(\roman*)}]
\item There is a universal constant $c_\Delta>0$ with $c_\Delta/s\le\Delta_{\ppc}(s)\le
c_\Delta^{-1}/s$ for all integers $s\ge1$.
\item There is a universal constant $C_\Delta\ge1$ such that $\Delta_\pp(s)\le C_\Delta/s$ for every
$\pp\in[\frac13,\frac34]$ and every integer $1\le s\le L(\pp)$.
\item For every $\lambda<\infty$, put $C_\Delta(\lambda):=C_\Delta\max\{1,\lambda/c_L\}$. Then
$\Delta_\pp(s)\le C_\Delta(\lambda)/s$ whenever $N\ge1$, $\pp\in[\frac13,\frac34]$, $|\pp-\ppc|\le
\lambda/N$ and $s$ is an integer with $1\le s\le N$.
\end{enumerate}
\noindent (For real $s\ge1$, with the convention $\Delta_\pp(s)=\Delta_\pp(\lfloor s\rfloor)$ above,
\textup{(ii)}--\textup{(iii)} hold with $2C_\Delta$ and $2C_\Delta(\lambda)$, since $\lfloor
s\rfloor\ge s/2$.)
\end{lemma}
\begin{proof}
(i) This is \cite[Thm.~1.13, (1.31)]{DCM} at $\pp=\ppc$, since $L(\ppc)=\infty$.

(ii) By the stability of the mixing rate \cite[Thm.~1.6(iii), (1.16)]{DCM}, there is $C_{\rm
D}<\infty$ with $\Delta_\pp(s)\le C_{\rm D}\Delta_{\ppc}(s)$ for every $\pp\in[\frac13,\frac34]$ and
$1\le s\le L(\pp)$. With $C_\Delta:=\max\{1,C_{\rm D}/c_\Delta\}$, the upper bound in (i) gives
$\Delta_\pp(s)\le C_\Delta/s$.

(iii) If $s\le L(\pp)$, then (ii) applies, since $C_\Delta\le C_\Delta(\lambda)$. Otherwise,
$\pp\ne\ppc$ and, by \eqref{eq:Lbound}, $L(\pp)\ge c_L/|\pp-\ppc|\ge c_LN/\lambda$. Since
$\Delta_\pp(\cdot)$ is nonincreasing, (ii) at $L(\pp)$ and $s\le N$ give
\[
\Delta_\pp(s)\le\Delta_\pp(L(\pp))\le C_\Delta/L(\pp)\le C_\Delta\lambda/(c_LN)\le C_\Delta\lambda/(c_Ls).\qedhere
\]
\end{proof}

\section{Proof of Theorem~\ref{thm:main} modulo four estimates}\label{sec:renyi} In this section, we
prove Theorem~\ref{thm:main}, assuming the four estimates \tagref{E1}{(E1)}--\tagref{E4}{(E4)} of
\eqref{eq:inputs}, which will be proved in Sections~\ref{sec:geom}--\ref{sec:mumu}. We first reduce
Theorem~\ref{thm:main}(i) to these estimates, using only exact identities. In \S\ref{ss:renyiid}, we
will compare $R_\delta^{\,a}=(R_{\delta})^a$ with the Radon--Nikodym derivative at mass $am$. This
cancels the term involving the first derivative of $\log P_\beta(\gamma=g)$, which has no useful
pathwise bound, and leaves the integral of the variance decrease $\cV_\beta$ of \S\ref{ss:strategy}
against an explicit kernel of mass proportional to $(m\delta)^2$. In \S\ref{ss:Dob}, the
Fortuin--Kasteleyn (random-cluster) representation of the Dobrushin and plus partition functions,
$Z^\Dob_\beta=Z^+_\beta\,\phi^1_\beta(\Dab)$, splits $\cV_\beta$ into a plus-domain part and the
boundary disorder correlation. As we will show in \S\ref{ss:split}, the latter enters only through a
second difference, hence, by monotonicity, only through an upper bound on the ratio
$\phi^1_\beta(\Dab)/\phi^1_{\beta_c}(\Dab)$ for $\beta\le\beta_c$, with no bound on its derivatives
in $\beta$. In \S\ref{ss:reduction}, we will compare both variances with a common bulk reference;
this reduces the proof to the four estimates \tagref{E1}{(E1)}--\tagref{E4}{(E4)}; the last three
give the pathwise bound \eqref{eq:pathwise}, and \tagref{E1}{(E1)} controls its exponential moments.
Finally, assuming \tagref{E1}{(E1)}--\tagref{E4}{(E4)}, we will deduce the bounds at the lattice
level in \S\ref{ss:assembly}, and we will pass to the scaling limit in \S\ref{ss:limit}.

\subsection{The R\'enyi and Kullback--Leibler identities}\label{ss:renyiid} In this subsection, we
show that every moment of $R_\delta$, as well as $\KL(P^\delta_c\|P^\delta_m)$, can be expressed
exactly in terms of the variance $\cV_\beta$ removed by conditioning on the interface.

For an admissible path $g$, let
\[
f_g(\beta):=\log P_\beta(\gamma=g)=\log\sum_{\sigma:\gamma(\sigma)=g}e^{\beta H(\sigma)}-\log Z^\Dob_\beta.
\]
Formula \eqref{eq:H} defines $P_\beta$ for every real $\beta$, so $f_g$ is smooth on $\R$, and
differentiating the two log-partition functions gives
\begin{equation}\label{eq:score}
\begin{gathered}
f_g'(\beta)=\E_\beta[H\mid\gamma=g]-\E_\beta[H],\qquad f_g''(\beta)=-\cV_\beta(g),\\
\cV_\beta(g)=\Var^\Dob_\beta(H)-\Var_\beta(H\mid\gamma=g),
\end{gathered}
\end{equation}
where $\E^\Dob_\beta:=\E_\beta$ and $\Var^\Dob_\beta:=\Var_\beta$ denote the expectation and the
variance under $P_\beta$. For $v\in\R$, we put
\[
R^{(v)}_\delta(g):=P_{\beta_c-v\delta}(\gamma=g)/P_{\beta_c}(\gamma=g),
\]
so that $R_\delta=R^{(m)}_\delta$, and we write $\E_{(v)}$ for the expectation under
$P_{\beta_c-v\delta}$. Let $k_v(\beta):=|v|\delta-|\beta-\beta_c|$ if $\beta$ lies between
$\beta_c-v\delta$ and $\beta_c$, and $k_v(\beta):=0$ otherwise. This is the weight in Taylor's
formula with integral remainder for $f_g$ between $\beta_c$ and $\beta_c-v\delta$; it is
nonnegative, and $\int k_v=\frac12(v\delta)^2$. Taylor's formula, together with \eqref{eq:score},
therefore gives, for every $v\in\R$,
\begin{equation}\label{eq:taylor}
\log R^{(v)}_\delta(g)=-v\delta\,f_g'(\beta_c)-\int k_v(\beta)\,\cV_\beta(g)\,d\beta .
\end{equation}
For $a\in\R$, we define the kernel $K_a:=k_{am}-a\,k_m$, and we put $K_{\rm KL}:=k_m$.

Recall from \S\ref{ss:mainresult} that $\KL(P^\delta_c\|P^\delta_m)=-\E_c[\log R_\delta]$ and
$\KL(P^\delta_m\|P^\delta_c)=\E_c[R_\delta\log R_\delta]$.

\begin{proposition}[R\'enyi and Kullback--Leibler identities]\label{prop:renyi}
Let $a\in\R$, and let $g$ be an admissible path.
\begin{enumerate}[label=\textup{(\roman*)},leftmargin=2em]
\item The following pathwise and moment identities hold:
\begin{align}
a\log R_\delta(g)-\log R^{(am)}_\delta(g)&=\int K_a(\beta)\,\cV_\beta(g)\,d\beta,\label{eq:renyipath}\\
\E_c\bigl[R_\delta^{\,a}\bigr]&=\E_{(am)}\Bigl[\exp\int K_a(\beta)\,\cV_\beta(\gamma)\,d\beta\Bigr].\label{eq:renyimom}
\end{align}
\item The kernel $K_{\rm KL}$ is nonnegative, and so is $K_a$ if $a\notin[0,1]$. Their $L^1$ norms
are
\[\|K_{\rm KL}\|_1=\tfrac12(m\delta)^2\qquad\text{and}\qquad\|K_a\|_1=\tfrac{a(a-1)}2(m\delta)^2\quad\text{for }a\notin[0,1].\]
\item For every $p>1$,
\begin{align}
\KL(P^\delta_c\|P^\delta_m)&=\int K_{\rm KL}(\beta)\,\E_c\bigl[\cV_\beta(\gamma)\bigr]\,d\beta,\label{eq:KLid}\\
\KL(P^\delta_m\|P^\delta_c)&\le\frac1{p-1}\log\E_c\bigl[R_\delta^{\,p}\bigr].\label{eq:KLbound}
\end{align}
\end{enumerate}
\end{proposition}
\begin{proof}
Identity \eqref{eq:renyipath} follows from \eqref{eq:taylor}, and \eqref{eq:renyimom} follows from
\eqref{eq:renyipath}. The nonnegativity of $K_{\rm KL}$ and of $K_a$ for $a\notin[0,1]$, and the
formulas for their $L^1$ norms, follow from direct computations. For~\eqref{eq:KLid}, we take the
expectation $\E_c$ in \eqref{eq:taylor} at $v=m$. Note that the first derivative
$f'_\gamma(\beta_c)$ has zero expectation:
\[
\E_c[f'_\gamma(\beta_c)]=\E_c\bigl[\E_c[H\mid\gamma]\bigr]-\E_c[H]=0.
\]
Therefore,
\[\KL(P^\delta_c\|P^\delta_m)=\int k_m(\beta)\,\E_c[\cV_\beta(\gamma)]\,d\beta=\int K_{\rm KL}(\beta)\,\E_c[\cV_\beta(\gamma)]\,d\beta,\]
which gives~\eqref{eq:KLid}. Finally, the bound on $\KL(P^\delta_m\|P^\delta_c)$ follows from
Jensen's inequality:
\[
\log\E_c\bigl[R_\delta^{\,p}\bigr]=\log\E_{(m)}\bigl[R_\delta^{\,p-1}\bigr]\ge(p-1)\,\E_{(m)}[\log R_\delta]=(p-1)\KL(P^\delta_m\|P^\delta_c).\qedhere
\]
\end{proof}

\subsection{The Dobrushin partition function and the boundary disorder correlation}\label{ss:Dob}
The Dobrushin variance in $\cV_\beta$ cannot be compared directly with a bulk quantity, because the
boundary condition itself forces an interface. In this subsection, we separate its contribution
exactly. Lemma~\ref{lem:dob} shows that $Z^\Dob_\beta$ is $Z^+_\beta$ times the probability, under
the wired FK measure, that no open path joins the two boundary arcs. The logarithmically growing
specific heat is then carried by the plus model alone, and the boundary condition enters only
through $\ell(\beta)=\log\phi^1_\beta(\Dab)$.

Let $\phi^1_\beta:=\phi^1_{G,\pp(\beta)}$ be the FK-Ising measure on
$G:=(\Omega_n\cup\partial\Omega_n,E)$ with the wired boundary condition, and let $\Dab$ be the
decreasing event that no open path connects the $+$ boundary arc to the $-$ boundary arc. Let
$Z^+_\beta:=\sum_\sigma e^{\beta H(\sigma)}$ be the partition function of the Ising model on $G$ at
inverse temperature $\beta$ with the $+$ boundary condition, where $\sigma$ ranges over
$\{\pm1\}^{\Omega_n}$ and is extended by $+1$ on $\partial\Omega_n$, and let $P^+_\beta$,
$\E^+_\beta$ and $\Var^+_\beta$ denote the corresponding probability measure, expectation and
variance. In particular, $\E^+_\beta H=\partial_\beta\log Z^+_\beta$ and
$\Var^+_\beta(H)=\partial_\beta^2\log Z^+_\beta$.

\begin{lemma}\label{lem:dob}
For every $\beta>0$: \textup{(i)} $Z^\Dob_\beta=Z^+_\beta\,\phi^1_\beta(\Dab)$; \textup{(ii)} with
$\ell(\beta):=\log\phi^1_\beta(\Dab)$,
\[
\ell'(\beta)=\E^\Dob_\beta H-\E^+_\beta H\qquad\text{and}\qquad\Var^\Dob_\beta(H)=\Var^+_\beta(H)+\ell''(\beta);
\]
\textup{(iii)} $\ell$ is nonincreasing.
\end{lemma}
\begin{proof}
We start with (i), which is the Fortuin--Kasteleyn (random-cluster) representation of the partition
function \cite{FK72}, carried out for both boundary conditions. Let $\pp=\pp(\beta)$ and note that
\[
e^{\beta\sigma_x\sigma_y}=e^\beta[(1-\pp)+\pp\one\{\sigma_x=\sigma_y\}].
\]
For $\eta\in\{+,\Dob\}$, expanding the product over $E$ gives
\[Z^\eta_\beta=e^{\beta|E|}\sum_{\omega\in\{0,1\}^E}\pp^{o(\omega)}(1-\pp)^{c(\omega)}N_\eta(\omega),\]
where $N_\eta(\omega)$ is the number of $\sigma$ that are constant on the clusters of
$(\Omega_n\cup\partial\Omega_n,\omega)$ and agree with the boundary condition $\eta$ on the outer
boundary. Let us count these configurations. A cluster without outer boundary vertices has two
admissible signs. A cluster containing outer boundary vertices has one admissible sign if its outer
boundary vertices carry a single $\eta$-value, and none otherwise. Hence, $N_+=2^{k_f(\omega)}$ and
$N_\Dob=2^{k_f(\omega)}\one_{\Dab}(\omega)$, where $k_f$ is the number of clusters without outer
boundary vertices. Once the outer boundary is identified with a single vertex, the number of
clusters is $k_f+1$, so the weight $\pp^{o}(1-\pp)^{c}2^{k_f}$ is proportional to
$\phi^1_\beta(\omega)$. Dividing the two expansions gives (i).

To prove (ii), we take logarithms in (i) and differentiate in $\beta$, once for the first identity
and twice for the second. Finally, (iii) holds because $\Dab$ is decreasing, $\phi^1_{G,\pp}$ is
stochastically increasing in $\pp$ by \tagref{CBC}{(CBC)} and $\pp(\beta)$ is increasing in $\beta$.
\end{proof}
Recall from \S\ref{ss:strategy} that $\phi^1_\beta(\Dab)=Z^\Dob_\beta/Z^+_\beta$ is the boundary
disorder correlation. Note that $R_\delta(g)$ is not expected to be bounded pathwise: $\log
P_\beta(\gamma=g)$ contains the term $-2\beta|g^*|$, where $|g^*|$ is the number of edges crossed by
$g$, and this term contributes to $R_\delta(g)$ the factor $e^{2m\delta|g^*|}$, which is
heuristically $\exp(n^{3/8+o(1)})$. However, since $-2\beta|g^*|$ is affine in $\beta$, it
disappears from the second differences used below.
\subsection{Second differences and the decomposition of $\cV_\beta$}\label{ss:split} In this
subsection, we show that the boundary disorder correlation enters the identities of
\S\ref{ss:renyiid} only through second differences, and that its monotonicity reduces these to
differences of the values of $\ell$, with no bound on its derivatives.

For a function $u$ on $W_M$, differentiable at $\beta_c$, and $a\in\R$ with $|a|m\le M$, we put
\begin{equation}\label{eq:seconddiff}
\begin{split}
\Delta^2_au&:=u(\beta_c-am\delta)-a\,u(\beta_c-m\delta)+(a-1)\,u(\beta_c),\\
\Delta^2_{\rm KL}u&:=u(\beta_c-m\delta)-u(\beta_c)+m\delta\,u'(\beta_c),
\end{split}
\end{equation}
so that $\Delta^2_{\rm KL}u=\lim_{a\to0}\Delta^2_au/(-a)$.

\begin{lemma}[second differences]\label{lem:seconddiff}
Let $u\in C^2(W_M)$ and $|a|m\le M$. \textup{(i)} $\int K_a(\beta)\,u''(\beta)\,d\beta=\Delta^2_au$
and $\int K_{\rm KL}(\beta)\,u''(\beta)\,d\beta=\Delta^2_{\rm KL}u$. \textup{(ii)} If $u$ is
nonincreasing, then
\[
\Delta^2_au\le\begin{cases}u(\beta_c-am\delta)-u(\beta_c)&\text{for }a>1,\\|a|\,[u(\beta_c-m\delta)-u(\beta_c)]&\text{for }a<0,\end{cases}\qquad\text{and}\qquad\Delta^2_{\rm KL}u\le u(\beta_c-m\delta)-u(\beta_c).
\]
\end{lemma}
\begin{proof}
Part (i) is Taylor's formula with integral remainder, as in \eqref{eq:taylor}. Part (ii) follows
readily from the definitions of $\Delta^2_au$ and $\Delta^2_{\rm KL}u$ in~\eqref{eq:seconddiff}.
\end{proof}

Recall that $p>1$ and $s>0$ are the orders of the positive and negative moments in
Theorem~\ref{thm:main}(i), and that $M=\max(p,s,1)\,m$. We use $a=p$ for the positive moments,
$a=-s$ for the negative moments and $a={\rm KL}$, which refers to the kernel $K_{\rm KL}$, for the
Kullback--Leibler divergence.

\begin{proposition}[decomposition of $\cV_\beta$]\label{prop:split}
Let $a\in\{p,-s,{\rm KL}\}$. For every admissible $g$,
\[
\int K_a(\beta)\,\cV_\beta(g)\,d\beta=\int K_a(\beta)\bigl[\Var^+_\beta(H)-\Var_\beta(H\mid\gamma=g)\bigr]\,d\beta+\Delta^2_a\ell.
\]
If moreover
\begin{equation}\label{eq:Smu}
C_{\rm dis}:=\sup_{n\ge n_0(M)}\ \sup_{\beta\in[\beta_c-M\delta,\beta_c]}\frac{\phi^1_\beta(\Dab)}{\phi^1_{\beta_c}(\Dab)}<\infty,
\end{equation}
then $\Delta^2_a\ell\le\kappa_a\log C_{\rm dis}$, where $\kappa_p:=1$, $\kappa_{-s}:=s$ and
$\kappa_{\rm KL}:=1$.
\end{proposition}
\begin{proof}
First, we prove the identity. By \eqref{eq:score} and Lemma~\ref{lem:dob}(ii),
\[\cV_\beta(g)=\Var^+_\beta(H)+\ell''(\beta)-\Var_\beta(H\mid\gamma=g).\]
Moreover, $\ell$ is smooth, as the logarithm of a positive rational function of $\pp(\beta)$. We
integrate against $K_a$, and Lemma~\ref{lem:seconddiff}(i) turns the integral of $\ell''$ into
$\Delta^2_a\ell$.

Next, we prove the upper bound on $\Delta^2_a\ell$. By Lemma~\ref{lem:dob}(iii), $\ell$ is
nonincreasing, so Lemma~\ref{lem:seconddiff}(ii) applies. It bounds $\Delta^2_a\ell$ by
$\ell(\beta)-\ell(\beta_c)$ at $\beta=\beta_c-pm\delta$ for $a=p$, by $s[\ell(\beta)-\ell(\beta_c)]$
at $\beta=\beta_c-m\delta$ for $a=-s$ and by $\ell(\beta)-\ell(\beta_c)$ at $\beta=\beta_c-m\delta$
for $a={\rm KL}$. In each case, $\ell(\beta)-\ell(\beta_c)$ is the logarithm of the ratio in
\eqref{eq:Smu} at a point $\beta$ of $[\beta_c-M\delta,\beta_c]$, hence it is at most $\log C_{\rm
dis}$.
\end{proof}
Thus, the term $\Delta^2_a\ell$ requires only an upper bound on the ratio
$\phi^1_\beta(\Dab)/\phi^1_{\beta_c}(\Dab)$, on the high-temperature side
$\beta\in[\beta_c-M\delta,\beta_c]$. We will prove~\eqref{eq:Smu} in Section~\ref{sec:mumu}.

\subsection{The reduction to four estimates}\label{ss:reduction} The two variances in
Proposition~\ref{prop:split} both grow like $n^2\log n$. In this subsection, we therefore compare
each of them with a common reference, up to additive errors $O(n^2)$ and $O(n^2\Pi_\delta(g))$,
respectively, and we state the four estimates to which Theorem~\ref{thm:main}(i) reduces. The
reference is local: each edge contributes its infinite-volume local specific heat up to a third of
its distance to the boundary. For $\beta\in W_M$, let
\[\chi_\beta(e;r):=\sum_{f:\,|f-e|\le r}\langle\varepsilon_e;\varepsilon_f\rangle_{\Z^2,\beta},\qquad X_\beta:=\sum_{e\in E}\chi_\beta\bigl(e;r_\Omega(e)\bigr),\]
where $r_\Omega(e)=1\vee\frac13d_\Omega(e)$, with $d_\Omega(e)$ the Euclidean distance from the
midpoint of $e$ to $\Z^2\setminus\Omega_n$ (see \S\ref{ss:notation}), and
$\langle\cdot\rangle_{\Z^2,\beta}$ is the expectation with respect to the infinite-volume Ising
measure, that is, the plus state if $\beta>\beta_c$ and the unique Gibbs state if $\beta\le\beta_c$.
For every admissible $g$ and $a\in\{p,-s,{\rm KL}\}$, adding and subtracting $X_\beta$ in the
identity of Proposition~\ref{prop:split} gives the exact decomposition \eqref{eq:reduction} of
\S\ref{ss:strategy}:
\[
\begin{split}
\int K_a(\beta)\,\cV_\beta(g)\,d\beta&=\int K_a(\beta)\bigl[\Var^+_\beta(H)-X_\beta\bigr]\,d\beta+\Delta^2_a\ell\\
&\quad+\int K_a(\beta)\bigl[X_\beta-\Var_\beta(H\mid\gamma=g)\bigr]\,d\beta .
\end{split}
\]
The proof of Theorem~\ref{thm:main}(i) uses the following precise form of the four estimates
\eqref{eq:inputs}, uniform in $n\ge n_0(M)$, $\beta\in W_M$ and admissible $g$, where $\Pi_\delta$
is the proximity functional of \S\ref{ss:strategy} (see also \eqref{eq:Lambda}):
\[
\begin{aligned}
&\text{\tagref{E1}{(E1)}}\quad \sup_{\beta\in W_M}\E_\beta\bigl[e^{\eta\Pi_\delta(\gamma)}\bigr]\le C(M,\eta),\quad\text{ for every }\eta>0 \quad&&\text{(Proposition~\ref{prop:expmom});}\\
&\text{\tagref{E2}{(E2)}}\quad X_\beta-\Var_\beta(H\mid\gamma=g)\le C_{\Pi}\,n^2\,\Pi_\delta(g) &&\text{(Proposition~\ref{prop:Bminus});}\\
&\text{\tagref{E3}{(E3)}}\quad \Var^+_\beta(H)-X_\beta\le C_Vn^2 &&\text{(Theorem~\ref{thm:VP});}\\
&\text{\tagref{E4}{(E4)}}\quad C_{\rm dis}<\infty\text{, that is, \eqref{eq:Smu}} &&\text{(Theorem~\ref{thm:mumu}).}
\end{aligned}
\]
The plus domain is compared with $X_\beta$ from above, by \tagref{E3}{(E3)}, and the slit domain
from below, by \tagref{E2}{(E2)}; in this way, the bulk logarithms cancel edge by edge. We insert
\tagref{E2}{(E2)}--\tagref{E4}{(E4)} and Proposition~\ref{prop:split} in \eqref{eq:reduction}, and
we use $\|K_a\|_1n^2=\frac{a(a-1)}2m^2$ for $a\in\{p,-s\}$ and $\|K_{\rm KL}\|_1n^2=\frac{m^2}2$.
This gives the following precise form of the pathwise bound of \S\ref{ss:strategy}: for every
admissible $g$ and $a\in\{p,-s,{\rm KL}\}$, with $\kappa_a$ as in Proposition~\ref{prop:split},
\begin{equation}\label{eq:pathwise}
\begin{gathered}
\int K_a(\beta)\,\cV_\beta(g)\,d\beta\le\|K_a\|_1\,n^2\bigl(C_V+C_{\Pi}\Pi_\delta(g)\bigr)+\kappa_a\log C_{\rm dis},\\
\|K_a\|_1\,n^2\in\bigl\{\tfrac{p(p-1)}2m^2,\tfrac{s(s+1)}2m^2,\tfrac{m^2}2\bigr\}.
\end{gathered}
\end{equation}
Together with Proposition~\ref{prop:renyi} and \tagref{E1}{(E1)}, this gives Theorem~\ref{thm:main},
as we will show in \S\ref{ss:assembly} and \S\ref{ss:limit}. The proofs of
\tagref{E1}{(E1)}--\tagref{E4}{(E4)} will occupy Sections~\ref{sec:geom}--\ref{sec:mumu}.

\subsection{The discrete bounds, assuming the four estimates}\label{ss:assembly} In this subsection,
we prove the bounds at the lattice level for $m>0$, assuming \tagref{E1}{(E1)}--\tagref{E4}{(E4)},
which will be proved in Sections~\ref{sec:geom}--\ref{sec:mumu}; for $m<0$, see the discussion at
the end of \S\ref{ss:strategy}. We fix $m>0$, $p>1$, $s>0$ and $M=\max(p,s,1)m$, and we let $n\ge
n_0(M)$, since the bounds for the finitely many smaller $n$ are immediate (\S\ref{ss:notation}).
Then \eqref{eq:pathwise} holds.

\begin{proof}[Proof of Theorem~\ref{thm:main}\textup{(i)}, assuming \textup{\tagref{E1}{(E1)}--\tagref{E4}{(E4)}}]
The idea is to write each moment of $R_\delta$ as an exponential moment under a near-critical
measure, by the R\'enyi identity, and to bound the exponent by \eqref{eq:pathwise}; what remains is
an exponential moment of $\Pi_\delta$, which is controlled by \tagref{E1}{(E1)}.

Let $a\in\{p,-s\}$. By Proposition~\ref{prop:renyi} and \eqref{eq:pathwise},
\[
\E_c[R_\delta^{\,a}]=\E_{(am)}\Bigl[\exp\int K_a(\beta)\,\cV_\beta(\gamma)\,d\beta\Bigr]\le C_{\rm dis}^{\kappa_a}\,e^{C_V\frac{a(a-1)}2m^2}\,
\E_{(am)}\Bigl[e^{C_{\Pi}\frac{a(a-1)}2m^2\,\Pi_\delta(\gamma)}\Bigr],
\]
where the equality is \eqref{eq:renyimom}, and the inequality follows by exponentiating
\eqref{eq:pathwise}, in which $\|K_a\|_1n^2=\frac{a(a-1)}2m^2$. Since $\beta_c-am\delta\in W_M$, the
right-hand side is bounded uniformly in $n$ by \tagref{E1}{(E1)} with
$\eta=C_{\Pi}\frac{a(a-1)}2m^2$. This proves the bounds on the positive and negative moments in
Theorem~\ref{thm:main}(i).

We now turn to the Kullback--Leibler divergences. By \eqref{eq:KLid} and \eqref{eq:pathwise} with
$a={\rm KL}$, we have
\[\KL(P^\delta_c\|P^\delta_m)\le\frac{m^2}2\bigl(C_V+C_{\Pi}\E_c[\Pi_\delta(\gamma)]\bigr)+\log C_{\rm dis} ,\]
and, by Jensen's inequality and \tagref{E1}{(E1)}, $\E_c[\Pi_\delta]\le\log\E_c[e^{\Pi_\delta}]\le
C$. Finally, $\KL(P^\delta_m\|P^\delta_c)$ is bounded by a moment already controlled, since
\[
\KL(P^\delta_m\|P^\delta_c)\le(p-1)^{-1}\log\E_c[R_\delta^p]
\]
by \eqref{eq:KLbound}.
\end{proof}

\subsection{Passage to the scaling limit}\label{ss:limit} It remains to pass to the limit, still
assuming \tagref{E1}{(E1)}--\tagref{E4}{(E4)} through the discrete bounds of \S\ref{ss:assembly}. By
Theorem~\ref{thm:CDHKS}, $P^\delta_c\to\nu$ weakly in $(\cX,d)$, and since $(\cX,d)$ is Polish,
$\{P^\delta_c\}$ is tight (\S\ref{ss:model}). For $q\in(1,\infty)$, let $q'=q/(q-1)$. For each $n$,
$P^\delta_m$ and $P^\delta_c$ assign positive probability to exactly the same finite set, namely the
images of the admissible paths under $g\mapsto\delta g$. Hence, $\E_{(m)}[F]=\E_c[R_\delta F]$ and
$\E_c[F]=\E_{(m)}[R_\delta^{-1}F]$ for every function $F$ of the curve. We use the following
standard fact.

\begin{lemma}[$L^q$ duality]\label{lem:limitLp}
Let $\mu$ and $\nu$ be Borel probability measures on a metric space, and let $C<\infty$ and $q>1$.
If $\mu(\varphi)\le C\,\nu(\varphi^{q'})^{1/q'}$ for every bounded continuous $\varphi\ge0$, then
$\mu\ll\nu$ and $\|d\mu/d\nu\|_{L^q(\nu)}\le C$.
\end{lemma}
\begin{proof}
First, we extend the hypothesis to bounded Borel functions $\varphi\ge0$. Bounded continuous
functions are dense in $L^1(\mu+\nu)$, so we may approximate such a $\varphi$ in $L^1(\mu+\nu)$ by
continuous $\varphi_j$ with $0\le\varphi_j\le\sup\varphi$. Then $\mu(\varphi_j)\to\mu(\varphi)$, and
also $\nu(\varphi_j^{q'})\to\nu(\varphi^{q'})$, since
\[
|x^{q'}-y^{q'}|\le q'(\sup\varphi)^{q'-1}|x-y|
\]
on $[0,\sup\varphi]$. With $\varphi=\one_B$, the extended hypothesis reads $\mu(B)\le
C\nu(B)^{1/q'}$, which gives $\mu\ll\nu$. Now let $f:=d\mu/d\nu$, let $T>0$ and take
$\varphi=(f\wedge T)^{q-1}$. Since $(q-1)q'=q$, we obtain
\[
\int(f\wedge T)^q\,d\nu\le\int f\,(f\wedge T)^{q-1}\,d\nu=\mu(\varphi)\le C\Bigl(\int(f\wedge T)^q\,d\nu\Bigr)^{1/q'}.
\]
The integral is finite, so $\|f\wedge T\|_{L^q(\nu)}\le C$, and we let $T\to\infty$ by monotone
convergence.
\end{proof}

\begin{proof}[Proof of Theorem~\ref{thm:main}\textup{(ii)}]
The uniform bounds on the positive and negative moments in part~(i) pass to the limit through
H\"older's inequality and Lemma~\ref{lem:limitLp}. Let $C_q:=\sup_n\E_c[R_\delta^q]$ and
$C_{-s}:=\sup_n\E_c[R_\delta^{-s}]$, which are finite by part (i).

\emph{Tightness.} For compact $\mathcal Z\subseteq\cX$, H\"older's inequality gives
\[P^\delta_m(\mathcal Z^c)=\E_c[R_\delta\one_{\mathcal Z^c}]\le C_q^{1/q}P^\delta_c(\mathcal Z^c)^{1/q'},\]
and $\{P^\delta_c\}$ is tight. Hence, $\{P^\delta_m\}$ is tight.

\emph{Absolute continuity: $\mu\ll\nu$.} Let $\delta_k\to0$ be such that $P^{\delta_k}_m\to\mu$
weakly. For bounded continuous $\varphi\ge0$, we have
\[\mu(\varphi)=\lim_k\E_c[R_{\delta_k}\varphi(\delta_k\gamma)]\le C_q^{1/q}\limsup_k\E_c[\varphi(\delta_k\gamma)^{q'}]^{1/q'}=C_q^{1/q}\,\nu(\varphi^{q'})^{1/q'},\]
where the equality uses the weak convergence and $\E_{(m)}[F]=\E_c[R_\delta F]$, the inequality is
H\"older's inequality and the last step uses the weak convergence $P^\delta_c\to\nu$, applied to the
bounded continuous function $\varphi^{q'}$. Hence, by Lemma~\ref{lem:limitLp}, $\mu\ll\nu$ and
$\|d\mu/d\nu\|_{L^q(\nu)}\le C_q^{1/q}$, for every $q\in(1,\infty)$.

\emph{Absolute continuity: $\nu\ll\mu$.} Let $q:=1+s$. Since
\[
\E_{(m)}[R_\delta^{-q}]=\E_c[R_\delta^{1-q}]=\E_c[R_\delta^{-s}],
\]
H\"older's inequality under $P^{\delta_k}_m$ gives, for bounded continuous $\varphi\ge0$,
\[
\nu(\varphi)=\lim_k\E_{(m)}[R_{\delta_k}^{-1}\varphi(\delta_k\gamma)]\le(C_{-s})^{1/q}\,\mu(\varphi^{q'})^{1/q'},
\]
where we used $\E_c[F]=\E_{(m)}[R_\delta^{-1}F]$ and the weak convergence of both sequences. Hence,
by Lemma~\ref{lem:limitLp} with the roles of $\mu$ and $\nu$ exchanged, $\nu\ll\mu$ and
$\|d\nu/d\mu\|_{L^{1+s}(\mu)}\le(C_{-s})^{1/(1+s)}$, for every $s>0$.

\emph{Kullback--Leibler divergences.} With $f=d\mu/d\nu$, we have
\[
\KL(\mu\|\nu)=\E_\nu[f\log f]\le\frac1{q-1}\E_\nu[f^q]<\infty,
\]
since $x\log x\le x^q/(q-1)$ for $x>0$. In the same way, $\KL(\nu\|\mu)<\infty$.
\end{proof}

\section{Exponential moments of the proximity functional: Proof of~(E1)}\label{sec:geom} In this
section, we prove the estimate \tagref{E1}{(E1)} (Proposition~\ref{prop:expmom}). Recall from
\S\ref{ss:strategy} that, for an admissible path $g$, the proximity functional is
\begin{equation}\label{eq:Lambda}
\Pi_\delta(g)=1+\Phi(g),\qquad\Phi(g)=\delta^2\sum_{e\in E}\log_+\frac{d_\Omega(e)}{1\vee d_g(e)} .
\end{equation}
Here, $d_\Omega(e)$ and $d_g(e)$ denote the Euclidean distances from $\mathrm{mid}(e)$ to
$\Z^2\setminus\Omega_n$ and to $g$, respectively (see \S\ref{ss:notation}).

\begin{proposition}\label{prop:expmom}
For every $\eta>0$, there is $C(M,\eta)<\infty$ such that
\[
\sup_{n\ge n_0(M)}\ \sup_{\beta\in W_M}\ \E_\beta\bigl[e^{\eta\Pi_\delta(\gamma)}\bigr]\le C(M,\eta).
\]
\end{proposition}

Throughout this section, $\beta\in W_M$ and $n\ge n_0(M)$. The next lemma bounds the probability
that the interface comes close to several edges at once.
\begin{lemma}[multipoint proximity bound]\label{lem:multipoint}
There are $C_{\rm mc}=C_{\rm mc}(M)$ and $\varepsilon_{\rm mc}=\varepsilon_{\rm
mc}(M)\in(0,\frac12]$ such that, for all $\beta\in W_M$, all distinct $e_1,\dots,e_j\in E$ and all
$r_1,\dots,r_j\ge1$,
\[
P_\beta\bigl(d_\gamma(e_i)\le r_i,\  1\leq i\leq j\bigr)\le\prod_{i=1}^jC_{\rm mc}\Bigl(\frac{r_i}{s_i}\Bigr)^{\varepsilon_{\rm mc}},\qquad
s_i:=\min\Bigl(d_\Omega(e_i),\ \min_{\substack{1\le l\le j\\ l\ne i}}|e_i-e_l|\Bigr)
\]
\textup{(}for $j=1$, the inner minimum is over the empty set, and $s_1=d_\Omega(e_1)$\textup{)}.
\end{lemma}
Lemma~\ref{lem:multipoint} is a standard application of the RSW estimates of Theorem~\ref{thm:RSW}
to monochromatic circuits around the $e_i$, which the interface cannot cross; we omit the proof. The
only input on the spins is the following: given the spins outside an annulus
$\Lambda_{8t}(x)\setminus\Lambda_t(x)$ contained in $\Omega_n$, a nearest-neighbor circuit of plus
spins in this annulus surrounds $\Lambda_t(x)$ with probability at least $c=c(M)>0$. Indeed, by the
monotonicity of the Ising model in the boundary spins, it suffices to take minus boundary spins, and
then, in the Edwards--Sokal coupling \cite{ES88}, \tagref{Ra}{(Ra)} gives an open circuit in the
middle dyadic annulus that is separated from the boundary of the annulus by closed dual circuits in
the two others, and the cluster of this open circuit is plus with probability $\frac12$. The
factorized bound in Lemma~\ref{lem:multipoint} controls the binomial moments of each dyadic scale of
$\Phi$ (Lemma~\ref{lem:factorial}), and the exponential moments follow.

We now turn to the proof of Proposition~\ref{prop:expmom}. Since
\[
\log_+y\le\log2\cdot\#\{k\ge1:2^{k-1}\le y\}
\]
for $y>0$, we have
\begin{equation}\label{eq:dyadicPhi}
\Phi(\gamma)\le\log2\sum_{k\ge1}\Phi_k,\qquad \Phi_k:=\delta^2\mathcal N_k,\qquad \mathcal N_k:=\#\{e\in E:\ 1\vee d_\gamma(e)\le2^{1-k}d_\Omega(e)\}.
\end{equation}
Note that $\Phi_k\le\delta^2|E|\le3$, and that $\Phi_k=0$ for $2^{k-1}>n$, since $d_\Omega\le n$.
The key estimate is the following bound on the binomial moments of $\mathcal N_k$. Here and below,
we use the convention that $\binom{\mathcal N_k}{j}=0$ if $j>\mathcal N_k$.

\begin{lemma}[binomial moments of one scale]\label{lem:factorial}
There is $C_{\mathcal N}=C_{\mathcal N}(M)\ge6C_{\rm mc}$ such that, with $\mathcal N_k$ the random
variable of \eqref{eq:dyadicPhi} built from the interface $\gamma$, we have
$\E_\beta\bigl[\binom{\mathcal N_k}{j}\bigr]\le\bigl(C_{\mathcal N}2^{-k\varepsilon_{\rm
mc}}n^2\bigr)^j$ for all $n\ge n_0(M)$, $\beta\in W_M$, $k\ge1$ and $j\ge1$, where $C_{\rm
mc}=C_{\rm mc}(M)$ and $\varepsilon_{\rm mc}\in (0,\frac{1}{2}]$ are the constants of
Lemma~\ref{lem:multipoint}.
\end{lemma}

We first prove Proposition~\ref{prop:expmom} modulo Lemma~\ref{lem:factorial}, whose proof we
postpone to the end of this section.

\begin{proof}[Proof of Proposition~\ref{prop:expmom}]
We fix $\eta>0$. Since $\Pi_\delta=1+\Phi$, it follows from \eqref{eq:dyadicPhi} that
\begin{equation}\label{eq:expscales}
\E_\beta\bigl[e^{\eta\Pi_\delta(\gamma)}\bigr]\le e^{\eta}\,\E_\beta\Bigl[\exp\Bigl(\eta\log2\sum_{k\ge1}\Phi_k\Bigr)\Bigr],
\end{equation}
and it suffices to bound the right-hand side uniformly in $\beta\in W_M$ and $n\ge n_0(M)$.

\emph{Step 1: one scale.} We first show that, for all $k\ge1$ and $\eta'>0$ with $\eta'\delta^2\le1$
and $4\eta'C_{\mathcal N}2^{-k\varepsilon_{\rm mc}}\le1$, we have
\begin{equation}\label{eq:onescale}
\E_\beta\bigl[e^{\eta'\Phi_k}\bigr]\le2.
\end{equation}
Let $u>0$ be defined by $1+u=e^{\eta'\delta^2}$. Since $e^y-1\le2y$ for $y\in[0,1]$, and since
$\eta'\delta^2\le1$, we have $u\le2\eta'\delta^2$. Using $(1+u)^{\mathcal
N_k}=\sum_{j\ge0}\binom{\mathcal N_k}ju^j$, we have
\[
\begin{aligned}
\E_\beta\bigl[e^{\eta'\Phi_k}\bigr]&=\E_\beta\bigl[(1+u)^{\mathcal N_k}\bigr]=\sum_{j\ge0}u^j\,\E_\beta\Bigl[\binom{\mathcal N_k}{j}\Bigr]
\le\sum_{j\ge0}\bigl(2\eta'\delta^2\,C_{\mathcal N}2^{-k\varepsilon_{\rm mc}}n^2\bigr)^j\\
&=\sum_{j\ge0}\bigl(2\eta'C_{\mathcal N}2^{-k\varepsilon_{\rm mc}}\bigr)^j\le\sum_{j\ge0}2^{-j}=2,
\end{aligned}
\]
where the inequality in the first line follows from Lemma~\ref{lem:factorial} and
$u\le2\eta'\delta^2$, the equality in the second line from $\delta n=1$ and the last inequality from
$4\eta'C_{\mathcal N}2^{-k\varepsilon_{\rm mc}}\le1$.

\emph{Step 2: the weights.} Let $w_*:=1-2^{-\varepsilon_{\rm mc}/2}$, and let $k_*$ be the least
integer $k_*\ge0$ with $w_*\,2^{(k_*+1)\varepsilon_{\rm mc}}\ge4\eta C_{\mathcal N}\log2$. For
$k>k_*$, we define
\[
w_k:=w_*\,2^{-(k-k_*-1)\varepsilon_{\rm mc}/2},\qquad\eta'_k:=\frac{\eta\log2}{w_k},
\]
so that $w_k\eta'_k=\eta\log2$. Note that $w_*$ and $k_*$ depend only on $M$ and $\eta$. We let
$n_\eta:=\eta\log2/w_*$, and we claim that the following hold:
\begin{enumerate}[label=\textup{(\roman*)}]
\item $\sum_{k>k_*}w_k=1$;
\item $4\eta'_kC_{\mathcal N}2^{-k\varepsilon_{\rm mc}}\le1$ for all $k>k_*$;
\item if $n\ge n_\eta$, then $\eta'_k\delta^2\le1$ for all $k>k_*$ with $2^{k-1}\le n$.
\end{enumerate}
Property (i) follows from the geometric series. For (ii), let $k>k_*$. Then
\[
4\eta'_kC_{\mathcal N}2^{-k\varepsilon_{\rm mc}}=\frac{4\eta C_{\mathcal N}\log2}{w_*}\,2^{-(k+k_*+1)\varepsilon_{\rm mc}/2}
\le\frac{4\eta C_{\mathcal N}\log2}{w_*}\,2^{-(k_*+1)\varepsilon_{\rm mc}}\le1,
\]
where the first inequality holds because $k\ge k_*+1$, and the second by the choice of $k_*$. For
(iii), let $k>k_*$ with $2^{k-1}\le n$. Since $k_*\ge0$, we have
\[
w_k\ge w_*2^{-(k-1)\varepsilon_{\rm mc}/2}\ge w_*n^{-\varepsilon_{\rm mc}/2},
\]
and therefore
\[
\eta'_k\delta^2\le\frac{\eta\log2}{w_*}\,n^{\varepsilon_{\rm mc}/2-2}\le\frac{n_\eta}{n}\le1,
\]
where the second inequality holds because $\varepsilon_{\rm mc}\le\frac12$ and $n\ge1$, and the last
because $n\ge n_\eta$.

\emph{Step 3: the case $n\ge n_\eta$.} The scales $k\le k_*$ are bounded deterministically: since
$\Phi_k\le3$ for every $k$, we have $\eta\log2\sum_{k\le k_*}\Phi_k\le3\eta k_*\log2$. For the
scales $k>k_*$, we have
\[
\E_\beta\Bigl[\exp\Bigl(\eta\log2\sum_{k>k_*}\Phi_k\Bigr)\Bigr]=\E_\beta\Bigl[\exp\Bigl(\sum_{k>k_*}w_k\eta'_k\Phi_k\Bigr)\Bigr]
\le\sum_{k>k_*}w_k\,\E_\beta\bigl[e^{\eta'_k\Phi_k}\bigr]\le2.
\]
Here, the equality follows from $w_k\eta'_k=\eta\log2$, and the first inequality from (i) and the
convexity of the exponential function. The second inequality follows from (i), together with the
bound $\E_\beta[e^{\eta'_k\Phi_k}]\le2$ for every $k>k_*$: if $2^{k-1}\le n$, this bound follows
from \eqref{eq:onescale}, whose assumptions hold by (ii) and (iii); if $2^{k-1}>n$, then $\Phi_k=0$,
and the bound is trivial. Together with the deterministic bound on the scales $k\le k_*$ and with
\eqref{eq:expscales}, this gives
\[
\E_\beta[e^{\eta\Pi_\delta(\gamma)}]\le2\exp(\eta+3\eta k_*\log2).
\]

\emph{Step 4: the case $n<n_\eta$.} It remains to consider the finitely many $n$ with $n_0(M)\le
n<n_\eta$. Since $\Phi_k\le3$ for every $k$, and $\Phi_k=0$ for $2^{k-1}>n$, it follows from
\eqref{eq:dyadicPhi} that $\Phi(\gamma)\le3\log2\,(1+\log_2n)$ deterministically. Hence,
\[
\E_\beta[e^{\eta\Pi_\delta(\gamma)}]\le\exp\bigl(\eta+3\eta\log2\,(1+\log_2n_\eta)\bigr)
\]
for these $n$.

Since $k_*$ and $n_\eta$ depend only on $M$ and $\eta$, the bounds of Steps~3 and~4 are uniform in
$\beta\in W_M$ and $n\ge n_0(M)$, and the proof is complete.
\end{proof}

It remains to prove Lemma~\ref{lem:factorial}. The bound of Lemma~\ref{lem:multipoint} is
unfavorable when the edges $e_i$ are close to one another, and the number of such configurations is
controlled by summing out the edges one at a time, following their nearest neighbors.

\begin{proof}[Proof of Lemma~\ref{lem:factorial}]
We fix $k$ and $j\ge1$, and we put $r(e):=2^{1-k}d_\Omega(e)$. An edge counted by $\mathcal N_k$
satisfies $r(e)\ge1$ and $d_\gamma(e)\le r(e)$. Each $j$-element set of counted edges corresponds to
$j!$ ordered $j$-tuples of distinct counted edges, so that $j!\,\binom{\mathcal N_k}{j}$ is the
number of such tuples, and only tuples with $r(e_i)\ge1$ for all $i$ contribute.
Lemma~\ref{lem:multipoint} applies to these tuples with $r_i=r(e_i)$, and, all terms being
nonnegative, we obtain
\begin{equation}\label{eq:factstart}
\begin{split}
j!\,\E_\beta\Bigl[\binom{\mathcal N_k}{j}\Bigr]&\le\sum_{\substack{e_1,\dots,e_j\ \rm distinct\\ r(e_i)\ge1,\  1\leq i\leq j}}P_\beta\bigl(d_\gamma(e_i)\le r(e_i),\ 1\le i\le j\bigr)\\
&\le\sum_{e_1,\dots,e_j\ \rm distinct}\ \prod_{i=1}^jC_{\rm mc}\Bigl(\frac{r(e_i)}{s_i}\Bigr)^{\varepsilon_{\rm mc}}.
\end{split}
\end{equation}
For $j=1$, we have $s_1=d_\Omega(e_1)$, so that $r(e_1)/s_1=2^{1-k}$, and, using $|E|\le3n^2$ and
$2^{\varepsilon_{\rm mc}}\le2$, we obtain
\[
\E_\beta[\mathcal N_k]\le|E|\,C_{\rm mc}2^{(1-k)\varepsilon_{\rm mc}}\le6C_{\rm mc}2^{-k\varepsilon_{\rm mc}}n^2.
\]
Now let $j\ge2$. We proceed in three steps.

\emph{Step 1: nearest neighbors.} For each $i$, let $\mathrm{nb}(i)\ne i$ minimize $|e_i-e_l|$ (if
there are multiple such $l$, let $\mathrm{nb}(i)$ be the smallest one), so that
$s_i=\min(d_\Omega(e_i),|e_i-e_{\mathrm{nb}(i)}|)$. Since $\max(1,a)^{\varepsilon_{\rm
mc}}\le1+a^{\varepsilon_{\rm mc}}$ and $d_\Omega\le n$, we have
\[
\Bigl(\frac{r(e_i)}{s_i}\Bigr)^{\varepsilon_{\rm mc}}=2^{(1-k)\varepsilon_{\rm mc}}\max\Bigl(1,\frac{d_\Omega(e_i)}{|e_i-e_{\mathrm{nb}(i)}|}\Bigr)^{\varepsilon_{\rm mc}}
\le2^{(1-k)\varepsilon_{\rm mc}}\bigl(1+n^{\varepsilon_{\rm mc}}|e_i-e_{\mathrm{nb}(i)}|^{-\varepsilon_{\rm mc}}\bigr).
\]
Inserting this into \eqref{eq:factstart} and expanding the product over $i$, we obtain
\[
j!\,\E_\beta\Bigl[\binom{\mathcal N_k}{j}\Bigr]\le\bigl(C_{\rm mc}2^{(1-k)\varepsilon_{\rm mc}}\bigr)^j\sum_{S\subseteq\{1,\dots,j\}}\ \sum_{e_1,\dots,e_j\ \rm distinct}\ \prod_{i\in S}n^{\varepsilon_{\rm mc}}|e_i-e_{\mathrm{nb}(i)}|^{-\varepsilon_{\rm mc}}.
\]

\emph{Step 2: two bounds on sums of $\mathcal K$.} For distinct edges $e,e'\in E$, we put
\[
\mathcal K(e,e'):=n^{\varepsilon_{\rm mc}}|e-e'|^{-\varepsilon_{\rm mc}},
\]
and we put $\mathcal K(e,e):=0$. In this notation, the inner sum above is the sum of $\prod_{i\in
S}\mathcal K(e_i,e_{\mathrm{nb}(i)})$ over all distinct $e_1,\dots,e_j$. We claim that there is a
constant $C_{\mathcal K}\ge3$ such that, for all $a,b\in E$,
\begin{equation}\label{eq:Kbounds}
\sum_{e\in E}\mathcal K(a,e)\le C_{\mathcal K}n^2\qquad\text{and}\qquad\sum_{e\in E}\mathcal K(a,e)\,\mathcal K(e,b)\le C_{\mathcal K}n^2.
\end{equation}
In words, summing over the position of one edge gives at most a constant multiple of the number
$|E|\le3n^2$ of edges, both for one factor $\mathcal K$ and for a product of two factors $\mathcal
K$ that share this edge. The second bound is needed because two edges may be the nearest neighbors
of each other, in which case the factor $\mathcal K$ of this pair appears twice in the product.

We now prove \eqref{eq:Kbounds}. Let $\vartheta\in(0,2)$ and $a\in E$. For every integer $t\ge0$,
there are at most $C(t+1)$ edges $e$ with $t\le|a-e|<t+1$; moreover, $|a-e|\ge\frac12$ for $e\ne a$,
and $|a-e|\le2n$ for $e\in E$. Hence, $\sum_{e\in E\setminus\{a\}}|a-e|^{-\vartheta}\le C_\vartheta
n^{2-\vartheta}$, where $C_\vartheta$ depends only on $\vartheta$. The case
$\vartheta=\varepsilon_{\rm mc}$ gives the first bound in \eqref{eq:Kbounds}. For the second bound,
the Cauchy--Schwarz inequality and the symmetry of $\mathcal K$ give
\[
\sum_e\mathcal K(a,e)\mathcal K(e,b)\le\max_{c\in E}\sum_e\mathcal K(c,e)^2,
\]
and we use the case $\vartheta=2\varepsilon_{\rm mc}<2$.

\emph{Step 3: summing out the edges one at a time.} The map $\mathrm{nb}$ depends on the tuple
$(e_1,\dots,e_j)$, and we first remove this dependence. For every tuple, $\mathrm{nb}$ is one of the
$(j-1)^j$ maps $\psi:\{1,\dots,j\}\to\{1,\dots,j\}$ with $\psi(i)\ne i$ for all $i$. Since all terms
are nonnegative, the inner sum is therefore at most
\[
\sum_{\psi}\ \sum_{e_1,\dots,e_j\in E}\ \prod_{i\in S}\mathcal K(e_i,e_{\psi(i)}),
\]
where we have also dropped the condition that the edges be distinct. We fix $S$ and $\psi$. For
$i,l\in\{1,\dots,j\}$, we say that $l$ \emph{links} to $i$ if $l\in S$ and $\psi(l)=i$. Each factor
$\mathcal K(e_l,e_{\psi(l)})$ corresponds to a link, and each index sends at most one link. We claim
that, for every $R\subseteq\{1,\dots,j\}$,
\begin{equation}\label{eq:sumout}
\sum_{(e_i)_{i\in R}\in E^R}\ \prod_{l\in S\cap R:\ \psi(l)\in R}\mathcal K(e_l,e_{\psi(l)})\le\bigl(C_{\mathcal K}n^2\bigr)^{|R|},
\end{equation}
that is, despite the factors $\mathcal K$, each edge contributes at most a factor $C_{\mathcal
K}n^2$. We argue by induction on $|R|$. If $R=\emptyset$, both sides are equal to $1$. Let
$R\ne\emptyset$. We first sum over $e_i$, for a suitable index $i\in R$, with the other edges fixed.
The factors that do not contain $e_i$ are those of \eqref{eq:sumout} with $R\setminus\{i\}$ in place
of $R$. Hence, by the induction hypothesis, it suffices to find $i\in R$ such that the sum over
$e_i$ of the product of the factors that contain $e_i$ is at most $C_{\mathcal K}n^2$, uniformly in
the other edges. We distinguish two cases.

\emph{Case 1: some $i\in R$ receives no link from an element of $R$.} We choose such an $i$. Then
$e_i$ appears in at most one factor, namely $\mathcal K(e_i,e_{\psi(i)})$ if $i\in S$ and
$\psi(i)\in R$. Hence, the sum over $e_i$ is at most $C_{\mathcal K}n^2$ by the first bound in
\eqref{eq:Kbounds}, since $\mathcal K$ is symmetric, or at most $|E|\le3n^2$ if $e_i$ appears in no
factor.

\emph{Case 2: every element of $R$ receives a link from an element of $R$.} Since each element sends
at most one link, there are at most $|R|$ such links, and therefore every element of $R$ receives
exactly one link and sends exactly one link. We choose any $i\in R$, and we let $i'\in R$ be the
element that links to $i$. Then $e_i$ appears in exactly two factors, $\mathcal K(e_{i'},e_i)$ and
$\mathcal K(e_i,e_{\psi(i)})$, and the sum over $e_i$ is at most $C_{\mathcal K}n^2$ by the second
bound in \eqref{eq:Kbounds}. This proves \eqref{eq:sumout}.

We now conclude. By \eqref{eq:sumout} with $R=\{1,\dots,j\}$, each pair $(S,\psi)$ contributes at
most $(C_{\mathcal K}n^2)^j$, and there are at most $2^j(j-1)^j\le(2e)^jj!$ pairs $(S,\psi)$, since
$j^j\le e^jj!$. Hence, by Step~1,
\[
j!\,\E_\beta\Bigl[\binom{\mathcal N_k}{j}\Bigr]\le\bigl(C_{\rm mc}2^{(1-k)\varepsilon_{\rm mc}}\bigr)^j\,(2e)^j\,j!\,\bigl(C_{\mathcal K}n^2\bigr)^j.
\]
Dividing by $j!$ and using $2^{\varepsilon_{\rm mc}}\le2$, we obtain the lemma for $j\ge2$ with
$C_{\mathcal N}:=4eC_{\rm mc}C_{\mathcal K}$. Since $C_{\mathcal K}\ge3$, we have $C_{\mathcal
N}\ge6C_{\rm mc}$, so that the case $j=1$ treated above is also covered, and the proof is complete.
\end{proof}
\section{The specific heat in the domain slit along the interface: Proof of~(E2)}\label{sec:cov} In
this section, we prove \tagref{E2}{(E2)}. The main step is Theorem~\ref{thm:LSS}, which compares, up
to an additive constant, the sum of the energy covariances between an edge $e$ and the edges within
distance $r$ of $e$, in an arbitrary domain with plus boundary condition, with the same sum
$\chi_\beta(e;r)$ under the infinite-volume measure (\S\ref{ss:reduction}). Theorem~\ref{thm:LSS}
will be used twice: in the slit domain, to prove \tagref{E2}{(E2)} (Proposition~\ref{prop:Bminus}),
and in the square, to prove the bound on the pairs of nearby edges in \tagref{E3}{(E3)}
(Proposition~\ref{prop:near}). The comparison must be additive, since a comparison up to constant
factors would produce a total error as large as $n^2\log n$, and it must be uniform over domains
with irregular boundaries. By Lemma~\ref{lem:ES}, the energy covariances
$\langle\varepsilon_e;\varepsilon_f\rangle$ are equal to $4/\pp^2$ times those of the FK-Ising edge
variables $\omega_e$ and $\omega_f$. We therefore state and prove the estimates of
\S\S\ref{ss:covsetting}--\ref{ss:stability} for the FK-Ising model, and we return to the spin
language in Corollary~\ref{cor:chi} and Theorem~\ref{thm:LSS}. The proof has two stages. First, we
prove the a priori bound $\Cov(\omega_e,\omega_f)\lesssim\rho^{-2}$, uniformly in the domain, where
$\rho$ is the distance between $e$ and $f$ (\S\S\ref{ss:covsetting}--\ref{ss:apriori}). Second, we
show by a coupling argument that boundary conditions at distance $s$ change this covariance by at
most a fraction $(\rho/s)^\theta$ of $\rho^{-2}$, for some $\theta>0$ (\S\ref{ss:stability}).
Summing over $f$ gives Theorem~\ref{thm:LSS}, from which we deduce \tagref{E2}{(E2)}
(\S\ref{ss:LSS}).

\subsection{Setting and four-arm stability}\label{ss:covsetting} In this subsection, we state the
four-arm bound that underlies all covariance estimates, and we fix the setting of the section.

\begin{theorem}[four-arm stability]\label{thm:R1}
There are a universal constant $\varepsilon_1>0$ and, for every $\lambda<\infty$, a constant
$C_4(\lambda)$ such that, for every $N\ge1$, every $\pp$ with $|\pp-\ppc|\le \lambda/N$ \textup{(}on
either side of $\ppc$\textup{)} and all $1\le r\le R\le N$, the four-arm probability of
Definition~\ref{def:A4} satisfies
\[
\pi_4(\pp;r,R)\le C_4(\lambda)\,(r/R)^{1+\varepsilon_1}.
\]
\end{theorem}
Note that, by \eqref{eq:Lbound}, Theorem~\ref{thm:R1} concerns scales at most a constant multiple of
the characteristic length. At $\pp=\ppc$, it states that the alternating four-arm exponent of the
critical FK-Ising model is larger than $1$. We postpone the proof of Theorem~\ref{thm:R1} to
Appendix~\ref{app:R1}. It follows from results of \cite{DCM}.

We briefly explain the importance of the fact that the four-arm exponent is strictly larger than
$1$. In the proof of Lemma~\ref{lem:edgeevent} below, the covariance of an edge with an increasing
boundary-partition event of the box of radius $r$ around it (\S\ref{ss:apriori}) is bounded by a sum
over dyadic scales $2^k\le r$ of the mixing rate $\Delta_\pp(2^k)\lesssim2^{-k}$ times the four-arm
probability $\pi_4(\pp;2^k,r)$ (Definition~\ref{def:A4}). If
$\pi_4(\pp;2^k,r)\lesssim(2^k/r)^{1+\varepsilon}$, the sum is $O(1/r)$. With exponent~$1$, it would
be as large as $r^{-1}\log r$, and this is precisely the type of logarithm that the proof must
avoid.

\emph{Setting.} Recall from \eqref{eq:N} that $N=4n$. Throughout this section,
$\pp\in[\frac13,\frac34]$ is any parameter with $|\pp-\ppc|\le 8M/N$. By \eqref{eq:window}, this
covers $\pp(\beta)$ for every $\beta\in W_M$ and $n\ge n_0(M)$, on both sides of $\ppc$. We use the
following inputs: Lemma~\ref{lem:Delta}(iii) with the constant $C_\Delta(8M)$, which we abbreviate
to $C_\Delta$ in this section, together with the convention $\Delta_\pp(s)=\Delta_\pp(\lfloor
s\rfloor)$ for real $s\ge1$; property \tagref{Ra}{(Ra)} with $\lambda=40M$ and the constant
$c_1:=c_1(40M)$; property \tagref{Mb}{(Mix)} with its constant $\Kmix$; and the four-arm bound of
Theorem~\ref{thm:R1} with $\lambda=8M$, which reads
\begin{equation}\label{eq:R1use}
\pi_4(\pp;r,R)\le C_4(8M)\,(r/R)^{1+\varepsilon_1}\qquad(1\le r\le R\le N).
\end{equation}
Recall from \S\ref{ss:FK} the relation $G\supseteq\Lambda_s(e)$. We fix two edges $e\ne f$ with
$\rho:=|e-f|_\infty\ge1$. Recall from \S\ref{ss:notation} that $x_e$ denotes the left endpoint of
$e$ if $e$ is horizontal and its bottom endpoint if $e$ is vertical. For each $e$, exactly one edge
$f\ne e$ of $\Z^2$ has $|e-f|_\infty=0$, namely the edge with $x_f=x_e$, and it will be treated
separately in the proof of Theorem~\ref{thm:LSS}. For $G\supseteq\Lambda_s(e)$ and a boundary
condition $\xi$ on $G$, we let
\begin{align*}
F(G,\xi)&:=\Cov^\xi_{G,\pp}(\omega_e,\omega_f),\\
F(s,\xi)&:=F(\Lambda_s(e),\xi),\\
F(\Z^2)&:=\Cov_{\phi_\pp}(\omega_e,\omega_f),\\
\Osc(s)&:=\sup_{\xi,\xi'}|F(s,\xi)-F(s,\xi')| .
\end{align*}
By \tagref{FKG}{(FKG)}, all these covariances are nonnegative. Moreover,
$F(\Z^2)=\lim_{S\to\infty}F(\Lambda_S,1)$, since $\omega_e$ and $\omega_f$ are local.

\subsection{Elementary lemmas}\label{ss:elem} In this subsection, we collect four elementary lemmas.
The first bounds a covariance through the ranges of conditional expectations, the second bounds
these ranges, the third reduces a covariance in a general domain to covariances in a box and the
fourth bounds the covariance of an edge with an arbitrary event, scale by scale.

\begin{lemma}[total covariance]\label{lem:totcov}
Let $X$ and $Y$ be bounded random variables, and let $\cG$ be a $\sigma$-algebra. Then
\[
\Cov(X,Y)=\E[\Cov(X,Y\mid\cG)]+\Cov(\E[X\mid\cG],\E[Y\mid\cG]).
\]
If $\E[X\mid\cG]$ and $\E[Y\mid\cG]$ take values in intervals of lengths $a$ and $b$, then we have
\[\left|\Cov(\E[X\mid\cG],\E[Y\mid\cG])\right| \le \frac{ab}{4}.\]

\end{lemma}
\begin{proof}
The identity is the law of total covariance. A random variable $Z$ with values in an interval
$[c,c+a]$ satisfies
\[
\Var(Z)\le\E[(Z-c-\frac a2)^2]\le\frac{a^2}4,
\]
and the bound follows from the Cauchy--Schwarz inequality.
\end{proof}

\begin{lemma}[comparison with free and wired boundary conditions]\label{lem:sandwich}
Let $s\ge2\rho$, let $G'\supseteq\Lambda_s(e)$ and let $\zeta$ be a boundary condition on $G'$. Then
$\phi^\zeta_{G',\pp}[\omega_e]$ lies in the interval
\[
I_e(s):=[\phi^0_{\Lambda_s(e),\pp}[\omega_e],\phi^1_{\Lambda_s(e),\pp}[\omega_e]],
\]
whose length is $\Delta_\pp(s)$, and $\phi^\zeta_{G',\pp}[\omega_f]$ lies in the interval
\[
I_f(s-\rho):=[\phi^0_{\Lambda_{s-\rho}(f),\pp}[\omega_f],\phi^1_{\Lambda_{s-\rho}(f),\pp}[\omega_f]],
\]
whose length is $\Delta_\pp(s-\rho)\le\Delta_\pp(s/2)$. The same holds for $\phi_\pp$. More
generally, if $\cG$ is a $\sigma$-algebra contained in the $\sigma$-algebra generated by the states
of the edges outside $E(\Lambda_s(e))$, then $\phi^\zeta_{G',\pp}[\omega_e\mid\cG]$ takes values in
$I_e(s)$, and the same holds for $\omega_f$, with $\Lambda_{s-\rho}(f)$ and $I_f(s-\rho)$ in place
of $\Lambda_s(e)$ and $I_e(s)$.
\end{lemma}
Lemma~\ref{lem:sandwich} is a standard application of \tagref{DMP}{(DMP)} and \tagref{CBC}{(CBC)};
we omit the proof.

\begin{lemma}[reduction to boxes]\label{lem:boxes}
Let $s\ge2\rho$, let $G\supseteq\Lambda_s(e)$ with boundary condition $\xi$ and let $\zeta(\omega)$
be the boundary condition induced on $\Lambda_s(e)$ by $\omega|_{E(G)\setminus E(\Lambda_s(e))}$ and
$\xi$.
\begin{enumerate}[label=\textup{(\roman*)}]
\item $0\le F(G,\xi)-\phi^\xi_{G,\pp}\bigl[F(s,\zeta)\bigr]\le\frac14\Delta_\pp(s)\Delta_\pp(s/2)$.
\item For every $G'\supseteq\Lambda_s(e)$ with boundary condition $\xi'$, we have
\[
|F(G,\xi)-F(G',\xi')|\le\Osc(s)+\frac14\Delta_\pp(s)\Delta_\pp(s/2),
\]
and the same holds with $F(\Z^2)$ in place of $F(G',\xi')$. In particular,
\[
\Osc(s')\le\Osc(s)+\frac14\Delta_\pp(s)\Delta_\pp(s/2)
\]
for $s'\ge s$.
\end{enumerate}
\end{lemma}
\begin{proof}
First, we prove (i). We apply Lemma~\ref{lem:totcov} under $\phi^\xi_{G,\pp}$, with $\cG$ generated
by the states of the edges of $E(G)\setminus E(\Lambda_s(e))$. By \tagref{DMP}{(DMP)}, the
conditional covariance is $F(s,\zeta)$, and the conditional expectations are
$\bar\omega_e:=\phi^\zeta_{\Lambda_s(e),\pp}[\omega_e]$ and
$\bar\omega_f:=\phi^\zeta_{\Lambda_s(e),\pp}[\omega_f]$. Hence,
\[
F(G,\xi)-\phi^\xi_{G,\pp}[F(s,\zeta)]=\Cov(\bar\omega_e,\bar\omega_f).
\]
By \tagref{CBC}{(CBC)}, $\bar\omega_e$ and $\bar\omega_f$ are increasing functions of the
configuration outside $\Lambda_s(e)$, so that their covariance is nonnegative by
\tagref{FKG}{(FKG)}. By Lemma~\ref{lem:sandwich}, they take values in intervals of lengths
$\Delta_\pp(s)$ and $\Delta_\pp(s/2)$, and Lemma~\ref{lem:totcov} bounds their covariance by
$\frac14\Delta_\pp(s)\Delta_\pp(s/2)$.

For (ii), we apply (i) to $(G,\xi)$ and to $(G',\xi')$. The two expectations of $F(s,\cdot)$ given
by (i) differ by at most $\Osc(s)$, and the two remainders lie in
$[0,\frac14\Delta_\pp(s)\Delta_\pp(s/2)]$, so that they differ by at most
$\frac14\Delta_\pp(s)\Delta_\pp(s/2)$. For $F(\Z^2)$, we take $G'=\Lambda_S$ and $\xi'=1$, and we
let $S\to\infty$. For the last claim, we take $G=G'=\Lambda_{s'}(e)$.
\end{proof}

The next lemma bounds the covariance of an edge with an arbitrary event, monotone or not. In this
bound, the scale $2^j$ around the edge contributes the mixing rate at that scale times the
probability that the event still depends on the configuration inside the box of radius $2^j$. The
lemma holds for every $\pp\in(0,1)$, and it will be used in the proofs of Lemmas~\ref{lem:edgeevent}
and~\ref{lem:J}.

\begin{lemma}[multiscale influence bound]\label{lem:influence}
Let $\pp\in(0,1)$, let $G$ be a finite graph with boundary condition $\xi$, let
$\phi:=\phi^\xi_{G,\pp}$, let $A$ be an event and let $e$ be an edge with endpoint $x$. Let $k\ge0$
be such that $G\supseteq\Lambda_{2^k}(x)$. Let $\cH_{-2}$ be the $\sigma$-algebra generated by the
states of all the edges of $E(G)$; let $\cH_{-1}$ be the $\sigma$-algebra generated by the states of
the edges of $E(G)\setminus\{e\}$; and, for $0\le j\le k$, let $\cH_j$ be the $\sigma$-algebra
generated by the states of the edges of $E(G)\setminus E(\Lambda_{2^j}(x))$, so that
$\cH_{-2}\supseteq\cH_{-1}\supseteq\cH_0\supseteq\dots\supseteq\cH_k$. Let
\[U_j:=\{\omega: \; \phi[A\mid\cH_j]\in(0,1)\}\in\cH_j,\]
for $-1\le j\le k$. Then
\begin{equation}\label{eq:influence}
|\Cov_\phi(\omega_e,\one_A)|\le2\sum_{j=-1}^k\bar\Delta_j\,\phi[A\cap U_j]+\Delta_\pp(2^k)\,\phi[A],
\end{equation}
where $\bar\Delta_{-1}=\bar\Delta_0:=1$ and $\bar\Delta_j:=\Delta_\pp(2^{j-1})$ for $j\ge1$.
Moreover, $U_{-1}\subseteq U_0\subseteq\dots\subseteq U_k$, and $\omega\in U_j$ if and only if $A$
is not determined by the restriction of $\omega$ to the edges whose states generate $\cH_j$.
\end{lemma}
\begin{proof}
The idea is to write $\one_A$ as a sum of martingale increments along the decreasing sequence of
$\sigma$-algebras $(\cH_j)$, and to bound the covariance of $\omega_e$ with each increment. Let
\[
D_j:=\phi[A\mid\cH_{j-1}]-\phi[A\mid\cH_j]
\]
for $-1\le j\le k$. Since $A$ is $\cH_{-2}$-measurable, we have $\phi[A\mid\cH_{-2}]=\one_A$, so
that the sum of the $D_j$ telescopes, and $\one_A=\sum_{j=-1}^kD_j+\phi[A\mid\cH_k]$. By
bilinearity,
\[
\Cov_\phi(\omega_e,\one_A)=\sum_{j=-1}^k\Cov_\phi(\omega_e,D_j)+\Cov_\phi\bigl(\omega_e,\phi[A\mid\cH_k]\bigr).
\]
We now rewrite each term $\Cov_\phi(\omega_e,D_j)$. Since $\cH_j\subseteq\cH_{j-1}$, the tower
property gives $\phi[D_j\mid\cH_j]=0$, and in particular $\phi[D_j]=0$, so that
$\Cov_\phi(\omega_e,D_j)=\phi[\omega_eD_j]$. Since $D_j$ is $\cH_{j-1}$-measurable, we have
$\phi[\omega_eD_j]=\phi\bigl[\phi[\omega_e\mid\cH_{j-1}]\,D_j\bigr]$. Since
$\phi[\omega_e\mid\cH_j]$ is $\cH_j$-measurable, we also have
\[
\phi\bigl[\phi[\omega_e\mid\cH_j]\,D_j\bigr]=\phi\bigl[\phi[\omega_e\mid\cH_j]\,\phi[D_j\mid\cH_j]\bigr]=0.
\]
Subtracting this identity from the previous one, we obtain
\[
\Cov_\phi(\omega_e,D_j)=\phi\bigl[(\phi[\omega_e\mid\cH_{j-1}]-\phi[\omega_e\mid\cH_j])\,D_j\bigr].
\]
In other words, only the change of the conditional expectation of $\omega_e$ from $\cH_{j-1}$ to
$\cH_j$ contributes to the covariance of $\omega_e$ with $D_j$. Bounding this change by its supremum
over all configurations, and using the triangle inequality, we obtain
\begin{equation}\label{eq:influencesplit}
|\Cov_\phi(\omega_e,\one_A)|\le\sum_{j=-1}^k\sup\big|\phi[\omega_e\mid\cH_{j-1}]-\phi[\omega_e\mid\cH_j]\big|\,\phi\bigl[|D_j|\bigr]
+\big|\Cov_\phi\bigl(\omega_e,\phi[A\mid\cH_k]\bigr)\big|.
\end{equation}
We bound the three ingredients of \eqref{eq:influencesplit} in turn.

\emph{Step 1: the increments of the conditional expectation of $\omega_e$.} We claim that
\[
|\phi[\omega_e\mid\cH_{j-1}]-\phi[\omega_e\mid\cH_j]|\le\bar\Delta_j.
\]
For $j\in\{-1,0\}$, this holds since both conditional expectations lie in $[0,1]$. Let $j\ge1$. By
\tagref{DMP}{(DMP)}, given $\cH_{j-1}$, the configuration on $E(\Lambda_{2^{j-1}}(x))$ is
distributed as an FK configuration on this box with the induced boundary condition. Since $e$ is
incident to the center $x$ of the box, \tagref{CBC}{(CBC)} puts $\phi[\omega_e\mid\cH_{j-1}]$
between the expectations of $\omega_e$ under the free and the wired FK measures on
$\Lambda_{2^{j-1}}(x)$, that is, in an interval of length $\Delta_\pp(2^{j-1})$, as for $I_e$ in
Lemma~\ref{lem:sandwich}. Since $\cH_j\subseteq\cH_{j-1}$, the tower property gives
\[
\phi[\omega_e\mid\cH_j]=\phi\bigl[\phi[\omega_e\mid\cH_{j-1}]\bigm|\cH_j\bigr],
\]
which lies in the same interval. This proves the claim.

\emph{Step 2: the increment $D_j$ vanishes on $(U_j)^c$.} On $(U_j)^c$, we have
$\phi[A\mid\cH_j]\in\{0,1\}$. Since $\phi[A\mid\cH_{j-1}]$ takes values in $[0,1]$ and has
conditional expectation $\phi[A\mid\cH_j]$ given $\cH_j$, it equals $\phi[A\mid\cH_j]$ there. Hence,
$D_j=0$ on $(U_j)^c$, and, since $U_j\in\cH_j\subseteq\cH_{j-1}$,
\[
\phi\bigl[|D_j|\bigr]\le\phi\bigl[(\phi[A\mid\cH_{j-1}]+\phi[A\mid\cH_j])\one_{U_j}\bigr]=2\phi[A\cap U_j].
\]
Note that we have also proved that $\phi[A\mid\cH_{j-1}]\in\{0,1\}$ on $(U_j)^c$, that is,
$U_{j-1}\subseteq U_j$ for $0\le j\le k$.

\emph{Step 3: the remainder term.} We have
\[
\Cov_\phi(\omega_e,\phi[A\mid\cH_k])=\phi\bigl[(\phi[\omega_e\mid\cH_k]-\phi[\omega_e])\,\phi[A\mid\cH_k]\bigr].
\]
By the argument of Step~1, applied to the box $\Lambda_{2^k}(x)$, both $\phi[\omega_e\mid\cH_k]$ and
its expectation $\phi[\omega_e]$ lie in an interval of length $\Delta_\pp(2^k)$, so that the
absolute value of this covariance is at most $\Delta_\pp(2^k)\,\phi[A]$.

Inserting these three bounds into \eqref{eq:influencesplit} gives \eqref{eq:influence}. Finally, by
\tagref{FE}{(FE)}, given $\cH_j$, every configuration of the remaining edges has positive
conditional probability. Hence, $\phi[A\mid\cH_j](\omega)\in\{0,1\}$ if and only if $A$ is
determined by the restriction of $\omega$ to the edges whose states generate $\cH_j$.
\end{proof}
\subsection{The a priori covariance bound}\label{ss:apriori} In this subsection, we prove that the
covariance of two edges at distance $\rho$ is at most a constant multiple of $\rho^{-2}$, in every
domain and for every boundary condition. The bound is the product of two factors $\rho^{-1}$.
Conditioning on the box of radius $\rho/2$ around $e$ turns the covariance of $\omega_e$ with
$\omega_f$ into an integral of covariances of $\omega_e$ with increasing boundary-partition events
of that box (defined below). Lemma~\ref{lem:edgeevent} bounds each of these covariances by a
constant divided by the radius of the box, which gives one factor $\rho^{-1}$. The other factor is
the length of the interval in which the conditional expectation of $\omega_f$ takes its values,
which is at most the mixing rate at scale $\rho/4$.

For $r\ge1$, let $\zeta_r(\omega)$ be the partition of $\partial\Lambda_r(e)$ in which two vertices
are in the same block if and only if they are joined by an open path with all its edges in
$E(\Lambda_r(e))$. An event $\mathcal B$ is a \emph{boundary-partition event of $\Lambda_r(e)$} if
it is measurable with respect to the $\sigma$-algebra generated by $\zeta_r$.

\begin{lemma}\label{lem:edgeevent}
There is $C_{\mathcal B}=C_{\mathcal B}(M)$ such that, for all $1\le r\le N/2$, every finite graph
$G\supseteq\Lambda_r(e)$, every boundary condition $\xi$ and every increasing boundary-partition
event $\mathcal B$ of $\Lambda_r(e)$, we have $|\Cov^\xi_{G,\pp}(\omega_e,\one_{\mathcal B})|\le
C_{\mathcal B}/r$.
\end{lemma}
\begin{proof}
We write $\phi:=\phi^\xi_{G,\pp}$. Recall from \S\ref{ss:inputs} that $\Kmix\ge4$ is the integer
constant of \tagref{Mb}{(Mix)}: if an event depends only on the states of the edges of a box
$\Lambda_{r'}(x)$, then its probabilities under $\phi_\pp$ and under any FK measure on a finite
graph containing $\Lambda_{\Kmix r'}(x)$, with any boundary condition, differ by a factor of at most
$\frac32$. If $r<4\Kmix$, the bound holds for every $C_{\mathcal B}\ge\Kmix$, since
$|\Cov_\phi(\omega_e,\one_{\mathcal B})|\le\frac14$. Hence, we may assume that $r\ge4\Kmix$. Let
$s:=\lfloor r/\Kmix\rfloor\ge4$, and let $k_r$ be the largest integer with $2^{k_r}\le r/2$.
Lemma~\ref{lem:influence} with $x=x_e$, $k=k_r$ and $A=\mathcal B$ gives
\[
|\Cov_\phi(\omega_e,\one_{\mathcal B})|\le2\sum_{j=-1}^{k_r}\bar\Delta_j\,\phi[\mathcal B\cap U_j]+\Delta_\pp(2^{k_r}).
\]
It remains to bound the probabilities $\phi[\mathcal B\cap U_j]$, and we do this in three steps.

\emph{Step 1: nondetermination forces four arms.} We claim that $\mathcal B\cap U_j\subseteq
A^{x_e}_4(2^j,r)$ for $0\le j\le k_r$; since $U_{-1}\subseteq U_0$, this also bounds the term
$j=-1$. Let $\omega\in\mathcal B\cap U_j$, and let $\omega^{[j]}$ be the configuration obtained from
$\omega$ by closing every edge of $E(\Lambda_{2^j}(e))$. Since $\mathcal B$ is increasing and, by
Lemma~\ref{lem:influence}, is not determined by the restriction of $\omega$ to the edges outside
$E(\Lambda_{2^j}(e))$, we have $\omega^{[j]}\notin\mathcal B$, and hence
$\zeta_r(\omega^{[j]})\ne\zeta_r(\omega)$. Since $\omega\ge\omega^{[j]}$, there are vertices
$y,y'\in\partial\Lambda_r(e)$ that are joined by an $\omega$-open path with all its edges in
$E(\Lambda_r(e))$ but lie in distinct $\omega^{[j]}$-clusters $\mathcal C\ne\mathcal C'$ of
$\Lambda_r(e)$. The path uses an edge of $E(\Lambda_{2^j}(e))$, since otherwise it would be
$\omega^{[j]}$-open and we would have $\mathcal C=\mathcal C'$. Hence, its initial and final
segments show that $\mathcal C$ and $\mathcal C'$ both meet $\partial\Lambda_{2^j}(e)$. Since every
edge of $E(\Lambda_{2^j}(e))$ is $\omega^{[j]}$-closed, $\mathcal C$ and $\mathcal C'$ are distinct
clusters of $\omega|_{E_{2^j,r}(x_e)}$, each of which meets $\partial\Lambda_{2^j}(e)$ and
$\partial\Lambda_r(e)$. By planar duality, two such clusters give four alternating arms, so that
$\omega\in A^{x_e}_4(2^j,r)$.

\emph{Step 2: arm bounds.} If $2^j\le s$, then
\[
A^{x_e}_4(2^j,r)\subseteq A^{x_e}_4(2^j,s)\in\cF(\Lambda_s(e))
\]
by \eqref{eq:A4mono}, and $G\supseteq\Lambda_r(e)\supseteq\Lambda_{\Kmix s}(e)$. Hence, by
\tagref{Mb}{(Mix)} and \eqref{eq:R1use},
\[
\phi[A^{x_e}_4(2^j,r)]\le\frac32C_4(8M)(2^j/s)^{1+\varepsilon_1}.
\]
If $2^j>s$, we use the trivial bound $1$.

\emph{Step 3: summation.} By Steps~1 and~2, together with $\bar\Delta_{-1}=\bar\Delta_0=1$,
$\bar\Delta_j=\Delta_\pp(2^{j-1})\le2C_\Delta2^{-j}$ for $j\ge1$ (Lemma~\ref{lem:Delta}(iii)) and
$\Delta_\pp(2^{k_r})\le C_\Delta2^{-k_r}\le4C_\Delta/r$, we obtain
\[
|\Cov_\phi(\omega_e,\one_{\mathcal B})|\le\frac{C}{s^{1+\varepsilon_1}}+C\sum_{j:\,2^j\le s}2^{-j}\Bigl(\frac{2^j}s\Bigr)^{1+\varepsilon_1}
+C\sum_{j:\,2^j>s}2^{-j}+\frac Cr\le\frac Cs+\frac Cr .
\]
The terms of the first sum are $2^{j\varepsilon_1}s^{-1-\varepsilon_1}$, so that this sum is at most
a constant times its largest term, which is at most $s^{-1}$. This is where the exponent
$1+\varepsilon_1>1$ is used. The second sum is at most $2/s$. Since $s\ge r/(2\Kmix)$, the claim
follows.
\end{proof}

We now prove the a priori bound.

\begin{proposition}[a priori bound, uniform in domains]\label{prop:apriori}
There is $C_{\rm ap}=C_{\rm ap}(M)$ such that, for all edges $e\ne f$ with
$1\le\rho=|e-f|_\infty\le\frac45N$, every finite graph $G\supseteq\Lambda_{\lceil5\rho/4\rceil}(e)$
and every boundary condition $\xi$,
\[
0\le\Cov^\xi_{G,\pp}(\omega_e,\omega_f)\le C_{\rm ap}\,\rho^{-2} .
\]
The same bound holds for $\phi_\pp$.
\end{proposition}
For the infinite-volume measure $\phi_\pp$ and edges at distance at most the characteristic length
$L(\pp)$ of \eqref{eq:Ldef}, the bound of Proposition~\ref{prop:apriori} is also stated in
\cite[Thm.~1.13, (1.32)]{DCM}, together with a matching lower bound. It follows from
\cite[Lemma~5.3, (5.3)]{DCM} applied to a single edge, as in \cite[\S8.2.4, p.~78]{DCM}, and the
proof of that lemma bounds the covariance by the same product of two factors as above. Duminil-Copin
and Manolescu also state \cite[Remark~5.5]{DCM} that their proof of \cite[Lemma~5.3]{DCM} adapts
readily to subgraphs of $\Z^2$ that contain a box twice as large as the quad of that lemma, with
arbitrary boundary conditions, for edges closer to the quad than to the complement of the subgraph.
We give a complete proof of the upper bound in this generality, which is what
\S\S\ref{ss:stability}--\ref{ss:LSS} and Section~\ref{sec:VP} require, with a different argument for
the factor coming from $e$, namely Lemma~\ref{lem:edgeevent}.
\begin{proof}
We write $\phi:=\phi^\xi_{G,\pp}$. The covariance is nonnegative by \tagref{FKG}{(FKG)}. If
$\rho<8$, the upper bound holds for every $C_{\rm ap}\ge16$, since
$\Cov^\xi_{G,\pp}(\omega_e,\omega_f)\le\frac14\le16\rho^{-2}$. Hence, we may assume that $\rho\ge8$.
Let $r:=\lfloor\rho/2\rfloor$, $U:=\Lambda_r(e)$ and $V:=\Lambda_{\lfloor\rho/4\rfloor}(f)$. Note
that $e\in E(U)$, $V\subseteq\Lambda_{\lceil5\rho/4\rceil}(e)$ and $V\cap U=\emptyset$.

\emph{Step 1: conditioning on $U$.} Given $\omega|_{E(U)}$, the configuration on $E(G)\setminus
E(U)$ has law $\phi^{\xi\sqcup\zeta_r}_{G\setminus U,\pp}$ by \tagref{DMP}{(DMP)}. Here, $G\setminus
U$ denotes the subgraph of $G$ with edge set $E(G)\setminus E(U)$, whose boundary vertices are the
boundary vertices of $G$ and the vertices of $\partial U$, and $\xi\sqcup\zeta_r$ denotes the
boundary condition on $G\setminus U$ that coincides with $\xi$ on the boundary vertices of $G$ and
with $\zeta_r$ on $\partial U$. Indeed, open paths inside $U$ only join vertices of $\partial U$,
and $\partial U$ contains no boundary vertex of $G$, since
$G\supseteq\Lambda_{\lceil5\rho/4\rceil}(e)$. Hence,
$\phi[\omega_f\mid\omega|_{E(U)}]=\Psi(\zeta_r)$ for a function $\Psi$ of the boundary partition. By
\tagref{CBC}{(CBC)}, $\Psi$ is increasing in the partition, and, as in Lemma~\ref{lem:sandwich}, it
takes values in the interval $[a,b]$ with $a:=\phi^0_{V,\pp}[\omega_f]$ and
$b:=\phi^1_{V,\pp}[\omega_f]$, whose length is
$b-a=\Delta_\pp(\lfloor\rho/4\rfloor)\le8C_\Delta\rho^{-1}$. This is the first factor $\rho^{-1}$.

\emph{Step 2: reduction to increasing boundary-partition events.} With $a$ and $b$ as in Step~1, we
use the elementary identity
\[
y=a+\int_0^{b-a}\one\{y>a+v\}\,dv,\qquad \text{for all } a\le y\le b,
\]
which holds since the integrand is equal to $1$ exactly for $v\in[0,y-a)$. Applied to
$y=\Psi(\zeta_r)$, which lies in $[a,b]$ by Step~1, it writes the increasing function
$\Psi(\zeta_r)$ as an integral of indicator functions of the events $\mathcal
B_v:=\{\Psi(\zeta_r)>a+v\}$, namely $\Psi(\zeta_r)=a+\int_0^{b-a}\one_{\mathcal B_v}\,dv$. Since
$\omega_e$ is $\omega|_{E(U)}$-measurable, the covariance of $\omega_e$ and $\omega_f$ therefore
becomes an integral of covariances of $\omega_e$ with these events:
\[
\Cov(\omega_e,\omega_f)=\Cov(\omega_e,\Psi(\zeta_r))=\int_0^{b-a}\Cov\bigl(\omega_e,\one_{\mathcal B_v}\bigr)\,dv.
\]
Since $\Psi$ is increasing, each $\mathcal B_v$ is an increasing boundary-partition event of
$\Lambda_r(e)$. Moreover, $G\supseteq\Lambda_r(e)$ and $r\le\frac25N$, so that
Lemma~\ref{lem:edgeevent} gives the second factor,
\[
|\Cov(\omega_e,\one_{\mathcal B_v})|\le C_{\mathcal B}/r\le4C_{\mathcal B}\rho^{-1}.
\]
Hence, $\Cov(\omega_e,\omega_f)\le32C_\Delta C_{\mathcal B}\rho^{-2}$. For $\phi_\pp$, we apply this
bound with $G=\Lambda_S$ and $\xi=1$, and we let $S\to\infty$.
\end{proof}

In the spin language, Proposition~\ref{prop:apriori} bounds the energy covariances by a constant
times $|e-f|^{-2}$, so that $\chi_\beta(e;r)$ grows at most logarithmically in $r$. Recall from
\S\ref{ss:reduction} that
\[
\chi_\beta(e;r)=\sum_{f:\,|f-e|\le r}\langle\varepsilon_e;\varepsilon_f\rangle_{\Z^2,\beta}.
\]
\begin{corollary}[infinite-volume energy covariances]\label{cor:chi}
There is $C_\chi=C_\chi(M)$ such that, for $\beta\in W_M$ and edges $e\ne f$ with
$|e-f|\le\frac12N$, we have
\[
0\le\langle\varepsilon_e;\varepsilon_f\rangle_{\Z^2,\beta}\le C_\chi(1\vee|e-f|)^{-2}.
\]
Consequently, for $1\le r\le R\le\frac12N$,
\begin{equation}\label{eq:chi}
0\le\chi_\beta(e;R)-\chi_\beta(e;r)\le C\Bigl(1+\log\frac Rr\Bigr),\qquad0\le\chi_\beta(e;R)\le C(1+\log R).
\end{equation}
\end{corollary}
\begin{proof}
By Lemma~\ref{lem:ES},
\[
\langle\varepsilon_e;\varepsilon_f\rangle_{\Z^2,\beta}=\frac4{\pp^2}\Cov_{\phi_\pp}(\omega_e,\omega_f),
\]
which is nonnegative by \tagref{FKG}{(FKG)}, and at most $1$ since the energies take values in
$\{\pm1\}$. The first claim follows from Proposition~\ref{prop:apriori} if $|e-f|_\infty\ge1$, since
$|e-f|\le\sqrt2|e-f|_\infty+1$, and from the bound $1$ otherwise, since then $|e-f|\le1$. For
\eqref{eq:chi}, note that, for every integer $j\ge1$, there are at most $Cj$ edges $f$ with
$j<|f-e|\le j+1$, and that the first claim gives
$\langle\varepsilon_e;\varepsilon_f\rangle_{\Z^2,\beta}\le C_\chi j^{-2}$ for each of them. Summing
over these edges and over $j$ gives \eqref{eq:chi}, since the edges $f$ with $|f-e|\le1$ contribute
at most a constant.
\end{proof}
\subsection{Stability of covariances under boundary changes}\label{ss:stability}
Proposition~\ref{prop:apriori} gives an upper bound of the correct size, but not the additive
precision that Theorem~\ref{thm:LSS} requires: an error by a constant factor in each of the
covariances of $\omega_e$ with the $\asymp r^2$ edges within distance $r$ of $e$ would add up to an
error $\asymp\log r$. In this subsection, we prove that boundary conditions at distance $s$ change
$\Cov(\omega_e,\omega_f)$ by at most a fraction $(\rho/s)^\theta$ of its natural size $\rho^{-2}$,
for some $\theta>0$ (Proposition~\ref{prop:stab}). The main step is Lemma~\ref{lem:contraction},
which shows that changing the boundary condition far away contracts the oscillation of the
covariance. The idea is to couple the FK measures with two boundary conditions on $\Lambda_s(e)$
monotonically, and to explore from $\partial\Lambda_s(e)$ the cluster of the boundary. If a closed
dual circuit around $\Lambda_r(e)$ stops the exploration, which happens in each dyadic annulus with
probability at least $c_1$, then both conditional laws inside the circuit are the same.

\begin{lemma}[one-step coupling]\label{lem:contraction}
Let $r\ge2\rho$ and $k\ge1$ be such that $s:=2^kr\le N$. Then, for every boundary condition $\xi$ on
$\Lambda_s(e)$,
\[
|F(s,\xi)-F(s,0)|\le(1-c_1)^k\,\Osc(r)+\tfrac12\Delta_\pp(r)\Delta_\pp(r/2).
\]
Consequently,
\[
\Osc(s)\le2(1-c_1)^k\Osc(r)+\Delta_\pp(r)\Delta_\pp(r/2).
\]
\end{lemma}

\begin{proof}
We write $B:=\Lambda_s(e)$ and
\[
V_A:=\{x:r<|x-x_e|_\infty\le s\}
\]
for the annulus between $\Lambda_r(e)$ and $\partial B$. We proceed in four steps: we construct an
exploration from $\partial B$ (Step~1); we show that, with probability at least $1-(1-c_1)^k$, a
closed dual circuit surrounds $\Lambda_r(e)$ (Step~2); we show that the two conditional laws inside
then coincide (Step~3); and we conclude with the law of total covariance (Step~4).

\emph{Step 1: the exploration.} By Lemma~\ref{lem:coupling}, we couple $\omega^0\sim\phi^0_{B,\pp}$
and $\omega\sim\phi^\xi_{B,\pp}$, under a probability measure $\Pp$ with expectation $\E$, so that
$\omega^0\le\omega$, along an exploration driven by the larger configuration $\omega$. Let $\mathcal
E$ be the event that an $\omega$-closed dual circuit that crosses only edges of $E(B)\setminus
E(\Lambda_r(e))$ surrounds $\Lambda_r(e)$. Starting from $\partial B$, the exploration moves inward,
toward $\Lambda_r(e)$: it reveals, one at a time, the edges of $E(B)$ with an endpoint in the set
$\mathcal C$ of vertices of $V_A$ that are joined to $\partial B$ by $\omega$-open paths in $V_A$,
in an order that depends only on the states of the edges already revealed. The exploration stops the
first time that the revealed closed edges contain a closed dual circuit surrounding $\Lambda_r(e)$.
This happens if and only if $\mathcal E$ occurs: on $\mathcal E$, the circuit of $\mathcal E$
separates $\mathcal C\supseteq\partial B$ from $\Lambda_r(e)$, so that all the edges between
$\mathcal C$ and the component of $B\setminus\mathcal C$ that contains $\Lambda_r(e)$ are closed,
and these edges contain such a circuit. On $\mathcal E^c$, once all the edges with an endpoint in
$\mathcal C$ have been revealed, the exploration reveals all the remaining edges of $E(B)\setminus
E(\Lambda_r(e))$. Let $\cG$ be the $\sigma$-algebra generated by the exploration; it corresponds to
a stopping time of the filtration of Lemma~\ref{lem:coupling}.

\emph{Step 2: the probability of $\mathcal E$.} By \tagref{Ra}{(Ra)}, applied successively in the
$k$ disjoint annuli $\Lambda_{2^{i+1}r}(e)\setminus\Lambda_{2^ir}(e)$, $0\le i<k$, each time
conditionally on the configuration outside the annulus, a closed dual circuit surrounding
$\Lambda_r(e)$ exists in at least one of these annuli with probability at least $1-(1-c_1)^k$.
Therefore,
\begin{equation}\label{eq:PEc}
\Pp(\mathcal E^c)\le(1-c_1)^k .
\end{equation}

\emph{Step 3: the conditional laws on $\mathcal E$ and on $\mathcal E^c$.} On $\mathcal E$, let $U$
be the set of vertices of $B$ inside the circuit found by the exploration. Since this circuit
surrounds $\Lambda_r(e)$ and crosses only edges of $E(B)\setminus E(\Lambda_r(e))$, we have
$\Lambda_r(e)\subseteq U$; in particular, $e,f\in E(U)$, since $r\ge2\rho$. Every edge between $U$
and $B\setminus U$ is crossed by the circuit; hence, it has been revealed and is $\omega$-closed,
and it is also $\omega^0$-closed, since $\omega^0\le\omega$. No edge of $E(U)$ has been revealed,
and $U$ contains no vertex of $\partial B$, so that neither boundary condition acts on $U$. By
Lemma~\ref{lem:coupling}, the conditional law of $\omega$ given $\cG$ is $\phi^\xi_{B,\pp}$
conditioned on the revealed states. Since $U$ is separated from the rest of $B$ by closed edges, the
number of clusters factorizes, and the conditional law of $\omega$ on $E(U)$ is the free measure
$\phi^0_{U,\pp}$. The same argument applies to $\omega^0$. Therefore, on $\mathcal E$, both
conditional laws on $E(U)$ are equal to $\phi^0_{U,\pp}$, and
\[
\Cov(\omega_e,\omega_f\mid\cG)=\Cov(\omega^0_e,\omega^0_f\mid\cG),\qquad\E[\omega_e\mid\cG]=\E[\omega^0_e\mid\cG],\quad\E[\omega_f\mid\cG]=\E[\omega^0_f\mid\cG].
\]
On $\mathcal E^c$, all the edges of $E(B)\setminus E(\Lambda_r(e))$ are revealed. By
\tagref{DMP}{(DMP)}, the conditional laws of $\omega$ and $\omega^0$ on $E(\Lambda_r(e))$ are
$\phi^\zeta_{\Lambda_r(e),\pp}$ and $\phi^{\zeta^0}_{\Lambda_r(e),\pp}$, where $\zeta$ and $\zeta^0$
are the induced boundary conditions. The conditional covariances are therefore $F(r,\zeta)$ and
$F(r,\zeta^0)$, which differ by at most $\Osc(r)$.

\emph{Step 4: conclusion.} For $e'\in\{e,f\}$, let $\bar\omega_{e'}:=\E[\omega_{e'}\mid\cG]$ and
$\bar\omega^0_{e'}:=\E[\omega^0_{e'}\mid\cG]$. By Lemma~\ref{lem:totcov}, applied to $\omega$ and to
$\omega^0$ with the same $\sigma$-algebra $\cG$,
\[
F(s,\xi)-F(s,0)=\E\bigl[\Cov(\omega_e,\omega_f\mid\cG)-\Cov(\omega^0_e,\omega^0_f\mid\cG)\bigr]+\bigl[\Cov(\bar\omega_e,\bar\omega_f)-\Cov(\bar\omega^0_e,\bar\omega^0_f)\bigr].
\]
By Step~3, the difference of the two conditional covariances in the first term vanishes on $\mathcal
E$, and its absolute value is at most $\Osc(r)$ on $\mathcal E^c$. Hence, by \eqref{eq:PEc}, the
first term is at most $(1-c_1)^k\Osc(r)$ in absolute value. For the second term, note that, in all
cases, the conditional laws on $E(\Lambda_r(e))$ are restrictions of FK measures on graphs
containing $\Lambda_r(e)$. Hence, $\bar\omega_e,\bar\omega^0_e\in I_e(r)$ and
$\bar\omega_f,\bar\omega^0_f\in I_f(r-\rho)$ by Lemma~\ref{lem:sandwich}, and, by
Lemma~\ref{lem:totcov}, each of the two covariances of conditional expectations is at most
$\frac14\Delta_\pp(r)\Delta_\pp(r/2)$ in absolute value. This proves the first claim. The second
claim follows from
\[
\Osc(s)\le2\sup_\xi|F(s,\xi)-F(s,0)|.\qedhere
\]
\end{proof}

We now iterate Lemma~\ref{lem:contraction} over geometrically growing scales, starting from the a
priori bound. Let
\[
k_0:=\min\{k\ge1:2(1-c_1)^k\le\frac12\},
\]
and let $\theta:=1/k_0\in(0,1]$.

\begin{proposition}[covariance stability]\label{prop:stab}
There is $C_{\rm st}=C_{\rm st}(M)$ such that, for all edges $e\ne f$ with $\rho=|e-f|_\infty\ge1$
and all $s$ with $4\rho\le s\le N$, the following hold.
\begin{enumerate}[label=\textup{(\roman*)}]
\item $\Osc(s)\le C_{\rm st}\rho^{-2}(\rho/s)^\theta$.
\item For every finite graph $G\supseteq\Lambda_s(e)$ and every boundary condition $\xi$,
\[
|\Cov^\xi_{G,\pp}(\omega_e,\omega_f)-\Cov_{\phi_\pp}(\omega_e,\omega_f)|\le C_{\rm st}\rho^{-2}(\rho/s)^\theta.
\]
\end{enumerate}
\end{proposition}
\begin{proof}
We iterate Lemma~\ref{lem:contraction} along the scales $s_0:=4\rho$ and $s_j:=2^{jk_0}s_0$,
starting from the a priori bound, and we interpolate between consecutive scales with
Lemma~\ref{lem:boxes}(ii).

\emph{Step 1: the initial scale.} Since $F\ge0$,
$\Lambda_{4\rho}(e)\supseteq\Lambda_{\lceil5\rho/4\rceil}(e)$ and $\rho\le N/4$,
Proposition~\ref{prop:apriori} gives
\[
\Osc(s_0)\le\sup_\xi F(s_0,\xi)\le C_{\rm ap}\rho^{-2}.
\]

\emph{Step 2: the recursion.} Let $j\ge0$ be such that $s_{j+1}\le N$. Lemma~\ref{lem:contraction}
with $r=s_j$ and $k=k_0$, together with the choice of $k_0$ and Lemma~\ref{lem:Delta}(iii), gives
\[
\Osc(s_{j+1})\le\tfrac12\Osc(s_j)+Cs_j^{-2},\qquad\text{where }s_j^{-2}\le s_0^{-2}4^{-j}.
\]
Hence, $a_j:=2^j\Osc(s_j)$ satisfies $a_{j+1}\le a_j+2Cs_0^{-2}2^{-j}$. Summing these increments, we
obtain $a_j\le a_0+4Cs_0^{-2}$, and therefore, by Step~1, $\Osc(s_j)\le C'\rho^{-2}2^{-j}$ for every
$j$ with $s_j\le N$.

\emph{Step 3: interpolation.} Let $j\ge0$ be such that $s_j\le s<s_{j+1}$; in particular, $s_j\le
N$. By Lemma~\ref{lem:boxes}(ii) and Lemma~\ref{lem:Delta}(iii), we have
\[
\Osc(s)\le\Osc(s_j)+Cs_j^{-2}\le C''\rho^{-2}2^{-j}.
\]
Moreover, by the definition of $\theta$ and since $s<s_{j+1}$, we have
\[
2^{-j}=2(s_0/s_{j+1})^\theta\le2(4\rho/s)^\theta\le8(\rho/s)^\theta.
\]
This proves (i). For (ii), Lemma~\ref{lem:boxes}(ii), applied with $F(\Z^2)$, bounds the left-hand
side by $\Osc(s)+Cs^{-2}$, and $s^{-2}\le\rho^{-2}(\rho/s)^\theta$, since $\rho\le s$ and
$\theta\le1$.
\end{proof}

\subsection{Local specific-heat stability and the cancellation of the bulk logarithm: Proof
of~(E2)}\label{ss:LSS} In this subsection, we sum Proposition~\ref{prop:stab} over the edges $f$
within distance $r$ of $e$, which gives the additive comparison of local specific heats
(Theorem~\ref{thm:LSS}), the main result of this section, and we deduce \tagref{E2}{(E2)} from it
(Proposition~\ref{prop:Bminus}). Since the domain enters Proposition~\ref{prop:stab} only through
the boundary condition that it induces on $\Lambda_s(e)$, no regularity of the boundary is needed,
and Theorem~\ref{thm:LSS} applies to the components of the slit domain, whose boundaries are
irregular. Recall from \S\ref{ss:FK} that $\langle\cdot\rangle^+_{D,\beta}$ denotes the expectation
with respect to the Ising measure on $D$ with plus boundary condition.
\begin{theorem}[local specific-heat stability]\label{thm:LSS}
There is $C_{\rm LS}=C_{\rm LS}(M)$ such that the following holds for all $n\ge n_0(M)$ and
$\beta\in W_M$. Let $D\subseteq\Omega_n$, let $e$ be an edge and let $r\ge1$ be such that every
vertex at Euclidean distance at most $2r$ from $\mathrm{mid}(e)$ belongs to $D$. Then
\begin{equation}\label{eq:LSS}
\Big|\sum_{f:\,|f-e|\le r}\langle\varepsilon_e;\varepsilon_f\rangle^+_{D,\beta}-\chi_\beta(e;r)\Big|\le C_{\rm LS} .
\end{equation}
\end{theorem}
\begin{proof}
We compare the two sums term by term. For the edges $f$ with $\rho=|f-e|_\infty$ up to a quarter of
the box size $s$ defined below, we use Proposition~\ref{prop:stab}, whose factor $(\rho/s)^\theta$
keeps the sum bounded uniformly in $s$; for the remaining edges $f$, which occupy a bounded number
of dyadic scales, the a priori bound suffices.

All the covariances involved lie in $[0,1]$: they are nonnegative by Lemma~\ref{lem:ES} and
\tagref{FKG}{(FKG)}, and at most $1$ since the energies take values in $\{\pm1\}$. Every $f$ with
$|f-e|\le r$ has both endpoints within distance $2r$ of $\mathrm{mid}(e)$, and hence in $D$. Let
$r_0$ be a large constant. If $r<r_0$, both sums have at most $Cr_0^2$ terms, and the claim holds
for every $C_{\rm LS}\ge Cr_0^2$. Hence, we may assume that $r\ge r_0$. Note also that $r<n$, since
all the vertices within distance $2r$ of $\mathrm{mid}(e)$ lie in $D\subseteq\Omega_n$.

Let $\pp:=\pp(\beta)$ and $s:=\lfloor4r/3\rfloor$. Since $4\sqrt2/3<2$ and $r\ge r_0$, every vertex
$x$ with $|x-x_e|_\infty\le s$ satisfies $|x-\mathrm{mid}(e)|\le\sqrt2s+1<2r$, so that
$G_D\supseteq\Lambda_s(e)$ (see \S\ref{ss:FK} for $G_D$). Similarly, since $\frac54<\frac43$ and
$r\ge r_0$, we have $\lceil\frac54(r+1)\rceil\le s$; in particular, $r+1\le s$. Moreover, $s\le N$,
since $r<n$. Every $f$ with $|f-e|\le r$ satisfies $\rho:=|f-e|_\infty\le r+1$, since
$|f-e|_\infty\le|f-e|+1$ (\S\ref{ss:notation}). By Lemma~\ref{lem:ES}, for $f\ne e$,
\[
\langle\varepsilon_e;\varepsilon_f\rangle^+_{D,\beta}-\langle\varepsilon_e;\varepsilon_f\rangle_{\Z^2,\beta}=\tfrac4{\pp^2}\bigl[F(G_D,1)-F(\Z^2)\bigr],\qquad\tfrac4{\pp^2}\le36 .
\]
We split the edges $f$ with $|f-e|\le r$ according to the value of $\rho$.

\emph{Step 1: the case $\rho=0$.} There are at most two such edges $f$, including $f=e$, and each of
them contributes at most $1$.

\emph{Step 2: the case $1\le\rho\le s/4$.} We use Proposition~\ref{prop:stab}(ii) with $G=G_D$ and
$\xi=1$. Since at most $16\rho$ edges $f$ satisfy $|f-e|_\infty=\rho$, the contribution of these
edges is at most
\[
C\sum_{\rho=1}^{s}\rho\cdot\rho^{-2}(\rho/s)^\theta\le C(1+\theta^{-1}).
\]

\emph{Step 3: the case $s/4<\rho\le r+1$.} Since $\rho\le r+1$, we have $\lceil\frac54\rho\rceil\le
s$. Hence, Proposition~\ref{prop:apriori} applies both to $G_D$ and to $\phi_\pp$, and it gives
\[
0\le F(G_D,1),F(\Z^2)\le C_{\rm ap}\rho^{-2}.
\]
Since $s\ge r+1$, the contribution of these edges is at most
\[
C\sum_{s/4<\rho\le s}\rho\cdot\rho^{-2}\le C(1+\log4).
\]

Summing the bounds of Steps~1--3 gives \eqref{eq:LSS}.
\end{proof}

We are now ready to prove \tagref{E2}{(E2)}, the estimate in which the bulk logarithm cancels. Given
$\gamma=g$, the components of the slit domain are independent and carry constant boundary spins
(\S\ref{ss:interface}), so that all energy covariances are nonnegative. A lower bound may therefore
discard, for each edge $e$, the terms $\langle\varepsilon_e;\varepsilon_f\rangle$ with $f$ beyond a
third of the distance from $e$ to the pinned vertices (defined below), and Theorem~\ref{thm:LSS}
compares what remains with $\chi_\beta(e;\cdot)$ at that radius. The resulting deficit relative to
$X_\beta$ is the increment of $\chi_\beta(e;\cdot)$ between the two radii, which is at most
$C(1+\log_+[d_\Omega(e)/(1\vee d_g(e))])$ by \eqref{eq:chi}. In other words, the bulk logarithm
$\log d_\Omega$ never appears, and summing over the edges gives a deficit of at most
$Cn^2\Pi_\delta(g)$ (Figure~\ref{fig:cancel}).

Recall from \S\ref{ss:notation} that $d_\Omega(e)$ and $d_g(e)$ denote the Euclidean distances from
$\mathrm{mid}(e)$ to $\Z^2\setminus\Omega_n$ and to $g$, and that
$r_\Omega(e)=1\vee\frac13d_\Omega(e)$; from \S\ref{ss:reduction} that $X_\beta=\sum_{e\in
E}\chi_\beta(e;r_\Omega(e))$; and from \eqref{eq:Lambda} that $\Pi_\delta(g)=1+\Phi(g)$. Recall also
from \S\ref{ss:interface} that $\Insp(g)$ is the set of inspected vertices, which depends only on
$g$, and that $\cD(g)$ is the set of nearest-neighbor components of $\Omega_n\setminus\Insp(g)$;
given $\gamma=g$, each $D\in\cD(g)$ has the constant boundary spin $\epsilon_D$. For an admissible
$g$, let $\mathrm{Pin}(g):=\Insp(g)\cup\partial\Omega_n$ be the set of \emph{pinned vertices}, that
is, of the vertices whose spins are fixed given $\gamma=g$, and let
\[
d^{\rm pin}_g(e):=\dist_2(\mathrm{mid}(e),\mathrm{Pin}(g)).
\]
Then $d^{\rm pin}_g(e)\le d_\Omega(e)$, since every vertex of $\Z^2\setminus\Omega_n$ nearest to
$\mathrm{mid}(e)$ is an outer boundary vertex. Moreover,
\begin{equation}\label{eq:dg}
d^{\rm pin}_g(e)\ge\min\bigl(d_g(e)-2,\ d_\Omega(e)\bigr),
\end{equation}
since every inspected vertex is within distance $2$ of $g$ (\S\ref{ss:interface}).

\begin{proposition}[lower bound on the conditional variance; \tagref{E2}{(E2)}]\label{prop:Bminus}
There is $C_{\Pi}=C_{\Pi}(M)$ such that, for all $n\ge n_0(M)$, $\beta\in W_M$ and admissible $g$,
\[
\Var_\beta(H\mid\gamma=g)\ \ge\ X_\beta-C_{\Pi}n^2\Pi_\delta(g).
\]
\end{proposition}
\begin{proof}
Let $r_0:=30$, $r_e:=\lfloor d^{\rm pin}_g(e)/3\rfloor$ and $E_{\rm far}:=\{e\in E:r_e\ge r_0\}$. We
proceed in three steps.

\emph{Step 1: Theorem~\ref{thm:LSS} applies to the edges of $E_{\rm far}$.} Let $e\in E_{\rm far}$.
No vertex within distance $2r_e<d^{\rm pin}_g(e)$ of $\mathrm{mid}(e)$ is pinned, and every such
vertex lies in $\Omega_n$, since every vertex outside $\Omega_n$ nearest to $\mathrm{mid}(e)$ is an
outer boundary vertex. These vertices form a lattice disk, which is nearest-neighbor connected and
contains both endpoints of $e$. Hence, they lie in the component $D_e\in\cD(g)$ that contains the
endpoints of $e$, and Theorem~\ref{thm:LSS} applies to $(D_e,e,r_e)$, after a global spin flip if
$\epsilon_{D_e}=-1$.

\emph{Step 2: discarding distant pairs.} By \S\ref{ss:interface},
\[
\Var_\beta(H\mid\gamma=g)=\sum_{D\in\cD(g)}\ \sum_{e,f\in E(G_D)}\langle\varepsilon_e;\varepsilon_f\rangle^{\epsilon_D}_{D,\beta},
\]
where $E(G_D)$ is the set of edges with an endpoint in $D$ (\S\ref{ss:FK}), and
$\langle\cdot\rangle^{\epsilon_D}_{D,\beta}$ denotes the expectation with respect to the Ising
measure on $D$ with boundary condition given by $\epsilon_D$. All the terms are nonnegative
(\S\ref{ss:interface}). For each $e\in E_{\rm far}$, we keep only the pairs $(e,f)$ with $|f-e|\le
r_e$, for which $e,f\in E(G_{D_e})$ by Step~1. By the nonnegativity of the discarded terms and
Theorem~\ref{thm:LSS}, we obtain
\[
\Var_\beta(H\mid\gamma=g)\ge\sum_{e\in E_{\rm far}}\bigl[\chi_\beta(e;r_e)-C_{\rm LS}\bigr].
\]
Since $X_\beta$ is the sum of $\chi_\beta(e;r_\Omega(e))$ over $e\in E$, it follows that
\[
X_\beta-\Var_\beta(H\mid\gamma=g)\le\sum_{e\in E_{\rm far}}\bigl[\chi_\beta(e;r_\Omega(e))-\chi_\beta(e;r_e)+C_{\rm LS}\bigr]+\sum_{e\in E\setminus E_{\rm far}}\chi_\beta(e;r_\Omega(e)).
\]

\emph{Step 3: the contribution of each edge.} All the radii involved are at most $n\le\frac12N$, so
that \eqref{eq:chi} applies.

\emph{Edges in $E_{\rm far}$.} If $r_e\ge r_\Omega(e)$, the term
\[
\chi_\beta(e;r_\Omega(e))-\chi_\beta(e;r_e)+C_{\rm LS}
\]
is at most $C_{\rm LS}$, since $\chi_\beta(e;\cdot)$ is nondecreasing. Otherwise, \eqref{eq:chi}
bounds this term by $C(1+\log(r_\Omega(e)/r_e))+C_{\rm LS}$. Here, $r_\Omega(e)\le d_\Omega(e)$, and
$r_e\ge d^{\rm pin}_g(e)/4$ because $d^{\rm pin}_g(e)\ge3r_0=90$. If $d^{\rm pin}_g(e)=d_\Omega(e)$,
then $r_\Omega(e)/r_e\le4$. Otherwise, $d^{\rm pin}_g(e)\ge d_g(e)-2$ by \eqref{eq:dg}, so that
$1\vee d_g(e)\le2d^{\rm pin}_g(e)$ and
\[
r_\Omega(e)/r_e\le8d_\Omega(e)/(1\vee d_g(e)).
\]

\emph{Edges not in $E_{\rm far}$.} These include the crossed edges and the edges with a pinned
endpoint. For them, $d^{\rm pin}_g(e)<90$, and
\[
\chi_\beta(e;r_\Omega(e))\le C(1+\log r_\Omega(e))
\]
by \eqref{eq:chi}. If $d^{\rm pin}_g(e)=d_\Omega(e)$, then $r_\Omega(e)\le30$. Otherwise,
$d_g(e)<92$ by \eqref{eq:dg}, and therefore
\[
\log r_\Omega(e)\le\log92+\log_+\frac{d_\Omega(e)}{1\vee d_g(e)}.
\]

In all cases, the contribution of $e$ is at most $C'(1+\log_+[d_\Omega(e)/(1\vee d_g(e))])$ for a
constant $C'=C'(M)$. Summing over $e\in E$ and using $|E|\le3n^2$ and the definition of $\Phi$ in
\eqref{eq:Lambda}, we obtain
\[
X_\beta-\Var_\beta(H\mid\gamma=g)\le3C'n^2+C'n^2\Phi(g)\le C_{\Pi}n^2\Pi_\delta(g)
\]
with $C_{\Pi}:=3C'$.
\end{proof}
\section{The specific heat with plus boundary conditions: Proof of~(E3)}\label{sec:VP} In this
section, we prove \tagref{E3}{(E3)} (Theorem~\ref{thm:VP}), that is, the upper bound in the
comparison of $\Var^+_\beta(H)$ with the bulk reference $X_\beta=\sum_e\chi_\beta(e;r_\Omega(e))$ in
the square with plus boundary conditions. We call a pair of edges $(e,f)$ \emph{near} if $|f-e|\le
r_\Omega(e)=1\vee\frac13d_\Omega(e)$, that is, if $f$ lies within a third of the distance from $e$
to the outer boundary, and \emph{far} if $|f-e|>r_\Omega(e)$. Near pairs are handled by
Theorem~\ref{thm:LSS}. The main difficulty lies in the analysis of the far pairs, in particular of
those with an edge close to the outer boundary $\partial\Omega_n$. We will overcome it by coupling
FK-Ising measures with different boundary conditions in boxes around such an edge, as in
Lemma~\ref{lem:contraction}; near the wired boundary, closed dual half-circuits ending on
$\partial\Omega_n$ take the place of the closed dual circuits (\S\ref{ss:geo}).

\begin{theorem}\label{thm:VP}
For every $M<\infty$, there is $C_V=C_V(M)$ such that $|\Var^+_\beta(H)-X_\beta|\le C_Vn^2$ for all
$n\ge n_0(M)$ and $\beta\in W_M$.
\end{theorem}
Only the upper bound is used in the proof of Theorem~\ref{thm:main}. Throughout this section, we let
$\beta\in W_M$, $n\ge n_0(M)$, $\pp=\pp(\beta)$ and $\phi:=\phi^1_\beta$ be the FK measure on
$G=(\Omega_n\cup\partial\Omega_n,E)$ with $\partial\Omega_n$ wired. Constants depend only on $M$.
Recall from \eqref{eq:N} that $N=4n$, and, as in Section~\ref{sec:cov}, we write $C_\Delta$ for the
constant $C_\Delta(8M)$ of Lemma~\ref{lem:Delta}(iii). By Lemma~\ref{lem:ES} with $D=\Omega_n$ and
by \tagref{FKG}{(FKG)}, writing
$\langle\cdot\rangle^+_\beta:=\langle\cdot\rangle^+_{\Omega_n,\beta}$, we have, for $e\ne f\in E$,
\begin{equation}\label{eq:ESplus}
\langle\varepsilon_e;\varepsilon_f\rangle^+_\beta=\tfrac4{\pp^2}\Cov_\phi(\omega_e,\omega_f)\ge0 ,
\end{equation}
including the edges joining $\Omega_n$ to the outer boundary, for which $\varepsilon_e=\sigma_x$,
with $x$ the endpoint of $e$ in $\Omega_n$.

\subsection{Near pairs and reduction to far pairs}\label{ss:near} In this subsection, we dispose of
the near pairs, and we reduce Theorem~\ref{thm:VP} to a bound on the covariances of the far pairs.

\begin{proposition}[near part]\label{prop:near}
For every $e\in E$,
\[
\Big|\sum_{f\in E:\,|f-e|\le r_\Omega(e)}\langle\varepsilon_e;\varepsilon_f\rangle^+_\beta-\chi_\beta(e;r_\Omega(e))\Big|\le C.
\]
\end{proposition}
\begin{proof}
If $d_\Omega(e)\ge3$, then $r_\Omega(e)=d_\Omega(e)/3$, and every vertex within distance
$2r_\Omega(e)<d_\Omega(e)$ of $\mathrm{mid}(e)$ lies in $\Omega_n$. Hence, Theorem~\ref{thm:LSS}
applies with $D=\Omega_n$ and $r=r_\Omega(e)$, and \eqref{eq:LSS} is the claim. If $d_\Omega(e)<3$,
then $r_\Omega(e)=1$, and both sums have $O(1)$ terms, each of which lies in $[0,1]$.
\end{proof}
We split $\Var^+_\beta(H)=\sum_{e,f\in E}\langle\varepsilon_e;\varepsilon_f\rangle^+_\beta$ into the
near part, with $|f-e|\le r_\Omega(e)$, and the far part. Since the far terms are nonnegative by
\eqref{eq:ESplus}, Proposition~\ref{prop:near} and $|E|\le3n^2$ give the lower bound of
Theorem~\ref{thm:VP}. For the upper bound, Proposition~\ref{prop:near} bounds the near part of the
sum over $f$ by $\chi_\beta(e;r_\Omega(e))+C$, and, by \eqref{eq:ESplus} and $\frac4{\pp^2}\le36$,
each far term is at most $36\Cov_\phi(\omega_e,\omega_f)$. Summing over $e\in E$, we obtain
\begin{equation}\label{eq:upred}
\Var^+_\beta(H)\le X_\beta+C|E|+36\,\mathcal F_\beta,\qquad\mathcal F_\beta:=\sum_{e\in E}\ \sum_{f\in E:\,|f-e|>r_\Omega(e)}\Cov_\phi(\omega_e,\omega_f).
\end{equation}
It therefore remains to show that $\mathcal F_\beta\le Cn^2$.

\subsection{Edges near the wired boundary}\label{ss:geo} In this subsection, we prove the analogues
of Lemmas~\ref{lem:sandwich} and~\ref{lem:edgeevent} for an edge $e$ near the wired boundary.

\begin{figure}[ht]\centering
\begin{tikzpicture}[scale=0.5]
\fill[gray!20,even odd rule] (-4,0) rectangle (4,4) (-2,0) rectangle (2,2);
\foreach \x in {-5,...,5} \draw[gray!60] (\x,0)--(\x,6);
\foreach \y in {1,...,6} \draw[gray!60] (-5,\y)--(5,\y);
\draw[line width=1.6pt] (-6.3,0)--(6.3,0);
\foreach \x in {-5,...,5} \fill (\x-0.14,-0.14) rectangle (\x+0.14,0.14);
\foreach \y in {1,...,6} {\filldraw[fill=white] (-5,\y) circle (3pt); \filldraw[fill=white] (5,\y) circle (3pt);}
\foreach \x in {-4,...,4} \filldraw[fill=white] (\x,6) circle (3pt);
\draw[line width=1.3pt,dashed] (-3.5,0)--(-3.5,2.5)--(-2.5,2.5)--(-2.5,3.5)--(2.5,3.5)--(2.5,2.5)--(3.5,2.5)--(3.5,0);
\draw[line width=2.2pt] (0,1)--(1,1); \fill (0,1) circle (3pt);
\node[above] at (0.55,1.05) {\small $e$};
\node[above left,inner sep=1pt] at (-0.05,1.05) {\small $x_e$};
\node[below] at (0,-0.2) {\small $y$};
\draw[<->,thin] (-0.6,0.12)--(-0.6,0.88); \node[left,fill=white,inner sep=1pt] at (-0.65,0.5) {\small $a_e$};
\node[above right] at (5,6) {\small $\Lambda_R(x_e)$};
\node[right] at (6.3,0) {\small $\partial\Omega_n$};
\end{tikzpicture}
\caption{An edge $e$ near the bottom side of the square, at $\ell^\infty$-distance $a_e$ from
$\Z^2\setminus\Omega_n$, and the edges of $E_e(R)$ (gray). Under $\xi_0$, the outer boundary
vertices of $E_e(R)$ (black squares) form one block, as under $\phi$ (thick line), and the other
vertices of $\partial E_e(R)$ (white circles), which lie on $\partial\Lambda_R(x_e)$, are
singletons. A closed dual half-circuit (dashed) in the half-annulus $\mathrm{Ann}^+_t(y)$ (shaded)
around the outer boundary vertex $y$ below $x_e$ separates $\Lambda_t(y)$, which contains $e$, from
$\partial\Lambda_R(x_e)$; given such a half-circuit, the configuration inside it has the same law
under $\xi_0$ as under the wired boundary condition.}\label{fig:wired}
\end{figure}
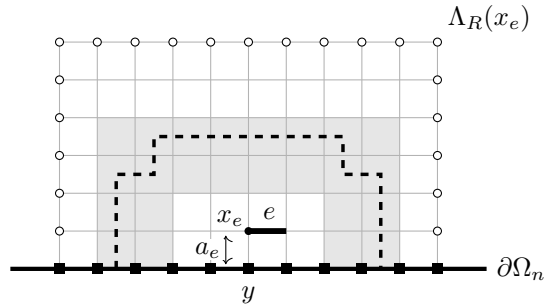

Recall from \S\ref{ss:notation} that $a_e=\dist_\infty(x_e,\Z^2\setminus\Omega_n)$ is the
$\ell^\infty$-distance from $x_e$ to the complement of $\Omega_n$, and that $\hat a_e=1\vee a_e$.
For integers $R\ge1$, let $E_e(R)$ be the set of edges of $E$ with both endpoints in
$\Lambda_R(x_e)$, viewed as a graph whose vertices are the endpoints of these edges, and let
$\partial E_e(R)$ be the set of its boundary vertices, that is, of its vertices incident to an edge
of $\Z^2$ not in $E_e(R)$ (\S\ref{ss:FK}). They are of two kinds: the outer boundary vertices, that
is, the vertices in $\partial\Omega_n$, and the vertices in $\Omega_n$ incident to an edge of
$E\setminus E_e(R)$, which lie on $\partial\Lambda_R(x_e)$. Let $\xi_0$ be the boundary condition on
$E_e(R)$ in which the outer boundary vertices form one block and every other vertex of $\partial
E_e(R)$ is a singleton. Thus, $\xi_0$ is wired on $\partial\Omega_n$, as $\phi$ is, and free
elsewhere (Figure~\ref{fig:wired}). Let
\[
J_e(R):=[\phi^{\xi_0}_{E_e(R),\pp}[\omega_e],\phi^1_{E_e(R),\pp}[\omega_e]],
\]
and we write $|J_e(R)|$ for its length. It follows from \tagref{DMP}{(DMP)}, \tagref{CBC}{(CBC)} and
the tower property that, under any FK measure that is defined on a set of edges containing $E_e(R)$
and wires the outer boundary, the conditional expectation of $\omega_e$ given the states of any set
of edges outside $E_e(R)$ takes values in $J_e(R)$; in particular, $J_e(R)\supseteq J_e(R')$ for
$R'\ge R$. If $R<a_e$, then $\Lambda_R(x_e)\subseteq\Omega_n$, so that $\xi_0$ is the free boundary
condition and $J_e(R)=I_e(R)$, with $I_e(R)$ as in Lemma~\ref{lem:sandwich}. Let $C_I:=4C_\Delta$.
Then we have
\begin{equation}\label{eq:Jelem}
|J_e(R)|\le C_I/\hat a_e\qquad\text{for }R\ge a_e/2 .
\end{equation}

\begin{lemma}[edges near the wired boundary]\label{lem:J}
Let $c_1=c_1(40M)$ be the constant of \tagref{Rb}{(Rb)} with $\lambda=40M$ (\S\ref{ss:inputs}),
which bounds from below the conditional probability of a closed dual path in a half-annulus, and let
$c_J:=-\log_2(1-c_1)\in(0,1]$. Let $e\in E$.
\begin{enumerate}[label=\textup{(\roman*)}]
\item $|J_e(R)|\le C\hat a_e^{-1}(\hat a_e/R)^{c_J}$ for every $R$ with $\hat a_e\le R\le N$.
\item Let $2\le r\le N/2$ with $4a_e<r$, and let $\zeta_e$ be the partition of $\partial E_e(r)$
induced by the FK-Ising configuration inside $E_e(r)$ and the wired boundary condition on $\partial
\Omega_n$. If $\mathcal B$ is an increasing event that is a function of $\zeta_e$, then
\[
|\Cov_\phi(\omega_e,\one_{\mathcal B})|\le C\hat a_e^{-1}(\hat a_e/r)^{c_J/2}.
\]
\end{enumerate}
\end{lemma}
\begin{proof}
Let $R_0:=4\hat a_e$, and let $y\in\Z^2\setminus\Omega_n$ with $|y-x_e|_\infty=a_e$. Here, closed
dual half-circuits, as in \tagref{Rb}{(Rb)}, have both ends on the outer boundary
(Figure~\ref{fig:wired}). Since $\phi$ is the restriction to $E$ of the wired FK measure on the box
$\{0,\dots,n+1\}^2$, whose inner boundary contains $y$, \tagref{Rb}{(Rb)}, with $\lambda=40M$,
applies to $\phi$ in the half-annuli
\[
\mathrm{Ann}^+_t(y)=(\Lambda_{2t}(y)\setminus\Lambda_t(y))\cap\{0,\dots,n+1\}^2,\qquad 2\le t\le n/8{:}
\]
conditionally on the edges outside $E(\mathrm{Ann}^+_t(y))$, with probability at least $c_1$, a
closed dual half-circuit in $\mathrm{Ann}^+_t(y)$ separates $\Lambda_t(y)\cap\Omega_n$ from
$\Omega_n\setminus\Lambda_{2t}(y)$. For $R\le N$, the half-annuli $\mathrm{Ann}^+_{r_i}(y)$ with
$r_i:=2^{i+1}R_0$, $i\ge0$ and $4r_i\le\min(R,n/2)$ are disjoint; there are at least
$\log_2(R/R_0)-9$ of them; and a closed dual half-circuit in any of them separates
$\Lambda_{R_0}(x_e)$ from $\partial\Lambda_R(x_e)$.

First, we prove (i). If $R\le2^9R_0$, the claim follows from \eqref{eq:Jelem}. Otherwise, we couple
$\omega\sim\phi^{\xi_0}_{E_e(R),\pp}$ and $\omega'\sim\phi^1_{E_e(R),\pp}$, with $\omega\le\omega'$,
as in the proof of Lemma~\ref{lem:contraction}: we explore from the outside the $\omega'$-open
cluster of $\partial\Lambda_R(x_e)$ in $E_e(R)\setminus E_e(R_0)$, and we stop when the revealed
closed edges contain a closed dual circuit or half-circuit separating $\Lambda_{R_0}(x_e)$ from
$\partial\Lambda_R(x_e)$. Here, the cluster is formed by $\omega'$-open paths only, without the
boundary wiring. Since an outer boundary vertex has a single edge in $E$, no open path passes
through such a vertex; hence, open paths cannot go around a closed dual half-circuit through the
outer boundary, and such a half-circuit separates as a closed dual circuit does. If the exploration
stops, the set $U$ of the vertices inside the circuit or half-circuit contains the endpoints of $e$,
and no edge of $E(U)$ is revealed. The edges joining $U$ to the other vertices cross the circuit or
half-circuit and are therefore closed, and the outer boundary vertices in $U$ are joined to the
other vertices only through the boundary wiring. Hence, under both measures, the conditional law of
the configuration on $E(U)$ is the FK measure on $U$ whose boundary condition, like $\xi_0$, has the
outer boundary vertices in $U$ as one block and every other boundary vertex as a singleton. In
particular, the conditional expectations of $\omega_e$ and $\omega'_e$ coincide. If the exploration
does not stop, no edge of $E_e(R_0)$ is revealed, and both conditional expectations lie in
$J_e(R_0)$, whose length is at most $C_I/\hat a_e$ by \eqref{eq:Jelem}. By \tagref{DMP}{(DMP)},
$\phi^1_{E_e(R),\pp}$ is the conditional law of $\phi$ on $E_e(R)$ given that all the edges of
$E\setminus E_e(R)$ are open, so that \tagref{Rb}{(Rb)} applies to $\omega'$ in the half-annuli
$\mathrm{Ann}^+_{r_i}(y)$, which lie in $\Lambda_R(x_e)$. Since the exploration stops as soon as one
of the half-annuli $\mathrm{Ann}^+_{r_i}(y)$ contains a closed dual half-circuit, \tagref{Rb}{(Rb)},
applied successively, gives
\[
|J_e(R)|=\E[\omega'_e-\omega_e]\le(1-c_1)^{\log_2(R/R_0)-9}C_I/\hat a_e,
\]
which is the claim.

Next, we prove (ii). We follow the proof of Lemma~\ref{lem:edgeevent}. Lemma~\ref{lem:influence}
holds with $E_e(2^j)$ in place of $E(\Lambda_{2^j}(x))$ and without the condition
$G\supseteq\Lambda_{2^k}(x)$, where now $\bar\Delta_j:=|J_e(2^{j-1})|$ for $j\ge1$ and the last term
has $|J_e(2^k)|$ in place of $\Delta_\pp(2^k)$. Indeed, its proof uses only that the conditional
expectation of $\omega_e$ given the edges outside $E_e(R)$ takes values in $J_e(R)$. We apply it
with $A=\mathcal B$ and $k$ maximal with $2^k\le r/2$, and we put $s:=\lfloor a_e/(2\Kmix)\rfloor$.

\emph{Step 1: the two clusters.} Let $0\le j\le k$ and $\omega\in\mathcal B\cap U_j$. As in Step~1
of the proof of Lemma~\ref{lem:edgeevent}, closing the edges of $E_e(2^j)$ changes $\zeta_e$, so
that there are two clusters of the resulting configuration, in distinct blocks of its partition
$\zeta_e$, that meet $\partial\Lambda_{2^j}(x_e)$ and contain a vertex of $\partial\Lambda_r(x_e)$
or an outer boundary vertex; they are not both joined to the outer boundary, since the outer
boundary vertices lie in one block of $\zeta_e$.

\emph{Step 2: the two factors.} We bound $\phi[\mathcal B\cap U_j]$ by the product of two factors,
which come from disjoint regions (Figure~\ref{fig:Jb}). For the first factor, let $2^j\le s$. The
two clusters leave $\Lambda_s(x_e)$, since $\partial\Lambda_r(x_e)$ and the outer boundary vertices
lie outside it. Hence, they cross the annulus between $\Lambda_{2^j}(x_e)$ and $\Lambda_s(x_e)$, so
that, by planar duality, $\omega\in A^{x_e}_4(2^j,s)$. This event depends only on the edges of
$E(\Lambda_s(x_e))$, and it is an event in the bulk: since $\Kmix s\le a_e/2$, we have
$\Lambda_{\Kmix s}(x_e)\subseteq\Omega_n$. Therefore, by \tagref{Mb}{(Mix)} and \eqref{eq:R1use}, as
in Step~2 of the proof of Lemma~\ref{lem:edgeevent}, its probability is at most
$\frac32C_4(8M)(2^j/s)^{1+\varepsilon_1}$. If $2^j>s$, we bound this factor by $1$.

For the second factor, note that one of the two clusters is not joined to the outer boundary; it
contains a vertex of $\partial\Lambda_r(x_e)$ and no outer boundary vertex. Hence, it crosses every
half-annulus $\mathrm{Ann}^+_{r_i}(y)$, $r_i=2^{i+1}R_0$, with $r_i\ge2^{j+1}$ and
$4r_i\le\min(r,n/2)$, so that none of these half-annuli contains a closed dual half-circuit
separating $\Lambda_{r_i}(y)\cap\Omega_n$ from $\Omega_n\setminus\Lambda_{2r_i}(y)$. There are at
least $m_j:=\log_2(r/(R_0\vee2^j))-9$ of them, and they are disjoint from one another and from
$\Lambda_s(x_e)$. By \tagref{Rb}{(Rb)}, each of them contains such a half-circuit with probability
at least $c_1$, conditionally on the edges outside it. Conditioning successively on the edges
outside each half-annulus, we obtain a factor $1-c_1$ per half-annulus, hence a factor at most
$(1-c_1)^{m_j}=2^{9c_J}((R_0\vee2^j)/r)^{c_J}$, since $c_J=-\log_2(1-c_1)$. Combining the two
factors, we obtain, for $0\le j\le k$,
\[
\phi[\mathcal B\cap U_j]\le C\Bigl(\frac{R_0\vee2^j}r\Bigr)^{c_J}\times\begin{cases}(2^j/s)^{1+\varepsilon_1}&\text{if }2^j\le s,\\1&\text{if }2^j>s.\end{cases}
\]
Since $U_{-1}\subseteq U_0$, this also bounds the term $j=-1$.

\emph{Step 3: summation.} It remains to sum over $j$. Recall that $\bar\Delta_{-1}=\bar\Delta_0=1$.
For $j\ge1$, we bound $\bar\Delta_j$ by $2C_\Delta2^{-j}$ if $2^{j-1}<a_e$, since then
$J_e(2^{j-1})=I_e(2^{j-1})$; by \eqref{eq:Jelem} if $2^{j-1}\ge a_e$ and $2^j\le R_0$; and by (i) if
$2^j>R_0$. For the scales $2^j\le R_0$, we have $R_0\vee2^j=R_0=4\hat a_e$, and, as in Step~3 of the
proof of Lemma~\ref{lem:edgeevent}, the sum of $\bar\Delta_j$ times the last factor of the display
is at most $C/\hat a_e$, since $s\ge a_e/(4\Kmix)$ if $s\ge1$, and $\hat a_e<2\Kmix$ if $s=0$; this
is where the exponent $1+\varepsilon_1>1$ is used. These scales therefore contribute at most $C\hat
a_e^{-1}(\hat a_e/r)^{c_J}$ in all. For a larger scale, $R_0<2^j\le2^k$, the decay
\[
\bar\Delta_j=|J_e(2^{j-1})|\le C\hat a_e^{-1}(\hat a_e/2^j)^{c_J}
\]
given by (i) exactly compensates for the growth of the factor $(2^j/r)^{c_J}$ of the display, so
that the term of index $j$ is at most $C\hat a_e^{-1}(\hat a_e/r)^{c_J}$; there are at most
$\log_2(r/\hat a_e)$ such scales. Finally, since $2^k>r/4>a_e$, (i) bounds the last term of
\eqref{eq:influence} by $|J_e(2^k)|\le C\hat a_e^{-1}(\hat a_e/r)^{c_J}$. Altogether,
\[
|\Cov_\phi(\omega_e,\one_{\mathcal B})|\le C\hat a_e^{-1}(\hat a_e/r)^{c_J}\bigl(1+\log(r/\hat a_e)\bigr),
\]
and the claim follows from $1+\log(r/\hat a_e)\le C(r/\hat a_e)^{c_J/2}$.
\end{proof}
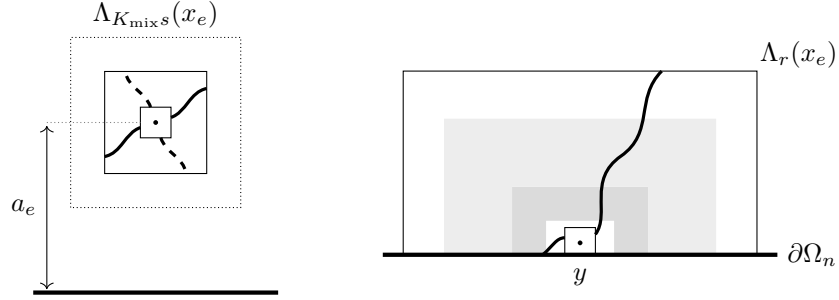
\begin{figure}[ht]\centering
\begin{tikzpicture}[scale=0.45]
\draw[line width=1.6pt] (-3.6,0)--(3.6,0);
\draw[densely dotted] (-2.5,2.5) rectangle (2.5,7.5);
\node[above] at (0,7.5) {\small $\Lambda_{\Kmix s}(x_e)$};
\draw (-1.5,3.5) rectangle (1.5,6.5);
\draw (-0.45,4.55) rectangle (0.45,5.45);
\draw[line width=1.2pt] (0.45,5.15) .. controls (0.9,5.3) and (1.0,5.9) .. (1.5,6.0);
\draw[line width=1.2pt] (-0.45,4.85) .. controls (-0.9,4.7) and (-1.0,4.1) .. (-1.5,4.0);
\draw[line width=1.2pt,dashed] (-0.15,5.45) .. controls (-0.3,5.9) and (-0.8,6.0) .. (-0.9,6.5);
\draw[line width=1.2pt,dashed] (0.15,4.55) .. controls (0.3,4.1) and (0.8,4.0) .. (0.9,3.5);
\fill (0,5) circle (2pt);
\draw[gray,thin,densely dotted] (-3.2,5)--(-0.45,5);
\draw[<->,thin] (-3.2,0.1)--(-3.2,4.9); \node[left] at (-3.2,2.5) {\small $a_e$};
\end{tikzpicture}\hspace{1.2cm}
\begin{tikzpicture}[scale=0.45]
\fill[gray!28,even odd rule] (-2,0) rectangle (2,2) (-1,0) rectangle (1,1);
\fill[gray!14,even odd rule] (-4,0) rectangle (4,4) (-2,0) rectangle (2,2);
\draw[line width=1.6pt] (-5.8,0)--(5.8,0);
\node[right] at (5.8,0) {\small $\partial\Omega_n$};
\draw (-5.2,0)--(-5.2,5.4)--(5.2,5.4)--(5.2,0);
\node[above right,inner sep=1pt] at (5.2,5.4) {\small $\Lambda_r(x_e)$};
\draw (-0.45,0)--(-0.45,0.8)--(0.45,0.8)--(0.45,0);
\fill (0,0.35) circle (2pt);
\node[below] at (0,-0.1) {\small $y$};
\draw[line width=1.2pt] (0.45,0.6) .. controls (0.9,1.2) and (0.2,2.2) .. (1.2,2.9) .. controls (2.2,3.6) and (1.6,4.6) .. (2.4,5.4);
\draw[line width=1.2pt] (-0.45,0.5) .. controls (-0.7,0.5) and (-0.8,0.25) .. (-1.1,0);
\end{tikzpicture}
\caption{The two factors in the proof of Lemma~\ref{lem:J}(ii), not to scale; the dot is $x_e$.
Left: for $2^j\le s$, the box $\Lambda_{\Kmix s}(x_e)$ (dotted) lies in $\Omega_n$, since $\Kmix
s\le a_e/2$. Hence, the two clusters given by Step~1 of that proof produce a four-arm event in the
bulk: open arms (thick) and dual-open arms (dashed) cross the annulus between the two solid squares
$\partial\Lambda_{2^j}(x_e)$ and $\partial\Lambda_s(x_e)$. Right: the two clusters leave the small
box $\Lambda_{2^j}(x_e)$; one of them is not joined to $\partial\Omega_n$ and reaches
$\partial\Lambda_r(x_e)$ (long thick path), so that, at the scales from $R_0\vee2^j$ to $r$, no
half-annulus $\mathrm{Ann}^+_{r_i}(y)$ (shaded) contains a closed dual half-circuit. No arm event is
used at these scales.}\label{fig:Jb}
\end{figure}

\subsection{The far pairs}\label{ss:far} In this subsection, we bound the covariances of all pairs
of edges and sum them over the far pairs. Besides \tagref{Rb}{(Rb)} and the fact that every outer
boundary vertex has a single edge in $E$, we use two properties of the square, which are the only
other features of the domain that this section uses.
\begin{enumerate}[leftmargin=3.2em]
\item[\tagdef{Geo}{(Geo)}] $\#\{e\in E:a_e\le s\}\le24ns$ for $s\ge1$.
\item[\tagdef{Geop}{(Geo$'$)}] $\#\{e\in E:|x_e-z|_\infty\le R,\ a_e\le B\}\le48RB$ for all $z$ and
$1\le B\le R$.
\end{enumerate}
Both hold because the boundary of the square consists of four segments of length $n+1$. In the next
proposition, the weight $w_{e'}$ is $1/\rho$ when $\hat a_{e'}\ge\rho$, and $\hat a_{e'}^{-1}(\hat
a_{e'}/\rho)^{c_J/2}$, the form of the bounds of Lemma~\ref{lem:J} at radius $\asymp\rho$,
otherwise.

\begin{proposition}[far pairs]\label{prop:far}
For all $e\ne f\in E$, with $\rho:=1\vee|e-f|_\infty$ and $c_J$ as in Lemma~\ref{lem:J},
\[
0\le\Cov_\phi(\omega_e,\omega_f)\le Cw_ew_f,\quad  \text{where }w_{e'}:=\frac1{\hat a_{e'}\wedge\rho}\Bigl(1\wedge\frac{\hat a_{e'}}\rho\Bigr)^{c_J/2}\text{ for }e'\in\{e,f\}.
\]
Consequently, we have
\[\sum_{f\in E:\,|f-e|>r_\Omega(e)}\Cov_\phi(\omega_e,\omega_f)\le C\bigl(1+(N/\hat a_e)^{1-c_J/2}\bigr)\]
for every $e\in E$, which implies that $\mathcal F_\beta\le Cn^2$.
\end{proposition}
\begin{proof}
Nonnegativity follows from \tagref{FKG}{(FKG)}, and, since $w_{e'}\ge1/\rho$, we may assume that
$\rho\ge64$. We follow the proof of Proposition~\ref{prop:apriori}, with $r:=\lfloor\rho/3\rfloor$.
Let $U_e:=E(\Lambda_{\lfloor r/16\rfloor}(x_e))$ if $a_e\ge r/4$ and $U_e:=E_e(r)$ otherwise, and
let $\Psi:=\phi[\omega_f\mid\omega|_{U_e}]$, so that $\Cov(\omega_e,\omega_f)=\Cov(\omega_e,\Psi)$.
By \tagref{DMP}{(DMP)} and \tagref{CBC}{(CBC)}, $\Psi$ is an increasing function of the boundary
partition of $U_e$, with the vertices of $\partial\Omega_n$ identified, and it takes values in
$J_f(r)$, since $U_e\cap E_f(r)=\emptyset$. If $a_f\ge r/4$, it takes values in $I_f(\lfloor
r/8\rfloor)$, as in Lemma~\ref{lem:sandwich}, whose length is at most $C/r\le Cw_f$; otherwise,
\[
|J_f(r)|\le C\hat a_f^{-1}(\hat a_f/r)^{c_J}\le Cw_f
\]
by Lemma~\ref{lem:J}(i). As in Step~2 of the proof of Proposition~\ref{prop:apriori},
$\Cov(\omega_e,\Psi)$ is the integral, over an interval of this length, of covariances of $\omega_e$
with increasing events of the boundary partition of $U_e$. These covariances are bounded by $Cw_e$,
by Lemma~\ref{lem:edgeevent} with radius $\lfloor r/16\rfloor$ if $a_e\ge r/4$, and by
Lemma~\ref{lem:J}(ii) otherwise. This proves the first claim.

For the second claim, note that every far pair has $|f-e|_\infty\ge\hat a_e/(3\sqrt2)-2$, since
$|f-e|\le\sqrt2|f-e|_\infty+1$, $d_\Omega(e)\ge a_e-\frac12$ and $r_\Omega(e)\ge\frac13d_\Omega(e)$.
In a dyadic shell $\{f:\rho\in[R,2R)\}$, we have $w_e\le(\hat a_e\wedge R)^{-1}(1\wedge\hat
a_e/R)^{c_J/2}$; at most $CR^2$ edges $f$ have $\hat a_f\ge R$, and they have $w_f\le R^{-1}$; and,
for dyadic $B<R$, at most $CRB$ edges $f$ have $\hat a_f\in[B,2B)$, by \tagref{Geop}{(Geo$'$)}, and
they have $w_f\le B^{-1}(2B/R)^{c_J/2}$. Hence, the shell contributes at most $Cw_eR$. The $O(1)$
shells with $R\le\hat a_e$ that contain far pairs contribute $O(1)$, and, since $c_J\le1$, the
shells with $\hat a_e<R\le N$ contribute at most
\[
C\sum_R(R/\hat a_e)^{1-c_J/2}\le C(N/\hat a_e)^{1-c_J/2}.
\]
Finally, by \tagref{Geo}{(Geo)}, $\#\{e:\hat a_e<2A\}\le CnA$ for $A\ge1$, so that, grouping the
edges according to the dyadic scale of $\hat a_e$,
\[
\sum_{e\in E}\Bigl(\frac N{\hat a_e}\Bigr)^{1-c_J/2}\le\sum_{A\text{ dyadic},\ A\le N}CnA\Bigl(\frac NA\Bigr)^{1-c_J/2}\le CnN\le Cn^2 .
\]
Together with $|E|\le3n^2$, this gives $\mathcal F_\beta\le Cn^2$.
\end{proof}

\begin{proof}[Proof of Theorem~\ref{thm:VP}]
The lower bound was proved after Proposition~\ref{prop:near}. The upper bound follows from
\eqref{eq:upred}, Proposition~\ref{prop:far} and $|E|\le3n^2$.
\end{proof}

\section{The boundary disorder correlation: Proof of~(E4)}\label{sec:mumu} In this section, we prove
\tagref{E4}{(E4)}. Recall from \S\ref{ss:Dob} that the boundary disorder correlation
$\phi^1_\beta(\Dab)=e^{\ell(\beta)}$ is the probability, under the wired FK measure, that no open
path joins the two boundary arcs. It is the only nonlocal term of the decomposition
\eqref{eq:reduction}. By Proposition~\ref{prop:split}, it enters only through the ratio
$\phi^1_\beta(\Dab)/\phi^1_{\beta_c}(\Dab)$, and it suffices to bound this ratio from above and
below by constants. We compare the two quantities in three steps. First, Kramers--Wannier duality
turns $\phi^1_\beta(\Dab)$ into the probability that the two marked boundary points of the dual box
are connected, at the dual parameter $\pp(\beta)^*$ (Lemma~\ref{lem:dual}). Second, near each of its
two ends, this connection is a one-arm event at a free boundary, so that, at a single parameter, a
bounded number of RSW gluings bound the probability from above and below by constant multiples of
squares of connection probabilities $A_h$ in strips of height $h\asymp n$ (defined in
\S\ref{ss:strip}; see \S\ref{ss:glue}). Third, in the strip, an identity based on the massive
observable (cf.\ \cite{BDC12,DGP14}) relates $A_h(\pp)$ and $A_h(\pp^*)$ \emph{exactly} across
parameters: their ratio is $\mathfrak r(\pp)^h$, where $\mathfrak r(\pp)$ is the explicit factor by
which the modulus of the observable changes from one row of the strip to the next \cite[\S3,
Step~1]{BDC12}, and this ratio is bounded away from $0$ and $\infty$ for the parameters that we use
(\S\ref{ss:DGP}). Recall that $W_M=[\beta_c-M\delta,\beta_c+M\delta]$.

\begin{theorem}\label{thm:mumu}
For every $M<\infty$, there is $C(M)<\infty$ such that
\[
C(M)^{-1}\le\phi^1_\beta(\Dab)/\phi^1_{\beta_c}(\Dab)\le C(M)
\]
for all $n\ge n_0(M)$ and $\beta\in W_M$. In particular, \eqref{eq:Smu} holds with $C_{\rm dis}\le
C(M)$.
\end{theorem}
Since $\ell$ is nonincreasing (Lemma~\ref{lem:dob}(iii)), on each side of $\beta_c$ one of the two
inequalities is immediate. The content of Theorem~\ref{thm:mumu} is therefore the upper bound for
$\beta<\beta_c$, which is used in the proof of Theorem~\ref{thm:main}, and the lower bound for
$\beta>\beta_c$, which is used only for $m<0$ (see the discussion at the end of
\S\ref{ss:strategy}).

\subsection{Duality and the strip}\label{ss:strip} In this subsection, we turn the boundary disorder
correlation into a connection probability in the dual box, and we introduce the strip in which the
comparison across parameters will be made in \S\ref{ss:DGP}. Let $\Omega^*_n$ be the box
$\{\frac12,\frac32,\dots,n+\frac12\}^2$ with all its nearest-neighbor edges, and let $x_1^*,x_2^*$
be its vertices $(\lfloor n/2\rfloor+\frac12,\frac12)$ and $(\lfloor n/2\rfloor+\frac12,n+\frac12)$,
which are the centers of the faces at which the boundary condition changes sign. Let
$\pp\mapsto\pp^*$ be the involution defined by
\[\frac{\pp\pp^*}{(1-\pp)(1-\pp^*)} =2.\]

\begin{lemma}[Kramers--Wannier duality; \cite{KW41}]\label{lem:dual}
For every $\beta>0$,
\[
\phi^1_\beta(\Dab)=\phi^0_{\Omega^*_n,\pp(\beta)^*}(x_1^*\leftrightarrow x_2^*).
\]
\end{lemma}
Lemma~\ref{lem:dual} is a standard application of the planar duality of random-cluster measures
(see, e.g., \cite[\S6.1, (6.6)]{Gri06}), after the identification of $\partial\Omega_n$ with a
single vertex; we omit the proof.

\emph{The parameter window and the RSW bound.} We put $\lambda:=3M$ and
\[
\mathcal W_\lambda(n):=[\ppc-\lambda/n,\ppc+\lambda/n].
\]
There is $n_1(M)\ge n_0(M)$ such that, for $n\ge n_1(M)$ and $\beta\in W_M$, both $\pp(\beta)$ and
$\pp(\beta)^*$ lie in $\mathcal W_\lambda(n)$, and we will enlarge $n_1(M)$ below. For $n_0(M)\le
n<n_1(M)$, Theorem~\ref{thm:mumu} is immediate, since $\beta\mapsto\phi^1_\beta(\Dab)$ is positive
and continuous on the compact set $W_M\subset(0,\infty)$. If $n\ge7\lambda$, then
\[
\mathcal W_\lambda(n)\subseteq[\ppc-\frac17,\ppc+\frac17]\subseteq[\frac13,\frac34].
\]
We take $n_1(M)\ge7\lambda$, and the statements of this section involving $\mathcal W_\lambda(n)$
are made for $n\ge7\lambda$ only. We use Theorem~\ref{thm:RSW}, with $2\lambda$ in place of
$\lambda$, in the following form. For $\varrho\ge1$, there is $c_{\rm RSW}(\lambda,\varrho)>0$ such
that every rectangle with sides $r,r'\le2n$ and $r'\in[r/\varrho,\varrho r]$ is crossed in either
direction with conditional probability in $[c_{\rm RSW}(\lambda,\varrho),1-c_{\rm
RSW}(\lambda,\varrho)]$ given any configuration outside it, for every $n\ge7\lambda$ and every
$\pp\in\mathcal W_\lambda(n)$, inside any domain.

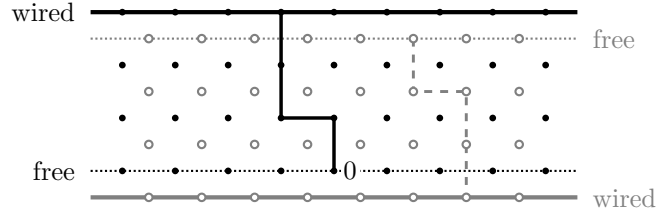
\begin{figure}[ht]\centering
\begin{tikzpicture}[scale=0.7]
\draw[gray,densely dotted,thick] (-4.6,2.5)--(4.6,2.5); \draw[gray,line width=1.6pt] (-4.6,-0.5)--(4.6,-0.5);
\draw[densely dotted,thick] (-4.6,0)--(4.6,0); \draw[line width=1.6pt] (-4.6,3)--(4.6,3);
\draw[line width=1.3pt] (0,0)--(0,1)--(-1,1)--(-1,3);
\draw[gray,line width=1.1pt,dashed] (1.5,2.5)--(1.5,1.5)--(2.5,1.5)--(2.5,-0.5);
\node[right,fill=white,inner sep=1pt] at (0.12,0) {\small $0$};
\foreach \x in {-4,...,3} \foreach \y in {-1,...,2} \filldraw[gray,thick,fill=white] (\x+0.5,\y+0.5) circle (2pt);
\foreach \x in {-4,...,4} \foreach \y in {0,...,3} \fill (\x,\y) circle (1.8pt);
\node[left] at (-4.7,3) {\small wired}; \node[left] at (-4.7,0) {\small free};
\node[right] at (4.7,2.5) {\small\color{gray}free}; \node[right] at (4.7,-0.5) {\small\color{gray}wired};
\end{tikzpicture}
\caption{The strip $S_h$ ($h=3$), with its vertices (black dots), free on the bottom row (dotted)
and wired on the top row (thick); $A_h(\pp)$ is the probability of an open path from $0$ to the top
row (black path). The dual vertices (gray circles) form the dual strip, of the same height, free on
its top row (gray, dotted) and wired on its bottom row (gray, thick), which lies below the bottom
row of $S_h$. After a reflection and a translation, the dual strip is $S_h$, with the same boundary
conditions and the parameter $\pp^*$; hence, the probability of a dual-open path from a vertex of
its free row to its wired row (dashed path) is $A_h(\pp^*)$.}\label{fig:strip}
\end{figure}

\emph{The strip.} Let $S_h:=\Z\times\{0,\dots,h\}$, and, for integers $a<b$, let
\[
\Rr_{a,b}:=\{a,\dots,b\}\times\{0,\dots,h\}.
\]
A boundary condition $\xi$ on $\Rr_{a,b}$ is \emph{admissible} if the top row $\{y=h\}$ lies in one
block of $\xi$ and every vertex $(x,0)$ with $a<x<b$ is a singleton of $\xi$; the blocks containing
vertices of the two sides are arbitrary. By Lemma~\ref{lem:stripmeasure}, as $a\to-\infty$ and
$b\to\infty$, the measures $\phi^\xi_{\Rr_{a,b},\pp}$ with admissible $\xi$ converge to a
probability measure $\phi^{0,1}_{S_h,\pp}$, free on the bottom row and wired on the top row
(Figure~\ref{fig:strip}). The convergence holds on every event depending on finitely many edges, at
an exponential rate that is uniform in $\xi$. The limit is invariant under horizontal translations
and stochastically increasing in $\pp$, and, as a limit of FK measures, it satisfies
\tagref{FKG}{(FKG)} for events depending on finitely many edges. It is the measure constructed in
\cite[\S1]{BDC12}, with the two rows exchanged, and denoted $\phi^{-\infty,\infty}_{S_n,\pp}$ in
\cite[\S2.2]{DGP14}. We let
\[A_h(\pp):=\phi^{0,1}_{S_h,\pp}(0\leftrightarrow\text{top row}).\]

\begin{lemma}\label{lem:Amono}
The probability $A_h(\pp)$ is nondecreasing in $\pp$ and nonincreasing in $h$.
\end{lemma}
\begin{proof}
For $J\ge2$, let $\Dd_J:=\Rr_{-J,J}$, with the admissible boundary condition $\xi_J$ that has one
block
\[
\Ww_J:=\{(x,y)\in\Dd_J:|x|=J\text{ or }y=h\},
\]
consisting of the top row and both sides. By Lemma~\ref{lem:stripmeasure}(ii),
\[
A_h(\pp)=\lim_{J\to\infty}\phi^{\xi_J}_{\Dd_J,\pp}(0\leftrightarrow\Ww_J).
\]
Both monotonicities then follow by standard applications of \tagref{CBC}{(CBC)} and
\tagref{DMP}{(DMP)} to the right-hand side; we omit the details.
\end{proof}
\subsection{The strip identity}\label{ss:DGP} In this subsection, we compare $A_h(\pp)$ across
parameters. From each row of the strip to the next, the modulus of the massive observable is
multiplied by the same explicit factor, which yields an exact relation between $A_h(\pp)$ and
$A_h(\pp^*)$ at the same height. Together with the monotonicity of Lemma~\ref{lem:Amono}, this gives
the two-sided comparison of Proposition~\ref{prop:strip}. For $0<\pp<1$, we put
\[
\mathfrak r(\pp):=\frac{\pp}{(1-\pp)(2-\pp)} .
\]

\begin{theorem}[strip identity; cf.\ {\cite[\S3, Step~1]{BDC12}} and {\cite[Proof of Prop.~2.5]{DGP14}}]\label{thm:strip}
For $0<\pp\le\ppc$ and $h\ge1$,
\begin{equation}\label{eq:DGPid}
A_h(\pp)=\mathfrak r(\pp)^h\,A_h(\pp^*).
\end{equation}
\end{theorem}
Theorem~\ref{thm:strip} combines two properties of the massive fermionic observable of \cite{BDC12}
in the strip, which we will recall in Appendix~\ref{app:strip}. First, translation and reflection
invariance turn its local relations into an exact relation between consecutive rows \cite[\S3,
Step~1]{BDC12}: from each row to the next, in the direction of the free row, its modulus is
multiplied by the same factor, which is $\mathfrak r(\pp)$ written in trigonometric form. Second,
its modulus is $A_h(\pp)$ on the free row \cite[Lemma~2.2]{BDC12} and, by planar duality,
$A_h(\pp^*)$ on the wired row \cite[p.~15]{DGP14}. The identity holds with the same height on both
sides, because the dual of the strip is again a strip of the same height (Figure~\ref{fig:strip}).

\begin{proposition}[strip comparison]\label{prop:strip}
For every $\lambda>0$, there is $h_0(\lambda)$ such that, with $c_S=c_S(\lambda):=e^{-12\lambda}$,
for all $h\ge h_0(\lambda)$ and all $\pp\le\ppc$ with $h(\ppc-\pp)\le2\lambda$,
\[
c_S\,A_h(\ppc)\le A_h(\pp)\le A_h(\ppc)\le A_h(\pp^*)\le c_S^{-1}A_h(\ppc).
\]
Consequently, $c_S\le A_h(\pp)/A_h(\ppc)\le c_S^{-1}$ for all $h\ge h_0(\lambda)$ and all $\pp$ with
$h|\pp-\ppc|\le\lambda$.
\end{proposition}
\begin{proof}
We choose $h_0(\lambda)$ such that $2\lambda/h_0\le\ppc-\frac13$ and $h|\pp^*-\ppc|\le2\lambda$
whenever $h\ge h_0$ and $h|\pp-\ppc|\le\lambda$. The second condition can be met because
$\pp\mapsto\pp^*$ is smooth near $\ppc$, with derivative $-1$ there. Now let $\pp\le\ppc$ with
$h(\ppc-\pp)\le2\lambda$. Then $\pp\in[\frac13,\ppc]$ and, since the involution is decreasing and
fixes $\ppc$, $\pp^*\ge\ppc$. The middle inequalities $A_h(\pp)\le A_h(\ppc)\le A_h(\pp^*)$
therefore follow from Lemma~\ref{lem:Amono}.

For the outer inequalities, we first bound $\mathfrak r(\pp)^h$. Note that $\mathfrak r(\ppc)=1$,
since $1-\ppc=\sqrt2-1$ and $2-\ppc=\sqrt2$, so that
\[
(1-\ppc)(2-\ppc)=2-\sqrt2=\ppc.
\]
Moreover,
\[
(\log\mathfrak r)'(\pp)=\frac1\pp+\frac1{1-\pp}+\frac1{2-\pp}
\]
is positive and convex on $(0,1)$. Its value is $5.1$ at $\frac13$ and $2+2\sqrt2<5$ at $\ppc$, so,
by convexity, it is at most $6$ on $[\frac13,\ppc]$. Integrating from $\pp$ to $\ppc$, we obtain
\[
0\le-\log\mathfrak r(\pp)\le6(\ppc-\pp)\le\frac{12\lambda}h .
\]
Hence, \eqref{eq:DGPid} gives
\[
A_h(\pp)=\mathfrak r(\pp)^hA_h(\pp^*)\ge e^{-12\lambda}A_h(\pp^*).
\]
Combined with the middle inequalities, this gives both outer inequalities.

For the consequence, let $h|\pp-\ppc|\le\lambda$. If $\pp\le\ppc$, we apply the first part at $\pp$.
If $\pp>\ppc$, then $\pp^*<\ppc$ and $h(\ppc-\pp^*)\le2\lambda$ by the choice of $h_0$, so we apply
the first part at $\pp^*$, using $(\pp^*)^*=\pp$. Note that the identity \eqref{eq:DGPid} is used
only at parameters that are smaller than or equal to $\ppc$.
\end{proof}
\subsection{From the strip to the dual box: Gluing at the macroscopic scale}\label{ss:glue} In this
subsection, we bound the connection probability in the dual box from above and below by constant
multiples of squares of the connection probabilities $A_h$. All comparisons here are at a
\emph{single} parameter $\pp\in\mathcal W_\lambda(n)$, with $n\ge7\lambda$. They change only the
geometry and the boundary conditions at scale $\asymp n$, using a bounded number of RSW events,
while the comparison across parameters has been made exactly, in \S\ref{ss:DGP}, by
\eqref{eq:DGPid}.

\emph{Half-boxes and half-circuits.} Let $x$ lie on a horizontal line $\mathcal L$ that is a free
boundary of the domain, and let $r\ge1$. Let $B^+(x,r)$, the \emph{half-box} of radius $r$ around
$x$, be the set of vertices within $\ell^\infty$-distance $r$ of $x$ on the domain side of $\mathcal
L$, and let $\partial_{\rm out}B^+(x,r)$ be the set of its vertices at $\ell^\infty$-distance $r$
from $x$. We put $\pi^{00}_\pp(r):=\phi^0_{B^+(x,r),\pp}(x\leftrightarrow\partial_{\rm out})$ and
$\pi^{10}_\pp(r):=\phi^{\rm w}_{B^+(x,r),\pp}(x\leftrightarrow\partial_{\rm out})$, where $\phi^{\rm
w}$ is wired on $\partial_{\rm out}$ and free on $\mathcal L$. These quantities do not depend on
$x$, and $\pi^{00}\le\pi^{10}$. An \emph{open half-circuit around $x$} in $B^+(x,r)\setminus
B^+(x,r/2)$ is an open path there from $\mathcal L$ on the left of $x$ to $\mathcal L$ on the right
of $x$. By planarity, every path from $B^+(x,r/2)$ to $\partial_{\rm out}B^+(x,r)$ meets it.

By a standard application of RSW and \tagref{FKG}{(FKG)}, we have
\begin{equation}\label{eq:halfcircuit}
\phi^\xi_{D,\pp}\bigl(\exists\text{ an open half-circuit around $x$ in }B^+(x,r)\setminus B^+(x,r/2)\bigr)\ge c_{\rm hc}(\lambda):=c_{\rm RSW}(\lambda,8)^3
\end{equation}
for every domain $D\supseteq B^+(x,r)$ having $\mathcal L$ as its boundary near $x$, every $\xi$,
real $r\in[4,n]$, $n\ge7\lambda$ and $\pp\in\mathcal W_\lambda(n)$.

The \emph{outermost} such half-circuit $\Gamma$ is found by exploring from $\partial_{\rm out}$
inward, so $\{\Gamma=\gamma\}$ depends only on the edges on and outside $\gamma$. Given
$\Gamma=\gamma$, the conditional law in the region $\mathrm{in}(\gamma)$ enclosed by $\gamma$ and
$\mathcal L$ \emph{equals} $\phi^{\gamma\text{-w}}_{\mathrm{in}(\gamma),\pp}$, the measure wired on
$\gamma$ with the original, free, boundary condition on $\mathcal L$. Indeed, $\gamma$ is open, so
all its vertices are connected, and the claim follows from \tagref{DMP}{(DMP)}. In
Lemma~\ref{lem:bc}, this will be applied under the strip measure $\phi^{0,1}_{S_h,\pp}$, which is
not a finite-volume measure. There, the event $\{\Gamma=\gamma\}$ and the boundary condition induced
on $\mathrm{in}(\gamma)$ are determined by finitely many edges, since $\gamma$ is open and the rest
of the boundary of $\mathrm{in}(\gamma)$ is the free bottom row. Consequently, \tagref{DMP}{(DMP)}
passes from the measures $\phi^\xi_{\Rr,\pp}$ to $\phi^{0,1}_{S_h,\pp}$ by
Lemma~\ref{lem:stripmeasure}(i), and so do \tagref{FKG}{(FKG)}, RSW and \eqref{eq:halfcircuit},
which concern events depending on finitely many edges. The next lemma compares $\pi^{00}$ and
$\pi^{10}$ with $A_h$ at comparable scales.
\begin{lemma}\label{lem:bc}
There is $c_G=c_G(\lambda)>0$ such that, for $n\ge7\lambda$, integers $8\le r\le n$ and
$\pp\in\mathcal W_\lambda(n)$: \textup{(i)} $A_{2r}(\pp)\ge c_G\,\pi^{10}_\pp(r)$; \textup{(ii)}
$\pi^{00}_\pp(r)\ge c_G\,A_{\lfloor r/2\rfloor}(\pp)$.
\end{lemma}
\begin{proof}
Both bounds rest on the same gluing. A path from $0$ to an outer half-circuit, together with a long
vertical crossing that meets this half-circuit, produces the connection on the left-hand side. The
connection to the half-circuit is then estimated by conditioning on the outermost half-circuit,
given which the measure inside it is wired on the half-circuit.

We start with (i). We work under $\phi:=\phi^{0,1}_{S_{2r},\pp}$, with $\mathcal L$ the bottom row
and $x=0$. Let $\mathcal E$ be the event that some open half-circuit around $0$ in
$B^+(0,r)\setminus B^+(0,r/2)$ is connected to $0$, and let $\mathcal V$ be the event that the
rectangle $[-\lceil\frac r8\rceil,\lceil\frac r8\rceil]\times[0,2r]$ is crossed vertically. As
$r\ge8$, the width of this rectangle lies in $[\frac r4,\frac r4+2]\subset[\frac r4,2r]$, so its
aspect ratio is at most $8$, and $\lceil\frac r8\rceil\le\frac r2$. Hence, a vertical crossing
starts in $B^+(0,r/2)$ and ends on the top row, and it meets every such half-circuit. This shows
that $\mathcal E\cap\mathcal V\subseteq\{0\leftrightarrow\text{top row}\}$, and, since both events
are increasing, \tagref{FKG}{(FKG)} gives
\[
A_{2r}(\pp)\ge\phi(\mathcal V)\,\phi(\mathcal E),\qquad\text{with }\phi(\mathcal V)\ge c_{\rm RSW}(\lambda,8)\text{ by RSW}.
\]
Next, we show that
\begin{equation}\label{eq:bcE}
\phi(\mathcal E)\ \ge\ \sum_\gamma\phi(\Gamma=\gamma)\,\phi^{\gamma\text{-w}}_{\mathrm{in}(\gamma),\pp}(0\leftrightarrow\gamma)\ \ge\ c_{\rm hc}\,\pi^{10}(r),
\end{equation}
where $\Gamma$ is the outermost half-circuit. The first inequality holds because, given
$\Gamma=\gamma$, the conditional law in $\mathrm{in}(\gamma)$ is
$\phi^{\gamma\text{-w}}_{\mathrm{in}(\gamma),\pp}$, and $\mathcal E$ occurs as soon as $0$ is
connected to $\gamma$. For the second inequality, note that $\mathrm{in}(\gamma)\subseteq B^+(0,r)$,
and that every path from $0$ to $\partial_{\rm out}B^+(0,r)$ contains a path from $0$ to $\gamma$
inside $\mathrm{in}(\gamma)$. Under $\phi^{\rm w}_{B^+(0,r),\pp}$, the conditional law in
$\mathrm{in}(\gamma)$ is dominated by the $\gamma$-wired one, both being free on $\mathcal L$.
Therefore, $\phi^{\gamma\text{-w}}_{\mathrm{in}(\gamma),\pp}(0\leftrightarrow\gamma)\ge\pi^{10}(r)$
for every $\gamma$, and it remains to use that
\[
\sum_\gamma\phi(\Gamma=\gamma)=\phi(\Gamma\text{ exists})\ge c_{\rm hc}
\]
by \eqref{eq:halfcircuit}. Thus, (i) holds for every $c_G\le c_{\rm RSW}(\lambda,8)\,c_{\rm hc}$,
and the same constant suffices for (ii).

The proof of (ii) is similar, one scale further in. We now work under $\phi:=\phi^0_{B^+(0,r),\pp}$.
Let $\mathcal E'$ be the event that some open half-circuit in $B^+(0,r/2)\setminus B^+(0,r/4)$ is
connected to $0$, and let $\mathcal V'$ be the event that $[-\lceil\frac r{16}\rceil,\lceil\frac
r{16}\rceil]\times[0,r]$ is crossed vertically. As $r\ge8$, the width of this rectangle lies in
$[\frac r8,r]$ and $\lceil\frac r{16}\rceil\le\frac r4$, so the crossing starts in $B^+(0,r/4)$ and
ends on $\partial_{\rm out}B^+(0,r)$. As in (i), $\pi^{00}(r)\ge\phi(\mathcal E')\phi(\mathcal V')$
and $\phi(\mathcal V')\ge c_{\rm RSW}(\lambda,8)$, and we use the half-circuit bound
\eqref{eq:halfcircuit} at the radius $r/2\ge4$. The only new point is the lower bound, given
$\Gamma=\gamma$,
\[
\phi^{\gamma\text{-w}}_{\mathrm{in}(\gamma),\pp}(0\leftrightarrow\gamma)\ge A_{\lfloor r/2\rfloor}(\pp).
\]
To see this, note that $\mathrm{in}(\gamma)\subseteq S_{\lfloor r/2\rfloor}$, that the restriction
of $\phi^{0,1}_{S_{\lfloor r/2\rfloor},\pp}$ to $\mathrm{in}(\gamma)$ is dominated by the
$\gamma$-wired measure and that every path from $0$ to the top row of $S_{\lfloor r/2\rfloor}$
contains a path from $0$ to $\gamma$ inside $\mathrm{in}(\gamma)$, since the vertices of
$\mathrm{in}(\gamma)$ off $\gamma$ lie below this row. Summing over $\gamma$ as in \eqref{eq:bcE}
gives $\phi(\mathcal E')\ge c_{\rm hc}A_{\lfloor r/2\rfloor}(\pp)$, and (ii) follows.
\end{proof}

\begin{figure}[ht]\centering
\begin{tikzpicture}[scale=0.38]
\fill[gray!20] (7,0) rectangle (9,16);
\draw[line width=1.1pt] (0,0) rectangle (16,16);
\draw[dashed] (4,0) rectangle (12,4); \draw[dashed] (4,12) rectangle (12,16);
\draw[line width=1.4pt] plot[smooth,tension=0.8] coordinates {(4.6,0) (4.8,2.2) (6.5,3.4) (9.5,3.3) (11.2,2.4) (11.4,0)};
\draw[line width=1.4pt] plot[smooth,tension=0.8] coordinates {(4.6,16) (5,13.8) (7,12.6) (9.5,12.7) (11.1,13.6) (11.4,16)};
\draw[line width=1.4pt] plot[smooth,tension=0.6] coordinates {(7.2,0) (7.6,4) (8.4,8) (7.8,12) (8.8,16)};
\draw[thin] plot[smooth,tension=0.7] coordinates {(8,0) (9.2,1.2) (10.1,2.9) (10.8,4)};
\draw[thin] plot[smooth,tension=0.7] coordinates {(8,16) (6.8,14.7) (5.6,13.2) (5.2,12)};
\fill (8,0) circle (4pt) node[below] {\small $x_1^*$};
\fill (8,16) circle (4pt) node[above] {\small $x_2^*$};
\node at (10.35,0.9) {\scriptsize $E_{x_1}$}; \node at (5.75,15.1) {\scriptsize $E_{x_2}$};
\node[right] at (12,2) {\small $B_{x_1}$}; \node[right] at (12,14) {\small $B_{x_2}$};
\node[right] at (16,8) {\small $\Omega^*_n$};
\end{tikzpicture}
\caption{The lower bound of Lemma~\ref{lem:factor}. The thin paths realize $E_{x_1}$ and $E_{x_2}$:
they join $x_1^*$ and $x_2^*$ to the outer boundaries of the half-boxes $B_{x_1}$ and $B_{x_2}$
(dashed). Each of them crosses the open half-circuit of its half-box (thick arcs, $\mathcal O_{x_1}$
and $\mathcal O_{x_2}$), and both half-circuits cross the vertical crossing $\mathcal V$ of the
central column (thick path in the gray column). Hence, on the intersection of the five increasing
events, $x_1^*$ is connected to $x_2^*$.}\label{fig:glue}
\end{figure}
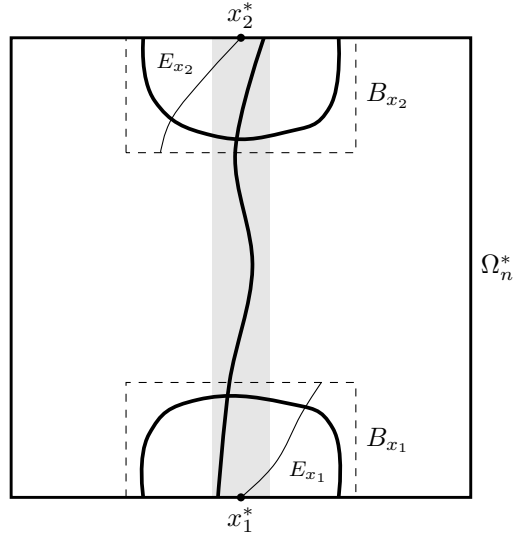

The second gluing step factorizes the connection in the dual box into two one-arm events at its free
boundary.
\begin{lemma}\label{lem:factor}
There is $c_F=c_F(\lambda)>0$ such that, for $n\ge\max(64,7\lambda)$, $\pp\in\mathcal W_\lambda(n)$
and $r=\lfloor n/4\rfloor$,
\[
c_F\,\pi^{00}_\pp(r)^2\le\phi^0_{\Omega^*_n,\pp}(x_1^*\leftrightarrow x_2^*)\le\pi^{10}_\pp(r)^2.
\]
\end{lemma}
\begin{proof}
The upper bound comes from the fact that a connection from $x_1^*$ to $x_2^*$ must leave both
half-boxes, and the lower bound from gluing two such one-arm events by two half-circuits and a
vertical crossing. The half-boxes $B_{x_1}:=B^+(x_1^*,r)$ and $B_{x_2}:=B^+(x_2^*,r)$, based on the
bottom and top sides of $\Omega^*_n$, are disjoint. We write $\phi:=\phi^0_{\Omega^*_n,\pp}$. Let
\[
E_{x_1}:=\{x_1^*\leftrightarrow\partial_{\rm out}B_{x_1}\text{ inside }B_{x_1}\},
\]
and define $E_{x_2}$ likewise.

\emph{Upper bound.} We have $\{x_1^*\leftrightarrow x_2^*\}\subseteq E_{x_1}\cap E_{x_2}$. Given
$\omega|_{B_{x_1}}$, the conditional law restricted to $B_{x_2}$ is free on the top side of
$\Omega^*_n$, and it is dominated by the measure that is, in addition, wired on $\partial_{\rm
out}B_{x_2}$. As $E_{x_2}$ is increasing, we obtain
$\phi(E_{x_2}\mid\omega|_{B_{x_1}})\le\pi^{10}(r)$. In the same way, $\phi(E_{x_1})\le\pi^{10}(r)$,
and the upper bound follows by conditioning on $\omega|_{B_{x_1}}$.

\emph{Lower bound.} As in Figure~\ref{fig:glue}, let $\mathcal O_{x_1}$ and $\mathcal O_{x_2}$ be
the events that there exist open half-circuits around $x_1^*$ and $x_2^*$ in $B_{x_1}\setminus
B^+(x_1^*,r/2)$ and $B_{x_2}\setminus B^+(x_2^*,r/2)$, respectively. Let $\bar x$ be the common
abscissa of $x_1^*,x_2^*$, and let $\mathcal V$ be the event that the column
\[
[\bar x-\lceil\frac n{32}\rceil,\bar x+\lceil\frac n{32}\rceil]\times[\frac12,n+\frac12]
\]
of $\Omega^*_n$ is crossed vertically. The aspect ratio of this column is at most $16$, and
$\lceil\frac n{32}\rceil\le r/2$. Hence, a vertical crossing starts in $B^+(x_1^*,r/2)$, ends in
$B^+(x_2^*,r/2)$ and meets both half-circuits. Consequently,
\[
E_{x_1}\cap\mathcal O_{x_1}\cap\mathcal V\cap\mathcal O_{x_2}\cap E_{x_2}\subseteq\{x_1^*\leftrightarrow x_2^*\}.
\]
All five events are increasing, so \tagref{FKG}{(FKG)} bounds the probability of the left-hand side
from below by the product of the five probabilities. Now
$\phi(E_{x_1}),\phi(E_{x_2})\ge\pi^{00}(r)$, because the restriction to a half-box dominates the
free measure on it; $\phi(\mathcal O_{x_1}),\phi(\mathcal O_{x_2})\ge c_{\rm hc}$ by
\eqref{eq:halfcircuit}; and $\phi(\mathcal V)\ge c_{\rm RSW}(\lambda,16)$ by RSW. This gives the
lower bound with $c_F:=c_{\rm hc}^2c_{\rm RSW}(\lambda,16)$.
\end{proof}

\begin{corollary}\label{cor:sandwich}
For $n\ge\max(64,7\lambda)$ and $\pp\in\mathcal W_\lambda(n)$,
\[
c_Fc_G^2\,A_{\lfloor n/8\rfloor}(\pp)^2\le\phi^0_{\Omega^*_n,\pp}(x_1^*\leftrightarrow x_2^*)\le c_G^{-2}\,A_{2\lfloor n/4\rfloor}(\pp)^2.
\]
\end{corollary}
\begin{proof}
This follows from Lemma~\ref{lem:factor} together with Lemma~\ref{lem:bc}(i) and~(ii) at $r=\lfloor
n/4\rfloor$, noting that $\lfloor\lfloor n/4\rfloor/2\rfloor=\lfloor n/8\rfloor$.
\end{proof}

We can now combine the gluing at a single parameter with the comparison across parameters in the
strip.

\begin{proof}[Proof of Theorem~\ref{thm:mumu}]
By duality, both boundary disorder correlations are connection probabilities in $\Omega^*_n$, at
$\pp(\beta)^*$ and at $\ppc$. By Corollary~\ref{cor:sandwich}, we bound each of them from above and
below by constant multiples of squares of connection probabilities $A_h$, and we compare these at
the two parameters by Proposition~\ref{prop:strip}.

We enlarge $n_1(M)$ so that $n_1\ge\max(64,7\lambda)$ and $\lfloor n_1/8\rfloor\ge h_0(\lambda)$,
where $\lambda=3M$. Let $n\ge n_1(M)$ and $\beta\in W_M$, and put $\pp:=\pp(\beta)^*$. By
Lemma~\ref{lem:dual} and $\ppc^*=\ppc$, we have
\[
\phi^1_\beta(\Dab)=\phi^0_{\Omega^*_n,\pp}(x_1^*\leftrightarrow x_2^*)\qquad\text{and}\qquad\phi^1_{\beta_c}(\Dab)=\phi^0_{\Omega^*_n,\ppc}(x_1^*\leftrightarrow x_2^*),
\]
with $\pp,\ppc\in\mathcal W_\lambda(n)$. Let $h_1:=\lfloor n/8\rfloor\le h_2:=2\lfloor n/4\rfloor$.
Both heights are at least $h_0(\lambda)$ and at most $n/2$, so that $h_i|\pp-\ppc|\le\lambda$.
Hence, Proposition~\ref{prop:strip} gives $c_S\le A_{h_i}(\pp)/A_{h_i}(\ppc)\le c_S^{-1}$ for
$i=1,2$, and Lemma~\ref{lem:Amono} gives $A_{h_2}\le A_{h_1}$. Applying Corollary~\ref{cor:sandwich}
at $\pp$ and at $\ppc$, we obtain
\begin{align*}
\phi^1_\beta(\Dab)&\le c_G^{-2}A_{h_2}(\pp)^2\le c_G^{-2}c_S^{-2}A_{h_2}(\ppc)^2\le c_G^{-2}c_S^{-2}A_{h_1}(\ppc)^2\le c_F^{-1}c_G^{-4}c_S^{-2}\,\phi^1_{\beta_c}(\Dab),\\
\phi^1_\beta(\Dab)&\ge c_Fc_G^2A_{h_1}(\pp)^2\ge c_Fc_G^2c_S^2A_{h_1}(\ppc)^2\ge c_Fc_G^2c_S^2A_{h_2}(\ppc)^2\ge c_Fc_G^4c_S^2\,\phi^1_{\beta_c}(\Dab).
\end{align*}
In each chain, the first and the last inequalities are Corollary~\ref{cor:sandwich}, at $\pp$ and at
$\ppc$, respectively; the second one is Proposition~\ref{prop:strip}; and the third one is the
monotonicity in $h$. Note that, in both chains, the change of scale runs in the direction allowed by
this monotonicity, so that no comparison $A_{h_1}\asymp A_{h_2}$ across scales is needed. Thus, the
bounds of the theorem hold for $n\ge n_1(M)$ with $(c_Fc_G^4c_S^2)^{-1}$ in place of $C(M)$;
together with the case $n_0(M)\le n<n_1(M)$ (\S\ref{ss:strip}), this proves the theorem.
\end{proof}

The estimates \tagref{E1}{(E1)}--\tagref{E4}{(E4)} are now proved: \tagref{E1}{(E1)} is
Proposition~\ref{prop:expmom}, \tagref{E2}{(E2)} is Proposition~\ref{prop:Bminus}, \tagref{E3}{(E3)}
is Theorem~\ref{thm:VP} and \tagref{E4}{(E4)} is Theorem~\ref{thm:mumu}. Hence, the discrete bounds
of \S\ref{ss:assembly} and the passage to the scaling limit in \S\ref{ss:limit} apply, and
Theorem~\ref{thm:main} is proved for $m>0$. For $m<0$, see the discussion at the end of
\S\ref{ss:strategy}.

\appendix
\section{Four-arm stability below the characteristic length: Proof of
Theorem~\ref{thm:R1}}\label{app:R1} In this appendix, we deduce Theorem~\ref{thm:R1} from results of
\cite{DCM}.

\emph{Inputs.} We use the four-arm events of \cite[\S2.2]{DCM}: $\AD(r,R)$ is the event that four
disjoint arms, of types open, dual-open, open, dual-open in counterclockwise order, join the inner
and outer boundaries of $\mathrm{Ann}(r,R):=\Lambda_R\setminus\Lambda_{r-1}$, and
$\piD(\pp;r,R):=\phi_\pp[\AD(r,R)]$. Restricting arms to a smaller annulus shows that
$\piD(\pp;r,R)$ is nondecreasing in $r$ and nonincreasing in $R$. Let $\rD$ be the smallest integer
with $\piD(\ppc;\rD,R)>0$ for all $R>\rD$, and put $\piD(R):=\piD(\ppc;\rD,R)$; this is the quantity
$\pi_4(R)$ of \cite[(1.17)]{DCM}. The constants of \cite{DCM} are uniform in
$\pp\in[\frac13,\frac34]$ (\S\ref{ss:inputs}), and $\asymp$ and $\lesssim$ below hold with such
constants.
\begin{enumerate}[label=\textup{(\roman*)},leftmargin=2.5em]
\item At $\ppc$, where $L(\ppc)=\infty$,
\[
\Delta_{\ppc}(1,r)\,\Delta_{\ppc}(r,R)\asymp\Delta_{\ppc}(1,R)
\]
for $1<r<R$ \cite[Thm.~1.6(ii), (1.15)]{DCM}, and $\Delta_{\ppc}(1,R)\asymp\Delta_{\ppc}(R)$ for
$R>2$ \cite[p.~45, after (4.13)]{DCM}, hence for $R>1$, since $\Delta_{\ppc}(1,2)$ and
$\Delta_{\ppc}(2)$ are positive.
\item There is $\varepsilon_1>0$ with
\[
\Delta_{\ppc}(r,R)\gtrsim(R/r)^{\varepsilon_1}\,\piD(R)/\piD(r)
\]
for $\rD<r<R$ \cite[Thm.~1.6(iv), (1.17)]{DCM}.
\item $\piD(r)\,\piD(\ppc;r,R)\asymp\piD(R)$ for $\rD<r<R$ \cite[Thm.~1.3]{CDH16}; see also
\cite[Prop.~2.3]{DCM}.
\item For $\pp\in[\frac13,\frac34]$ and $\rD<r<R\le L(\pp)$, we have, by \cite[Cor.~7.1]{DCM},
\[
\Bigl|\frac{d}{d\pp}\log\piD(\pp;r,R)\Bigr|\lesssim\frac{d}{d\pp}\phi_\pp[\mathcal C(\Lambda_R)],
\]
since the bound \cite[(7.2)]{DCM} has the right-hand side of \cite[(7.1)]{DCM} for the box
$\Lambda_R$.
\end{enumerate}

\begin{proof}[Proof of Theorem~\ref{thm:R1}]
\emph{Step 1: the bound at $\ppc$.} Let $\Delta:=\Delta_{\ppc}$ and $\rD<r<R$, so that $r\ge2$. By
(i) and Lemma~\ref{lem:Delta}(i),
\[
\Delta(r,R)\asymp\Delta(1,R)/\Delta(1,r)\asymp\Delta(R)/\Delta(r)\asymp r/R.
\]
Hence, by (iii) and (ii),
\[
\piD(\ppc;r,R)\asymp\frac{\piD(R)}{\piD(r)}\lesssim\Bigl(\frac rR\Bigr)^{\varepsilon_1}\Delta(r,R)\lesssim\Bigl(\frac rR\Bigr)^{1+\varepsilon_1} .
\]

\emph{Step 2: stability below the characteristic length.} We argue as in the proof of
\cite[Thm.~1.6(iii)]{DCM}. Let $\pp\in[\frac13,\frac34]$ and $\rD<r<R\le L(\pp)$. For every $u$
between $\ppc$ and $\pp$ and every $s$, $\phi_u[\mathcal C(\Lambda_s)]$ lies between its values at
$\ppc$ and at $\pp$, by \tagref{CBC}{(CBC)}, so that \eqref{eq:Ldef} and $L(\ppc)=\infty$ give
$L(u)\ge L(\pp)\ge R$. Hence, (iv) holds at every such $u$, and integrating it from $\ppc$ to $\pp$
gives
\[
\big|\log\piD(\pp;r,R)-\log\piD(\ppc;r,R)\big|\lesssim\big|\phi_\pp[\mathcal C(\Lambda_R)]-\phi_{\ppc}[\mathcal C(\Lambda_R)]\big|\le1 .
\]

\emph{Step 3: the bound below the characteristic length.} Restricting the four arms of $A_4(r,R)$ to
$\mathrm{Ann}(r+2,R-2)$, which absorbs the difference between the conventions for annuli, gives
$A_4(r,R)\subseteq\AD(r+2,R-2)$ for $R>r+4$. Now let $\pp\in[\frac13,\frac34]$ and $1\le r\le R\le
L(\pp)$, and put $r':=\max\{r,\rD\}+2$, so that $r'\le(\rD+5)\,r$. If $R\ge4r'$, then $R-2\ge R/2$
and
\[
\pi_4(\pp;r,R)\le\piD(\pp;r',R-2)\lesssim\piD(\ppc;r',R-2)\lesssim\Bigl(\frac{r'}{R-2}\Bigr)^{1+\varepsilon_1}\lesssim\Bigl(\frac rR\Bigr)^{1+\varepsilon_1},
\]
by the inclusion and the monotonicity in $r$, and then by Steps~2 and~1. If $R<4r'$, the same bound
holds trivially, since $\pi_4\le1$ and $r/R>1/(4\rD+20)$.

\emph{Step 4: the bound for $|\pp-\ppc|\le\lambda/N$.} Let $N\ge1$, $|\pp-\ppc|\le\lambda/N$ and
$1\le r\le R\le N$. If $\pp\notin[\frac13,\frac34]$, then $|\pp-\ppc|>\frac17$, so that
$r/R\ge1/N>1/(7\lambda)$, and the bound is trivial. Let $\pp\in[\frac13,\frac34]$. If $R\le L(\pp)$,
we apply Step~3. Otherwise, $\pp\ne\ppc$ and, by \eqref{eq:Lbound},
\[
L(\pp)\ge c_L/|\pp-\ppc|\ge c_LN/\lambda\ge c_LR/\lambda.
\]
If $r\le L(\pp)$, then \eqref{eq:A4mono} and Step~3 give
\[
\pi_4(\pp;r,R)\le\pi_4(\pp;r,L(\pp))\lesssim(r/L(\pp))^{1+\varepsilon_1}\le(\lambda/c_L)^{1+\varepsilon_1}(r/R)^{1+\varepsilon_1};
\]
if $r>L(\pp)$, then $r/R>c_L/\lambda$, and the bound is trivial.
\end{proof}

\section{The strip measure and the strip identity: Proof of
Theorem~\ref{thm:strip}}\label{app:strip} In this appendix, we sketch the construction of the strip
measure (Lemma~\ref{lem:stripmeasure}), and we prove Theorem~\ref{thm:strip} by combining the exact
relation between consecutive rows of the strip in \cite[\S3, Step~1]{BDC12} with the values of the
observable on the two boundary rows, as in the proof of \cite[Prop.~2.5]{DGP14}.

\subsection{The strip measure}\label{ss:stripmeasure} We identify $(x,y)\in\Z^2$ with
$x+iy\in\mathbb C$. For $k\in\Z$, let
\[
\mathcal C_k:=\{\{k+iy,k+1+iy\}:0\le y\le h-1\}
\]
be the \emph{column} of the $h$ horizontal edges crossing the line $\operatorname{Re}z=k+\frac12$
below the top row; we call it \emph{closed} if all its edges are closed. We put
$\varepsilon_\pp:=1-(1-\pp)^h\in(0,1)$. For an admissible boundary condition $\xi$
(\S\ref{ss:strip}), $u\overset{\xi}{\leftrightarrow}\text{top row}$ means that $u$ is joined to the
top row in $\omega^\xi$, that is, by a path made of open edges and of jumps between vertices of a
common block of $\xi$.

\begin{lemma}[the strip measure]\label{lem:stripmeasure}
Let $\pp\in(0,1)$.
\begin{enumerate}[label=\textup{(\roman*)},leftmargin=2em]
\item Let $x\in\Z$ and $1\le j<j'$, let $A$ be an event depending only on the edges of
$\Rr_{x-j,x+j}$ and let $\xi,\xi'$ be admissible boundary conditions on rectangles
$\Rr,\Rr'\supseteq\Rr_{x-j',x+j'}$. Then
$|\phi^\xi_{\Rr,\pp}(A)-\phi^{\xi'}_{\Rr',\pp}(A)|\le4\varepsilon_\pp^{j'-j}$. Consequently, for
every event $A$ depending on finitely many edges,
$\phi^{0,1}_{S_h,\pp}(A):=\lim\phi^\xi_{\Rr,\pp}(A)$ exists as $\Rr\uparrow S_h$ and does not depend
on the admissible $\xi$. These limits define a probability measure on $\{0,1\}^{E(S_h)}$, which is
invariant under horizontal translations and stochastically increasing in $\pp$, and
$|\phi^\xi_{\Rr,\pp}(A)-\phi^{0,1}_{S_h,\pp}(A)|\le4\varepsilon_\pp^{j'-j}$ in the above setting.
\item Let $x\in\Z$ and $j\ge1$, and let $\xi$ be admissible on a rectangle
$\Rr\supseteq\Rr_{x-2j,x+2j}$. Then
\[
|\phi^\xi_{\Rr,\pp}(x\overset{\xi}{\leftrightarrow}\textup{top row})-A_h(\pp)|\le8\varepsilon_\pp^j.
\]
\end{enumerate}
\end{lemma}

\begin{proof}[Sketch of the proof]
By \tagref{FE}{(FE)}, each column $\mathcal C_k$ is closed with conditional probability at least
$(1-\pp)^h$, given all other edges, and a closed column separates the two sides of a rectangle,
except through the wired top row. By \tagref{DMP}{(DMP)}, we may assume that
$\Rr=\Rr'=\Rr_{x-j',x+j'}$, since the boundary conditions induced there are admissible. Under the
monotone coupling of Lemma~\ref{lem:coupling} with the measure that wires both sides and the top
row, closed columns thus appear on both sides of $\Rr_{x-j,x+j}$ outside an event of probability at
most $2\varepsilon_\pp^{j'-j}$, and, by \tagref{DMP}{(DMP)}, the configurations agree between them.
This gives (i), and (ii) follows in the same way, from closed columns within distance $j$ of $x$.
\end{proof}

\subsection{The relation between consecutive rows}\label{ss:stripproof} Let $0<\pp<\ppc$, and let
\[
\alpha:=2\arctan(1-\pp)-\frac\pi4\in(0,\frac\pi4),
\]
so that
\[
e^{i\alpha}=\sqrt2\,(1+i(1-\pp))/(2-\pp+i\pp);
\]
this is the angle of \cite[(2.3)]{BDC12}. Let $F_\pp$ be the fermionic observable of
\cite{Smi10,BDC12} in the strip $S_h$, free on the bottom row and wired on the top row, defined as a
limit of finite-volume observables \cite[\S2]{BDC12}. In \cite{BDC12}, the strip is wired on the
bottom row and free on the top row; their results apply to $S_h$ after the reflection $z\mapsto\bar
z+ih$, which exchanges the two rows. We use the medial graph with its black and white diamonds, the
exploration path and its winding, as in \cite[\S2]{BDC12}. For $0\le y\le h$, let $e_y$ be the side
of the black diamond of the vertex $iy$ that borders the white diamond of the dual vertex
$\frac12+i(y-\frac12)$; after the reflection, it is the edge $e_{h-y}$ of \cite[\S3, Step~1]{BDC12}.
We put $f(y):=|F_\pp(e_y)|$ and use three properties of $f$.
\begin{enumerate}[label=\textup{(\roman*)},leftmargin=2.5em]
\item \emph{Bottom row.} $f(0)=\phi^{0,1}_{S_h,\pp}(0\leftrightarrow\text{top row})=A_h(\pp)$
\cite[Lemma~2.2]{BDC12}: the exploration path passes through $e_0$ exactly when $0$ is connected to
the top row, and its winding from $e_0$ is the same for all such configurations.
\item \emph{Top row.} $f(h)$ is the probability that the dual vertex $\frac12+i(h-\frac12)$, which
lies on the free row of the dual strip, is dual-connected to its wired row, formed by the dual
vertices below the bottom row. Indeed, as in (i), with the primal and dual clusters exchanged, the
exploration path passes through $e_h$ exactly when this occurs, and its winding from $e_h$ is then
deterministic. By planar duality, after a reflection and a translation, the dual strip is $S_h$,
with the same boundary conditions and the parameter $\pp^*$ (Figure~\ref{fig:strip}). Hence,
$f(h)=A_h(\pp^*)$ \cite[p.~15]{DGP14}.
\item \emph{Consecutive rows.} By \cite[\S3, Step~1]{BDC12}, for $1\le y\le h$,
\[
f(y-1)=\mathfrak r_\alpha\,f(y),\qquad\text{where}\qquad
\mathfrak r_\alpha:=\frac{\bigl[1+\cos(\frac\pi4+\alpha)\bigr]\cos(\frac\pi4+\alpha)}{\bigl[1+\cos(\frac\pi4-\alpha)\bigr]\cos(\frac\pi4-\alpha)} .
\]
This is the display after \cite[(3.4)]{BDC12}. It is derived there from the local relations
\cite[Lemmas~2.3 and~2.4]{BDC12} at two vertices of the medial graph between the two rows, which
have four adjacent medial edges even when $y=1$ or $y=h$, and from the invariance of the strip
measure under horizontal translations and under reflections in vertical lines
\cite[(3.3)--(3.4)]{BDC12}.
\end{enumerate}

\begin{proof}[Proof of Theorem~\ref{thm:strip}]
If $\pp=\ppc$, then $\pp^*=\pp$ and $\mathfrak r(\pp)=1$ (see the proof of
Proposition~\ref{prop:strip}), and there is nothing to prove. Let $0<\pp<\ppc$. By (iii),
$f(0)=\mathfrak r_\alpha^hf(h)$, and, by (i) and (ii), this reads $A_h(\pp)=\mathfrak
r_\alpha^hA_h(\pp^*)$. It remains to show that $\mathfrak r_\alpha=\mathfrak r(\pp)$. Put
$t:=1-\pp$. Since $\alpha=2\arctan t-\frac\pi4$, we have
\[\cos(\tfrac\pi4+\alpha)=\frac{1-t^2}{1+t^2}\qquad\text{and}\qquad\cos(\tfrac\pi4-\alpha)=\frac{2t}{1+t^2},
\]
so that $1+\cos(\frac\pi4+\alpha)=2/(1+t^2)$ and $1+\cos(\frac\pi4-\alpha)=(1+t)^2/(1+t^2)$. Hence,
\[
\mathfrak r_\alpha=\frac{2(1-t^2)}{2t(1+t)^2}=\frac{1-t}{t(1+t)}=\frac{\pp}{(1-\pp)(2-\pp)}=\mathfrak r(\pp).\qedhere
\]
\end{proof}
\section*{Acknowledgments} The author thanks Dmitry Chelkak for bringing the problem of scaling
limits of near-critical Ising interfaces in the thermal direction to his attention and for
illuminating discussions.

\emph{AI use statement.} The author used Claude Opus 5.5 to improve the exposition of the paper. The author takes full responsibility of this paper.

\renewcommand{\bibliofont}{\small}

\end{document}